\documentclass[10pt]{article}
\usepackage{fullpage,epsfig,graphics,cancel,slashed}
\usepackage{amsthm,amssymb,amsmath,amsfonts,amscd,amsbsy}%
\usepackage{lmodern}%
\usepackage{slantsc}%
\usepackage{psfrag,hyperref,color}
\usepackage[]{youngtab}
\usepackage{graphicx}
\usepackage{wasysym}
\usepackage{mathdots}
\usepackage{mathrsfs}

\usepackage{tikz}
\usetikzlibrary{arrows,chains,matrix,positioning,scopes,decorations.pathmorphing,bending,calc}
\usepackage{tikz-cd}
\tikzset{snake it/.style={decorate, decoration=snake}}

\makeatletter
\tikzset{join/.code=\tikzset{after node path={%
\ifx\tikzchainprevious\pgfutil@empty\else(\tikzchainprevious)%
edge[every join]#1(\tikzchaincurrent)\fi}}}
\makeatother

\tikzset{>=stealth',every on chain/.append style={join},
         every join/.style={->}}
\tikzset{
    >=stealth',
    punkt/.style={
           rectangle,
           rounded corners,
           draw=black, very thick,
           text width=6.5em,
           minimum height=2em,
           text centered},
    pil/.style={
           ->,
           thick,
           shorten <=2pt,
           shorten >=2pt,}
}

\def\centerarc[#1](#2)(#3:#4:#5){\draw[#1] ($(#2)+({#5*cos(#3)},{#5*sin(#3)})$) arc (#3:#4:#5)}

\newcommand{\BB}{\mathbb}

\newcommand{\FR}{\mathfrak}

\newcommand{\bea}{\begin{eqnarray}}
\newcommand{\eea}{\end{eqnarray}}
\newcommand{\nn}{\nonumber}

\newcommand{\sbullet}{\vcenter{\hbox{\tiny$\bullet$}}}
\newcommand{\udl}{\underline}
\newcommand{\bra}{\langle}
\newcommand{\ket}{\rangle}

\newcommand{\To}{\Rightarrow}

\newcommand{\ext}{\operatorname{Ext}}

\newcommand{\opn}{\operatorname}
\newcommand{\CUP}{\vcenter{\hbox{\large$\cup$}}}
\newcommand{\CAP}{\vcenter{\hbox{\large$\cap$}}}
\newcommand{\OC}{\begin{tikzpicture}[scale=0.5]
\draw[black] (0,0) -- (0.5,0.5);
\draw[black] (0.5,0) -- (0.34,0.16);
\draw[black] (0.16,0.34) -- (0,0.5);
\end{tikzpicture}}
\newcommand{\UC}{\begin{tikzpicture}[scale=0.5]
\draw[black] (0.5,0) -- (0,0.5);
\draw[black] (0,0) -- (0.16,0.16);
\draw[black] (0.34,0.34) -- (0.5,0.5);
\end{tikzpicture}}
\def\ot{\leftarrow}

\newcommand\hcancel[2][black]{\setbox0=\hbox{$#2$}%
\rlap{\raisebox{.45\ht0}{\textcolor{#1}{\rule{\wd0}{1pt}}}}#2}

\def\Gc{\Gamma}

\def\Gd{\Delta}
\def\ep{\epsilon}

\def\gs{\sigma}
\def\Gs{\Sigma}

\def\Gl{\Lambda}

\def\go{\omega}

\DeclareMathAlphabet{\mathpzc}{OT1}{pzc}{m}{it}

\newtheorem{theorem}{Theorem}[section]
\newtheorem{lemma}[theorem]{Lemma}
\newtheorem{proposition}[theorem]{Proposition}

\theoremstyle{definition}
\newtheorem{example}[theorem]{Example}
\newtheorem{remark}[theorem]{Remark}
\newtheorem{corollary}[theorem]{Corollary}
\newtheorem{definition}[theorem]{Definition}

\numberwithin{equation}{section}
\begin{document}
\begin{flushright} \small
UUITP-26/26
 \end{flushright}
\smallskip
\begin{center} \Large
{\bf Tangle Invariants Valued in Derived Categories of Sheaves and the Rouquier Complex}
 \\[12mm] \normalsize
{\bf Jian Qiu${}^{a,b}$ and Alex Takeda${}^{b}$} \\[8mm]
 {\small\it
    ${}^a$Department of Physics and Astronomy, Uppsala University,\\
    Uppsala, Sweden\\  
        \vspace{.3cm}
      ${}^b$ Instituto de Matem\'atica Pura e Aplicada,\\
    Rio de Janeiro, Brazil.
      \vspace{.3cm}
   }
\end{center}
\vspace{7mm}
\begin{abstract}
\noindent
  We develop effective calculation tools for constructing categorified invariants of tangles and knots, following the work of Cautis and Kamnitzer. These invariants, originally constructed by means of Fourier--Mukai transform functors between derived categories of coherent sheaves on certain iterated flag varieties, can be given an alternative algebraic description by using tilting bundles. We use this method and provide several examples of its application, in the hope of showing its utility to explicit pen-and-paper computations. This tilting-bundle method also allows us to lift the construction from derived to homotopy categories of complexes of modules, reducing the calculation of tangle invariants to a type of categorified matrix algebra, where the matrices in question have maps between their elements and calculations get turbocharged with chain homotopy.
\end{abstract}

\eject
\normalsize

\tableofcontents

\section{Introduction}
In Chern--Simons theory, an important correlation function to compute is that of the Wilson loops. These are links embedded in a 3-manifold 
with its components coloured with different representations of the gauge group. The correlation function is an invariant of the link, e.g. with $SU(2)$ gauge group and colouring with the fundamental representation, one gets the Jones polynomial \cite{witten1989}. Rozansky--Witten theory, a fermionic cousin to Chern--Simons theory, also naturally admits Wilson loops whose components are coloured with objects of the sheaf category of some specific type of variety $X$. This way, one constructs link invariants valued in the derived category $D^b(X)$ \cite{RobertsWillerton}, see also \cite{Kapustin:2008sc} \cite{Kapustin:2009uw}, as well as \cite{Khovanov2004MatrixFA} for a link invariant valued in matrix factorisation categories.

In their paper \cite{Cautis2007KnotHV}, Cautis and Kamnitzer also constructed a tangle invariant valued in some derived category of sheaves. This invariant turns out to be a categorification of the Jones polynomial, equivalent to the Khovanov homology, though the way this arises is from the category of sheaves and much more geometrical. Schematically, to a tangle $K$ with $m$ incoming and $n$ outgoing strands 
\bea 
\begin{tikzpicture}[scale=.8]
\draw [blue,-] (0,0) -- (2,0) -- (2,1) -- (0,1) -- (0,0);
\draw [blue,->] (0.3,-.5) -- (0.3,0);
\draw [blue,->] (0.5,-.5) -- (0.5,0);
\draw [blue,->] (1.7,-.5) -- (1.7,0);
\draw [blue,->] (1.5,-.5) -- (1.5,0);
\draw [blue,->] (0.4,1) -- (0.4,1.5);
\draw [blue,->] (1.6,1) -- (1.6,1.5);
\node at (1,-0.4) {$\cdots$};
\node at (1,1.4) {$\cdots$};
\node at (1,-.7) {$\underbrace{\phantom{XXXX}}_{m}$};
\node at (1,1.7) {$\overbrace{\phantom{XXX}}^{n}$};
\node at (1,0.5) {$K$};
\end{tikzpicture}\nn\eea 
they assign a functor between derived categories
\bea K_m^n:\, D^b(Y_m)\to D^b(Y_n), \nn\eea
where $\{Y_n\}_{n \ge 1}$ is a certain family of smooth projective varieties. 
These maps can be realised as Fourier--Mukai transforms with kernels 
\bea \Phi(K)\in D^b(Y_m\times Y_n),\nn\eea
the idea being that, if two tangles $K,K'$ are isotopic, the corresponding kernels $\Phi(K),\Phi(K')$ are quasi-isomorphic. 

The long-term purpose of our programme is to better understand the constructions of \cite{Cautis2007KnotHV} and study the deformation theory of their categorification, using the method developed by the first-named author in \cite{Def_qnt}, which relies on the techniques of \cite{BarmeierWang}. This method allows for a systematic study of the deformation theory of derived categories of coherent sheaves, starting from a choice of exceptional collection, and is conjectured to give a concrete realization of the brane quantisation procedure of \cite{Gukov:2008ve}. We also believe that this perspective would shed new light on the relation between these categorified invariants and the link invariants coming from Rozansky--Witten theory. Though having similar buzzwords, there seem to be some big differences between the two constructions: the latter colours each strand with an arbitrary object in a fixed sheaf category, while the former only uses the link diagram to encode morphisms between different categories.

In this paper, we take the first step in the development of this programme, which consists in translating the constructions and results of Cautis and Kamnitzer into the language that one would need for the brane quantisation procedure, i.e. using tilting bundles and minimal resolutions. What we lack in conceptual breakthrough we make up for in terms of concreteness: our calculations in the derived category are mostly done by passing to an equivalent module category thanks to tilting. For the varieties $\{Y_n\}$, we find an exceptional collection consisting of $2^n$ line bundles $\{L_I\}$. Tilting with respect to this exceptional collection allows us to pass to an equivalent module category
\bea D^b(Y_n)\simeq D^b(\opn{mod}_{\Gl_n}),\label{tilting_equ}\eea
where $\Gl_n$ is a quiver algebra attached to $Y_n$. 
The matter of understanding a map $D^b(Y_m)\to D^b(Y_n)$ can be settled on the exceptional objects and maps therebetween. 
As for the module category, one can exploit the equivalence 
\bea D^b(\opn{mod}_{\Gl_n})\simeq K^b(\opn{proj}_{\Gl_n}), \label{proj/htpy}\eea
where it is important that $\Gl_n$ have finite projective dimension.
From this view, kernels of FM transforms $D^b(\opn{mod}_{\Gl_m})\to D^b(\opn{mod}_{\Gl_n})$ are just (a complex of projective) $\Gl_m$-$\Gl_n$ bi-modules. An FM transform is then just the ordinary tensor product with its kernel.
This is essentially the idea of Rouquier complex, to be discussed in sec.\ref{sec_LthHC}. We give an algorithm for computing these kernels starting from the exceptional collections, and also a rule for computing the tensor of two kernels as a matrix multiplication, with the possibility of applying a chain homotopy at each step.

Below is a summary of the results and organisation of the paper. 
In sec.\ref{sec_KIViDCoS}, we recall the basic properties of the central geometric object, the variety $Y_n$, while leaving the more technical aspects to appendix sec.\ref{Sec_DCoY_n}.
We move then to the crossings and cup/caps in sec.\ref{sec_CCC}. 
The Kauffman distinguished triangle in $D^b(Y_n)$ (prop.\ref{prop_Kauffman}), e.g.
\bea 
\begin{tikzpicture}[scale=1]
\draw [blue,-] (2,0) -- (3,1);
\draw [white,-,line width=2mm] (3,0) -- (2,1);
\draw [blue,-] (3,0) -- (2,1);
\node at (1.3,0.5) {\scriptsize$\longrightarrow~~t^{-1/4}$};
\centerarc[blue](-0.5,0.5)(65:-65:0.5);
\centerarc[blue](0.8,0.5)(115:245:0.5);
\node at (3.8,0.5) {\scriptsize$\longrightarrow~~t^{-1/2}$};
\centerarc[blue](5,1.2)(-35:-145:0.5);
\centerarc[blue](5,-0.2)(35:145:0.5);
\node at (5.9,0.5) {\scriptsize{$\stackrel{+1}{\longrightarrow}$}};
\end{tikzpicture}\nn\eea
is a very handy tool for small calculations.
We also comment on Rouquier's construction in sec.\ref{sec_RC} by phrasing the twist functor in terms of adjunctions.
We claim no originality in these sections compared to what is already in \cite{Cautis2007KnotHV}, though some proofs are made more direct. Hopefully, the reader who is interested in practical calculations will benefit from our presentation. 

In sec.\ref{sec_CS} we give quite a few examples of computations. We make use of explicit formulae for various functors in the sheaf category, proved in the appendix. Some calculations are done with the exceptional collection of $D^b(Y_n)$, consisting of $2^n$ line bundles. This collection and the tilting is shown in appendix sec.\ref{sec_ECoDY}.

As an aside, if one takes a link and applies the above Kauffman DT to each crossing, then the link is expressed as repeated cones between different ways of resolving a crossing. This repeated cone can be organised into a cube, whose vertices correspond to disjoint unions of circles, resulting from smoothing all crossings in different ways. Computing the cohomology of this cube brings us to the Khovanov homology of the link. We illustrate this in sec.\ref{sec_RtKh}.
Incidentally, the same algorithm also allows us to prove certain isotopy (Reidemeister) moves rather quickly.

If we descend to the Grothendieck group, the Kauffman DT above de-categorifies to the usual Kauffman bracket relation e.g.
\bea 
\begin{tikzpicture}[scale=1]
\draw [blue,-] (0,0) -- (1,1);
\draw [white,-,line width=2mm] (1,0) -- (0,1);
\draw [blue,-] (1,0) -- (0,1);

\node at (1.1,0.5) {$=$};

\node at (1.6,0.55) {\scriptsize$t^{1/4}$};
\centerarc[blue](1.5,0.5)(65:-65:0.5);
\centerarc[blue](2.8,0.5)(115:245:0.5);

\node at (2.8,0.5) {$+$};

\node at (3.4,0.5) {\scriptsize$t^{-1/4}$};
\centerarc[blue](4,1.2)(-35:-145:0.5);
\centerarc[blue](4,-0.2)(35:145:0.5);
\node at (4.8,0.2) {$,$};

\draw[blue] (7,0.5) ellipse (0.4 and 0.55);
\node at (8.6,0.5) {\scriptsize$=-t^{-1/2}-t^{1/2}$};
\end{tikzpicture}\nn\eea
The resulting link invariant is the Jones polynomial, as expected. We point out that 
\cite{Cautis2007KnotHV} uses a nice base for the Grothendieck group which is unfortunately not exceptional, as there are higher Ext's in between. In this base, one can give short formulae for all caps, cups and crossings for all $n$. In fact, they are just matrices:\footnote{Though we have two equivariant parameters $s,t$, the $s$ never seems to feature in any final result, but it does help organise the calculation in the derived category.} 
\bea \OC,\UC  \in \opn{Mat}_{2^n\times2^n}(R),~~
\CAP \in \opn{Mat}_{2^n\times2^{n-2}}(R),~~\CUP \in \opn{Mat}_{2^n\times2^{n+2}}(R),~~~R=\BB{Z}[t^{\pm1/2},s^{\pm1}].\nn\eea
The same formulae are much more complicated with the exceptional collection. One can imagine combining the Cautis--Kamnitzer collection, for which the traditional tilting would not fly, 
with tilting in the DG or $A_{\infty}$ sense so as to have the best of both worlds. We leave this for future exploration. 

In sec.\ref{sec_LthHC}, we lift the previous construction to the homotopy category, in line with Rouquier's idea mentioned above \cite{rouquier2006categorification}. There are multiple benefits to doing this. First, the maps in the Kauffman DT above exist by the axioms of derived category or adjunction and are not always straightforward to describe. 
But in $K^b(\opn{proj}_{\Gl_n})$, the only manoeuvre possible is chain homotopy, thus we can always keep track of all maps by following these chain homotopies. Second, in $K^b(\opn{proj}_{\Gl_{\sbullet}})$, we can functorially assign to a braid $K$ the kernel 
\bea \Phi(K_m^n)\in K^b({}_{\Gl_m}\hspace{-1mm}\opn{proj}_{\Gl_n})\nn\eea
that is a complex of projective $\Gl_m$-$\Gl_n$ bi-modules effecting the FM transform. 
Once given the formulae for the crossings/cap/cup, one can tensor them together and get the kernel for any tangle valued in $K^b({}_{\Gl_m}\hspace{-1mm}\opn{proj}_{\Gl_n})$. The homotopy class (which is much easier to treat) of this kernel then serves as the tangle invariant.
For this part, we are only able to go up to $n=3$, because the calculation, though conceptually clear,  rather quickly snowballs out of hand. We plan to revisit this problem with some computer coding. In this work, we only give the formulae for various kernels for $m,n\leq 3$. We then demonstrate how certain isotopy moves result from simple matrix multiplications and chain homotopies.
\vspace{0.5cm}

\noindent {\bf Acknowledgements}: Though we have never met Sabin Cautis and Joel Kamnitzer, we feel deeply indebted to them and are amazed by the depths on their paper \cite{Cautis2007KnotHV} and the sequels. A.T. would like to thank Uppsala University and IMPA for the working conditions provided during this project. This project was partly funded by the Knut and Alice Wallenberg Foundation.
\vspace{0.5cm}

\noindent\textbf{AI use statement} In no part of the conception, materialisation or the write-up of the paper did we use AI. All calculations were done by hand. Though once the computation was done, we used AI assistance to write computer codes to check some of the messier tangles such as the reef or the bowline knot.

\section{Knot Invariants Valued in Derived Categories of Sheaves}\label{sec_KIViDCoS}
\subsection{Basics of \texorpdfstring{$Y_n$}{Yn}}\label{sec_BoY}
The family of varieties $\{Y_n\}$ with the requisite properties for making tangle invariants was constructed in \cite{Cautis2007KnotHV}. We gather the most essential information about $Y_n$ here, while leaving the more technical aspects to appendix sec. \ref{Sec_DCoY_n}.

Denote the standard basis of $V=\BB{C}^2$ as $\{e,f\}$. We will need $N$ copies of $V$, and we name the basis of the $n^{th}$ copy as $e_n,f_n$. Let $z$ be the up-shift operation that sends $V_n$ to $V_{n-1}$
\bea z \colon \begin{cases}
  V_n\stackrel{id}{\to}V_{n-1}, & n\geq2 \\
  V_1\to 0 &  \\
\end{cases} \eea
We also denote by $\gs$ the reverse shift
\bea \gs(e_i)=e_{i+1},~~\gs(f_i)=f_{i+1}.\nn\eea
Note that the preimage operation $z^{-1}$ is related to $\gs$ as $z^{-1}(\ell)=\gs(\ell)+V_1$. 

We define two $\BB{C}^*$-actions, with equivariant parameters $t,s$, for which the bases vectors have weights 
\bea e_n\colon\ t^{-n},~~~f_n \colon\ t^{-n}s^{-1}.\nn\eea
We will discuss this in more detail in appendix sec.\ref{sec_equivariance}. What is important for now is that $z$ be equivariant of  weight $t$, so that all constructions below are equivariant.

\smallskip

\noindent\textbf{Convention}: Before going further, we want to clarify our convention in dealing with equivariant/graded categories. The equivariance mentioned above gives rise to two gradings on all the vector spaces and categories involved. We use the convention that in a graded category all morphisms must have grading zero. But one can shift the gradings of objects: for example, if ${\cal O}_Y$ is the structure sheaf of a variety $Y$, then $s^2t{\cal O}_Y$ means shifting its grading by $s^2t$. For a graded vector space, we will write e.g. $1+s$ to mean a 2-dim vector space with two copies of $\BB{C}$, one of them having its degree shifted by $s$. With this convention, rather than writing $\BB{C}$, we will write 1 instead. Note that the cohomological grading shift is still written in the standard manner, so that $[1]$ shifts complexes down one degree.

\begin{definition}\label{def_Yn}
The variety $Y_n$ is a flag-variety in $\BB{C}^{2N}$ obeying the following conditions
\bea Y_n=\{{\cal L}_1 \subset {\cal L}_2 \cdots \subset {\cal L}_n\ |\ \dim {\cal L}_i=i,~z{\cal L}_i\subset {\cal L}_{i-1}\},\label{eq_def_Yn}\eea
where we take ${\cal L}_0=\{0\}$.
\end{definition}
For example, $Y_1=\BB{P}^1$ and $Y_2=\BB{F}_2$, the second Hirzebruch surface, explained in ex.\ref{ex_Y2} (it might be helpful to peep at these examples to get a feeling for the geometry).
We find the following `mixing-board' picture of $Y_n$ rather useful.
 \bea
\begin{tikzpicture}
  \matrix (m) [matrix of math nodes, row sep=2em, column sep=3em]
    { e_1 & \gs(\ell_1) & \gs(\ell_2) & \gs(\ell_3) & \gs(\ell_4) \\
      f_1 & {\cal L}_1^{\perp} & {\cal L}_2^{\perp} & {\cal L}_3^{\perp} & {\cal L}_4^{\perp} \\ };
  \path[-, font=\scriptsize]
  (m-1-1) edge (m-2-1)
  (m-1-2) edge (m-2-2)
  (m-1-3) edge (m-2-3)
  (m-1-4) edge (m-2-4)
  (m-1-5) edge (m-2-5);
  \node at (-11.25em,0) {$\bullet~\ell_1$};
  \node at (-5.98em,.5em) {$\bullet~\ell_2$};
  \node at (-0.08em,-0.9em) {$\bullet~\ell_3$};
  \node at (5.85em,0.2) {$\bullet~\ell_4$};
  \node at (11.79em,-0.2) {$\bullet~\ell_5$};
\end{tikzpicture}\label{mix_board}\eea
To read this picture, 
denote the generator of the line ${\cal L}_{n+1}/{\cal L}_n$ as $\ell_{n+1}$, it is the linear combination of $\gs(\ell_n)$ and ${\cal L}_n^{\perp}$ where 
\bea {\cal L}_n^{\perp}=z^{-1}({\cal L}_{n-1})/{\cal L}_n.\label{def_perp}\eea
This is a $\BB{P}^1$ worth of choice, which is depicted as a segment between $\gs(\ell_n),\,{\cal L}_n^{\perp}$. At one end of the segment, say $\gs(\ell_n)$ is sent isomorphically under $z$ to $\ell_n$, while at the other end ${\cal L}_n^{\perp}$ is sent to 0 mod ${\cal L}_{n-1}$.
We have that $Y_n$ is an iterated $\BB{P}^1$-fibration
\bea
\begin{tikzpicture}
  \matrix (m) [matrix of math nodes, row sep=0.6em, column sep=2em]
    { Y_n & \BB{P}_n^1 \\
      Y_{n-1} & \BB{P}_{n-1}^1 \\
      \vdots & \vdots \\
      Y_2 & \BB{P}_2^1 \\
      Y_1 & \BB{P}^1 \\ };
  \path[->, font=\scriptsize]
  (m-1-2) edge (m-1-1)
  (m-2-2) edge (m-2-1)
  (m-4-2) edge (m-4-1)
  (m-1-1) edge node[right] {$\pi_n$} (m-2-1)
  (m-2-1) edge node[right] {$\pi_{n-1}$} (m-3-1)
  (m-4-1) edge node[right] {$\pi_2$} (m-5-1);
  \draw[double](-0.5,-1.7)--(.5,-1.7);
 \end{tikzpicture}\label{tower}\eea
where we denote the successive projections as $\pi_k$. 

Each $\BB{P}^1$-fibre has a special point where the $k^{th}$ `dial' in the picture \ref{mix_board} is at the bottom, then
$\ell_k\in {\cal L}_{k-1}^{\perp}=z^{-1}({\cal L}_{k-2})/{\cal L}_{k-1}$, hence $z(\ell_k)\subset {\cal L}_{k-2}$.
This gives rise to a family of divisors $\{D_k\}$
\bea D_k=\{z{\cal L}_k\subset {\cal L}_{k-2}\},~~k\geq2;~~~~D_1:=\opn{span}\bra f_1\ket\label{divisor_D}\eea
that plays a crucial role in the construction of \cite{Cautis2007KnotHV}. Note that if we want to emphasise that $D_k$ be in $Y_n$, then we write $D_k^n$ as we do in part 2 of lem.\ref{lem_trivial_lb}. 
\begin{lemma}\label{lem_fibration}
The locus $D_k$ is also a $\BB{P}^1$-fibration over $Y_{n-2}$
\bea \begin{tikzpicture}
  \matrix (m) [matrix of math nodes, row sep=2em, column sep=2em]
    { D_k  & Y_n  \\
      Y_{n-2} & \\ };
  \path[>->,font =\scriptsize]
  (m-1-1) edge node[above]{$\iota_k$} (m-1-2);
  \path[->,font =\scriptsize]
  (m-1-1) edge node[right] {$q_k$} (m-2-1);
\end{tikzpicture}\label{Dk_fibre}\eea
with the projection $q_k:\,D_k\to Y_{n-2}$ defined as
\bea D_k\ni ({\cal L}_1\subset \cdots\subset {\cal L}_{k-1}\subset {\cal L}_k\subset {\cal L}_{k+1}\cdots)\to ({\cal L}_1\subset \cdots\subset {\cal L}_{k-2}\subset z{\cal L}_{k+1}\subset z{\cal L}_{k+2}\subset\cdots).\label{q_k}\eea
\end{lemma}
\begin{proof}
We need to prove first that the rhs does lie in $Y_{n-2}$ i.e. satisfies the incidence relation \eqref{eq_def_Yn}. %

At $D_k$, we have $z{\cal L}_k\subset {\cal L}_{k-2}$, but $\dim z{\cal L}_k\geq k-2=\dim {\cal L}_{k-2}$, so $z{\cal L}_k={\cal L}_{k-2}$ and ${\cal L}_k$ contains the dim 2 kernel of $z$ (so will ${\cal L}_i$, for all $i\geq k$, of course).
Thus for all $z{\cal L}_i$, $i\geq k$, we have $\dim z{\cal L}_i=i-2$, so the dimensionality is as required. Further ${\cal L}_{k-2}=z{\cal L}_k\subset z{\cal L}_{k+1}\subset z{\cal L}_{k+2} \cdots$ is also clear.
Finally $zz{\cal L}_{k+1}\subset z{\cal L}_k={\cal L}_{k-2}$, proving the incidence relation.

As ${\cal L}_k$ is determined entirely as $z^{-1}{\cal L}_{k-2}$, the projection map forgets the line ${\cal L}_{k-1}/{\cal L}_{k-2}$. We recall this latter has a 2-dimensional space worth of choices, so the fibre of $q_k$ is $\BB{P}^1$ and $q$ is flat. The key here is the single relation ${\cal L}_k=z^{-1}{\cal L}_{k-2}$, so that when one constructs ${\cal L}_{k+1}$ mod ${\cal L}_k$, one needs pay no regard to ${\cal L}_k,{\cal L}_{k-1}$. The projection $q_k$ essentially forgets these two.
\end{proof}
We summarise the main facts that went into the proof above, as it is useful to keep these in mind for the rest of the paper. First $z{\cal L}_k$ has dimension $k-1$ (in which case it is ${\cal L}_{k-1}$) or $k-2$ (when ${\cal L}_k$ contains the kernel of $z$). Further $z{\cal L}_k$ is included in ${\cal L}_{k-1}$ always, and $z{\cal L}_k+{\cal L}_{k-2}={\cal L}_{k-1}$ generically, except when $z{\cal L}_k\subset {\cal L}_{k-2}$ (which is the locus $D_k$), but it will not be a subspace of ${\cal L}_{k-3}$, due to dimension reasons. We gather a few more results about $Y_n$ and the divisors $D_k$ next.
\begin{lemma}\label{lem_trivial_lb}(c.f. lemma 4.3 \cite{Cautis2007KnotHV})
\begin{enumerate}
  \item The tautological line bundle ${\cal E}_k={\cal L}_k/{\cal L}_{k-1}$ over $Y_n$ is isomorphic to $t^{-k}{\cal O}(-D_1-\cdots-D_k)$.
  \footnote{Recall that equivariant varieties are defined using graded rings and homogeneous ideals. In this setting the divisors must be invariant under the two $\BB{C}^*$ actions. As for the grading shift $t^{-k}$, see the clarification of convention in the beginning of this section. So e.g. $st^{-1}{\cal O}_{Y_n}$ means that the grading of the ring of homogeneous rational functions is shifted by $st^{-1}$.}
  \item The projections $q_{k+1},q_{k+2}$ are isomorphisms when restricted to the locus $D_{k+1}\cap D_{k+2}$. 
  But for $l>k+1$, we have the pullback square
  \bea \begin{tikzpicture}
  \matrix (m) [matrix of math nodes, row sep=2em, column sep=2em]
    { D^{n+2}_{k+1}\cap D^{n+2}_{l+1}  & D^n_{l-1}  \\
      D^n_{k+1} & Y_{n-2} \\ };
  \path[->,font =\scriptsize]
  (m-1-1) edge node[above]{$q_{k+1}$} (m-1-2)
  (m-2-1) edge node[above]{$q_{k+1}$} (m-2-2);
  \path[->,font =\scriptsize]
  (m-1-2) edge node[right] {$q_{l-1}$} (m-2-2)
  (m-1-1) edge node[right] {$q_{l+1}$} (m-2-1);
\end{tikzpicture}\label{pullback_square}\eea
where the superscript $m$ in $D^m_i$ means $D_i$ is a divisor in $Y_m$.

  \item We have the following isomorphism of line bundles 
\bea {\cal O}(D_{k+1}+D_{k+2})|_{D_{k+1}\cap D_{k+2}}=t^{-1}{\cal O}_{D_{k+1}\cap D_{k+2}}\cong t^{-1}{\cal O}_{Y_{n-2}}.\nn\eea
  \item Let $\iota:\,D_{k+1}\hookrightarrow Y_n$ be the embedding, then $\iota^*{\cal O}(-D_{k+1})\cong t^ks^{-1}\iota^*{\cal O}(2D_1+\cdots+2D_k)$, which contains the previous one as a special case.
\end{enumerate}
\end{lemma}
As the proof is on the technical side, we give it in the appendix sec.\ref{sec_GoY}. From the third result of the lemma, we have a neat formula for the equivariant characters of the line bundles ${\cal O}(\sum_ia_iD_i)$, see prop.\ref{prop_sec_L}. The line bundles ${\cal O}(-D_1-\cdots-D_k)$ that appear in the lemma above will often be used in this paper, so we introduce the notation
\[ L_{[k]}={\cal O}(D_1+\cdots+D_{k}) \quad \text{and} \quad L_{-[k]}={\cal O}(-D_1-\cdots-D_k), \]
and for any sheaf $F$, we will also use the shorthand $F_{[k]}=F\otimes L_{[k]}$. 

\subsection{Crossings, Caps and Cups}\label{sec_CCC}
To construct the kernel $\Phi(K)$ associated to a tangle $K$, we need only the cases where $K$ is an over- or undercrossing, a cap and cup.
We denote with the hieroglyphs ${}_i\OC$ and ${}_i\UC$ the following overcrossing and undercrossing between $i^{th}$ and $(i+1)^{th}$ strand
\bea 
\begin{tikzpicture}[scale=1]

\draw[blue] plot[smooth, tension=.7] coordinates {(-0.2,-0.3) (-0.1,-0.1) (0.1,0.1) (0.2,0.3)};
\draw[blue] plot[smooth, tension=.7] coordinates {(0.2,-0.3) (0.1,-0.1)};
\draw[blue] plot[smooth, tension=.7] coordinates {(-0.1,0.1) (-0.2,0.3)};
\node at (-0.3,-0.5) {\scriptsize$i$};
\node at (0.4,-0.5) {\scriptsize$i+1$};

\end{tikzpicture}~~~~~
\begin{tikzpicture}[scale=1]
\draw[blue] plot[smooth, tension=.7] coordinates {(0.2,-0.3) (0.1,-0.1) (-0.1,0.1) (-0.2,0.3)};
\draw[blue] plot[smooth, tension=.7] coordinates {(-0.2,-0.3) (-0.1,-0.1)};
\draw[blue] plot[smooth, tension=.7] coordinates {(0.1,0.1) (0.2,0.3)};
\node at (-0.3,-0.5) {\scriptsize$i$};
\node at (0.4,-0.5) {\scriptsize$i+1$};

\end{tikzpicture}
\nn\eea 
We follow \cite{Cautis2007KnotHV} in constructing these two kernels. Consider the subvariety $Z^i_n\subset Y_n\times Y_n$ defined as
\bea Z_n^i=\{({\cal L}_1\subset\cdots\subset {\cal L}_n)\times ({\cal L}_1'\subset\cdots\subset {\cal L}_n')\,|\,{\cal L}_j={\cal L}_j',~~~j\neq i\}.\nn\eea
It has two branches, one branch is simply the diagonal when ${\cal L}_j={\cal L}_j'$ for all $j$, the other is
\bea V^i_n=D_{i+1}\times_{q_{i+1}} D_{i+1}\nn\eea
where $q_k$ is projection map defined in \eqref{Dk_fibre} and \eqref{q_k}. Indeed,
we mentioned above that generically $z{\cal L}_{i+1}+{\cal L}_{i-1}={\cal L}_i$. Thus ${\cal L}_j={\cal L}_j'$ for $j\neq i$ also implies ${\cal L}_i={\cal L}_i'$ generically, which is the diagonal. This reasoning is invalid precisely on $D_{i+1}$ and $D_{i+1}'$ where $z{\cal L}_{i+1}={\cal L}_{i-1}$, giving the other branch.

Following sec.4.2.2 of \cite{Cautis2007KnotHV}, we take the structure sheaf of $Z_n^k$ as the kernel for ${}_k\OC$, more precisely
\bea \Phi({}_k\OC)=t^{-1/4}{\cal O}_{Z_n^k},~~~\Phi({}_k\UC)=t^{1/4}L_{-[k]}\otimes {\cal O}_{Z_n^k}\otimes L_{[k+1]}.\nn\eea
Note the shifts $t^{\pm1/4}$ are for later convenience.

\begin{proposition}\label{prop_twists}
The crossings ${}_k\OC$ and ${}_k\UC$ send a sheaf $F$ to 
\bea {}_k\OC\,F&=&t^{-1/4}\opn{cone}\big(F\oplus \iota_*q^*Rq_{*}L\iota^*F\stackrel{({\rm res},\ep)}{\longrightarrow} F|_{D_{k+1}}\big)[-1],\nn\\
{}_k\UC\,F&=&t^{1/4}L_{-[k]}\otimes\opn{cone}\big(F_{[k+1]}\oplus \iota_*q^*Rq_{*}L\iota^*F_{[k+1]}\stackrel{({\rm res},\ep)}{\longrightarrow} F_{[k+1]}|_{D_{k+1}}\big)[-1].\label{braid}\eea
We also write $q=q_{k+1}$, $\iota=\iota_{k+1}$, $F|_{D_{k+1}}:=\iota_*L\iota^*F$ and $\ep$ is the canonical map $q^*Rq_*\to id$. 
\end{proposition}
We remark that we will often work with the formula without the untidy $t^{\pm1/4}$ factors, and only restore them to a tangle at the last step.

\begin{proof}
To avoid clutter, we will sometimes omit the left/right derivation of the various functors, but do note that $q^*$ and $\iota_*$ do not need derivation as $q$ is flat and $\iota$ is a closed embedding. Also we set aside $t^{\pm1/4}$.

The FM transform acting on a sheaf $F$ has the familiar pull-push formula.\footnote{Our main reference for manipulating the FM transform is \cite{FM_AG}.} 
  \bea {}_k\OC\,F=Rp_{1*}({\cal O}_{Z_n^k}\otimes^Lp_2^*F)\nn\eea
  where $p_{1,2}$ are the projections from $Y_n\times Y_n$ to the left and right $Y_n$ factor. As $Z_n^k=\Gd\cup V_n^k$, we have
  \bea 0\to {\cal O}_{Z_n^k}\to {\cal O}_{\Gd}\oplus {\cal O}_{V_n^k}\to  {\cal O}_{\Gd\cap V_n^k}\to 0.\nn\eea
  As $Y_n$ is smooth, we can assume $F$ is a perfect complex. Tensoring the above sequence with $p_2^*F$ 
  \bea 0\to {\cal O}_{Z_n^k}\otimes p_2^*F\to {\cal O}_{\Gd}\otimes p_2^*F\oplus {\cal O}_{V_n^k}\otimes p_2^*F\to  {\cal O}_{\Gd\cap V_n^k}\otimes p_2^*F\to 0.\nn\eea
  In $D^b(Y_n)$, the left term ${\cal O}_{Z_n^k}\otimes p_2^*F$ is isomorphic to the cone
  \bea {\cal O}_{Z_n^k}\otimes p_2^*F=\opn{cone}({\cal O}_{\Gd}\otimes p_2^*F\oplus {\cal O}_{V_n^k}\otimes p_2^*F\to  {\cal O}_{\Gd\cap V_n^k}\otimes p_2^*F)[-1].\nn\eea
  The pushdown under $Rp_{1*}$ of the left term becomes
  \bea Rp_{1*}({\cal O}_{Z_n^k}\otimes p_2^*F)=\opn{cone}(\underbrace{Rp_{1*}({\cal O}_{\Gd}\otimes p_2^*F)}_F\oplus \underbrace{Rp_{1*}({\cal O}_{V_n^k}\otimes p_2^*F)}_{\mathrm{I}} \to  \underbrace{Rp_{1*}({\cal O}_{\Gd\cap V_n^k}\otimes p_2^*F)}_{\mathrm{II}})[-1].\nn\eea
  We compute each of the pushdowns. For $\mathrm{II}$, we note that $\Gd\cap V_n^k$ is $D_{k+1}$ embedded diagonally in the product $Y_n\times Y_n$, thus
  \bea \mathrm{II}=Rp_{1*}(\Gd_*{\cal O}_{D_{k+1}}\otimes_{{\cal O}_{Y_n\times Y_n}}p_2^*F)\stackrel{1}{=}Rp_{1*}\Gd_*({\cal O}_{D_{k+1}}\otimes_{{\cal O}_{Y_n}} \Gd^*p_2^*F)\stackrel{2}{=}{\cal O}_{D_{k+1}}\otimes_{{\cal O}_{Y_n}}F=F|_{D_{k+1}},\nn\eea
  where (1) is the projection formula and (2) is because $p_{1,2}\Gd=id$. Now for $\mathrm{I}=Rp_{1*}({\cal O}_{X_n^k}\otimes_{{\cal O}_{Y_n\times Y_n}}p_2^*F)$. We let $i$ be the inclusion $i:\,V_n^k\to Y_n\times Y_n$, then
  \bea \mathrm{I}=Rp_{1*}i_*(Li^*p_2^*F)=\iota_*Rp_{1*}({\cal O}_{D_{k+1}}\otimes_{{\cal O}_{Y_{n-2}}}L\iota^*F)=\iota_*q^*R_{q*}(L\iota^*F),\nn\eea
  where we wrote $p_{1,2}i=\iota p_{1,2}$, and the last step is the flat base change associated with the pullback diagram 
  \bea \begin{tikzpicture}
  \matrix (m) [matrix of math nodes, row sep=2em, column sep=2em]
    { D_{k+1}\times_{Y_{n-2}}D_{k+1}  & D_{k+1}  \\
      D_{k+1} & Y_{n-2} \\ };
  \path[->,font =\scriptsize]
  (m-1-1) edge node[above]{$p_2$} (m-1-2)
  (m-1-1) edge node[right]{$p_1$} (m-2-1);
  \path[->>,font =\scriptsize]
  (m-1-2) edge node[right] {$q$} (m-2-2)
  (m-2-1) edge node[above] {$q$} (m-2-2);
\end{tikzpicture}\nn\eea
Finally, putting back the shift, we conclude
\bea {}_k\OC F=\opn{cone}(F\oplus \iota_*q^*Rq_*(F|_{D_{k+1}})\to F|_{D_k+1})[-1].\nn\eea
\end{proof}
Now for the cups and caps 
\bea {}^i\CUP:\;D^b(Y_{n-2}) \to D^b(Y_n),~~~~{}_i\CAP:\; D^b(Y_n) \to D^b(Y_{n-2}).\nn\eea
Following sec.4.2.1 of \cite{Cautis2007KnotHV}, we set the kernels for the cup and cap to be
\bea\Phi(
{}^k\CUP)=t^{-\frac{k}{2}}L_{-[k]}\otimes{\cal O}_{D_{k+1}},~~~\Phi({}_k\CAP)=t^{\frac{k+1}{2}}{\cal O}_{D_{k+1}}\otimes L_{[k+1]},\nn\eea
where the first ${\cal O}_{D_{k+1}}$ should be regarded as being pushed forward into $Y_{n+2}\times Y_n$ via the map $(\iota_{k+1},q_{k+1})$, while the second ${\cal O}_{D_{k+1}}$ is similar but with the $Y_{n+2}$ and $Y_n$ reversed. So they are essentially the structure sheaf of the graph of $q_{k+1}$.
\begin{proposition}
The FM transforms act on a sheaf $F$ as
\bea {}^k\CUP F=t^{-\frac{k}{2}}L_{-[k]}\otimes\iota_*q^*F,~~~{}_k\CAP F=t^{\frac{k+1}{2}}Rq_*L\iota^*(F\otimes L_{[k+1]}),\nn\eea
where $q=q_{k+1}$ and $\iota=\iota_{k+1}$.
\end{proposition}
\begin{proof}
  That the FM transform using the structure sheaf of the graph of a map $f$ gives $Lf^*$ or $Rf_*$ is standard, see ex.5.4 of \cite{FM_AG}.
\end{proof}
\begin{proposition}\label{prop_cap_cup}
The caps and cups satisfy 
\bea &&\CAP_k\;{}^k\CUP =id={}_{k+1}\CAP\;{}^k\CUP,\nn\\
&& {}_k\CAP\;{}^k\CUP=(t^{1/2}[1]\oplus t^{-1/2}[-1]).\nn\eea
In pictures
\bea 
\begin{tikzpicture}[scale=0.7,baseline={([yshift=-.5ex]current bounding box.center)}]
\draw[blue] (0.6,.6) -- (0.6,0);
\centerarc[blue](0.3,0)(-180:0:0.3);
\centerarc[blue](-0.3,0)(180:0:0.3);
\draw[blue] (-0.6,0) -- (-0.6,-0.6);
\end{tikzpicture}
\quad = \quad 
\begin{tikzpicture}[scale=0.7,baseline={([yshift=-.5ex]current bounding box.center)}]
\draw[blue] (0,0.6) -- (0,-0.6);
\end{tikzpicture}
\quad = \quad
\begin{tikzpicture}[scale=0.7,baseline={([yshift=-.5ex]current bounding box.center)}]
\draw[blue] (-0.6,.6) -- (-0.6,0);
\centerarc[blue](-0.3,0)(-180:0:0.3);
\centerarc[blue](0.3,0)(180:0:0.3);
\draw[blue] (0.6,0) -- (0.6,-0.6);
\end{tikzpicture}~, \qquad
\begin{tikzpicture}[scale=0.7,baseline={([yshift=-.5ex]current bounding box.center)}]
\centerarc[blue](0.3,0)(-180:0:0.3);
\centerarc[blue](0.3,0)(180:0:0.3);
\end{tikzpicture} \quad = \quad t^{1/2}[1]\oplus t^{-1/2}[-1].
\nn\eea
There are also the following relations
\bea  {}_k\CAP\;{}^l\CUP={}^{l-2}\CUP\;{}_k\CAP,~~~{}^l\CAP\;{}_k\CUP={}^k\CUP\;{}_{l-2}\CAP,~~~
{}_l\CAP\;{}_k\CAP={}_k\CAP\;{}_{l+2}\CAP,~~~{}^l\CUP\;{}^k\CUP={}^k\CUP\;{}^{l-2}\CUP,~~~l\geq k+2.\nn\eea
We leave it to the reader to draw the pictures and see why they are expected.
\end{proposition}
\begin{proof}
Using the formula for the action of cup and cap,
\bea \CAP_k\;{}^k\CUP F=Rq_{k,*}L\iota_k^*(L_{[k]}\otimes L_{-[k]}\otimes \iota_{k+1,*}q_{k+1}^*F)
=Rq_{k,*}L\iota_k^*\iota_{k+1,*}q_{k+1}^*F.\nn\eea
In the rhs, we pull $F$ back to $D_{k+1}$ then restrict to $D_{k+1}\cap D_k$. 
From lem.\ref{lem_trivial_lb}, the two projections $q_k,q_{k+1}$ from $D_k\cap D_{k+1}$ are both isomorphisms and so the pull and push gives $F$ back.

For ${}_{k+1}\CAP\;{}^k\CUP F$,
\bea {}_{k+1}\CAP\;{}^k\CUP F=tRq_{k+2,*}L\iota_{k+2}^*(L_{[k+2]}\otimes L_{-[k]}\otimes \iota_{k+1,*}q_{k+1}^*F).\nn\eea
The line bundle $L_{[k+2]}\otimes L_{-[k]}$ is trivial over $D_{k+1}\cap D_{k+2}$, more precisely, it is $t^{-1}{\cal O}_{D_{k+1}\cap D_{k+2}}$, as shown in lem.\ref{lem_trivial_lb}. Having replaced $L_{[k+2]}\otimes L_{-[k]}$, the rest of the argument is identical to the one above.

Finally ${}_k\CAP\;{}^k\CUP$ gives the unknot
\bea {}_k\CAP\;{}^k\CUP F=t^{1/2} q_{k+1,*}({\cal O}_{D_{k+1}}(D_{k+1})\otimes\iota_{k+1}^*\iota_{k+1,*}q_{k+1}^*F).\nn\eea
Note that $\iota_{k+1}^*\iota_{k+1,*}$ is given by tensoring over ${\cal O}_{D_{k+1}}$ with the complex ${\cal O}_{D_{k+1}}(-D_{k+1})\to  {\cal O}_{D_{k+1}}$ (the rightmost term has degree 0). Thus
\bea{}_k\CAP\;{}^k\CUP F&=&t^{1/2} q_{k+1,*}(({\cal O}_{D_{k+1}}\to{\cal O}_{D_{k+1}}(D_{k+1}))\otimes q_{k+1}^*F)\nn\\
&=&t^{1/2} q_{k+1,*}({\cal O}_{D_{k+1}}\to{\cal O}_{D_{k+1}}(D_{k+1}))\otimes F\nn\\
&=&t^{1/2} ({\cal O}_{Y_{n-2}}\to t^{-1}{\cal O}_{Y_{n-2}}[-1])\otimes F=(t^{1/2}F[1]\oplus t^{-1/2}F[-1]).\nn\eea
Note that to compute $q_{k+1,*}{\cal O}_{D_{k+1}}(D_{k+1})$, one can use the short exact sequence $0\to {\cal O}_{D_{k+1}}(D_{k+1})\to {\cal O}_{D_{k+1}}(D_k+D_{k+1})\to {\cal O}_{D_{k+1}}(D_k+D_{k+1})|_{D_k}\to0$. The middle term vanishes under $q_{k+1,*}$ as it has fibre degree $-1$ (see prop.\ref{prop_pushdown_q}), and so $q_{k+1,*}{\cal O}_{D_{k+1}}(D_{k+1})=q_{k+1,*}\iota_{k,*}{\cal O}_{D_k\cap D_{k+1}}(D_k+D_{k+1})[-1]=t^{-1}{\cal O}_{Y_{n-2}}$, where we used item 3 of lem.\ref{lem_trivial_lb}.

For the remaining relations, we prove one of them: ${}_k\CAP\;{}^l\CUP={}^{l-2}\CUP\;{}_k\CAP$, which looks like
\bea 
\begin{tikzpicture}[scale=0.7]
\draw[blue] (0.6,.3) -- (0.6,0);
\centerarc[blue](0.3,0)(-180:0:0.3);
\draw[blue] (0,.3) -- (0,0);

\draw[blue] (-0.8,-.3) -- (-0.8,0);
\centerarc[blue](-0.5,0)(0:180:0.3);
\draw[blue] (-0.2,-.3) -- (-0.2,0);

\node[blue] at (1.1,-0.1) {\small{$=$}};

\draw[blue] (3.1,.6) -- (3.1,0.3);
\centerarc[blue](2.8,0.3)(-180:0:0.3);
\draw[blue] (2.5,.6) -- (2.5,0.3);

\draw[blue] (1.7,-.8) -- (1.7,-.5);
\centerarc[blue](2,-0.5)(0:180:0.3);
\draw[blue] (2.3,-.8) -- (2.3,-0.5);
\end{tikzpicture}\nn\eea
Set $l\geq k+2$
\bea {}^{l-2}\CUP\; {}_k\CAP F=t^{(k-l+3)/2} L_{-[l-2]}\iota_{l-1,*}q_{l-1}^*q_{k+1,*}\iota_{k+1}^*(L_{[k+1]}\otimes F).\nn\eea
The sequence of maps starts from the top left in the following diagram and goes along the blue route to the bottom right.
\bea \begin{tikzpicture}
  \matrix (m) [matrix of math nodes, row sep=2em, column sep=2em]
    { Y_n & D_{l+1}^{n+2} & Y_{n+2} \\
      D^n_{k+1} & D^{n+2}_{k+1}\cap D^{n+2}_{l+1} & D^{n+2}_{k+1} \\ 
      Y_{n-2} & D^n_{l-1} & Y_n \\    };
  \path[->,font =\scriptsize]
  (m-1-2) edge[red] node[above]{$q_{l+1}$} (m-1-1)
  (m-2-2) edge node[above]{$q_{l+1}$} (m-2-1)
  (m-3-2) edge[blue] node[above]{$q_{l-1}$} (m-3-1)
  (m-2-1) edge[blue] node[left]{$q_{k+1}$} (m-3-1)
  (m-2-2) edge node[right]{$q_{k+1}$} (m-3-2)
  (m-2-3) edge[red] node[right]{$q_{k+1}$} (m-3-3);
  \path[>->,font =\scriptsize]
  (m-1-2) edge[red] node[above] {$\iota_{l+1}$} (m-1-3)
  (m-2-2) edge node[above] {$\iota_{l+1}$} (m-2-3)
  (m-2-1) edge[blue] node[left] {$\iota_{k+1}$} (m-1-1)
  (m-3-2) edge[blue] node[above] {$\iota_{l-1}$} (m-3-3)
  (m-2-2) edge node[left] {$\iota_{k+1}$} (m-1-2)
  (m-2-3) edge[red] node[left] {$\iota_{k+1}$} (m-1-3);
\node at (-1.15,-0.7) {\small{\eqref{pullback_square}}};
\end{tikzpicture}\nn\eea
As all squares are pullbacks, applying the flat base change formula along these squares, we can move from the blue route to the red one
\bea &&L_{-[l-2]}\iota_{l-1}^*q_{l-1}^*q_{k+1,*}\iota_{k+1}^*L_{[k+1]}F=L_{-[l-2]}q_{k+1,*}\iota_{k+1}^*\iota_{l+1,*}q_{l+1}^*L_{[k+1]}F\nn\\
&=&L_{-[l-2]}q_{k+1,*}\iota_{k+1}^*L_{[l]}L_{[k+1]}L_{-[l]}\iota_{l+1,*}q_{l+1}^*F
=L_{-[l-2]}q_{k+1,*}(\iota_{k+1}^*L_{[l]})\iota_{k+1}^*L_{[k+1]}L_{-[l]}\iota_{l+1,*}q_{l+1}^*F.\nn\eea
Now $\iota_{k+1}^*L_{[l]}=t^{-1}q_{k+1}^*L_{[l-2]}$ from prop.\ref{prop_pushdown_q}, then the projection formula allows us write the above as
$t^{-1}q_{k+1,*}\iota_{k+1}^*L_{[k+1]}L_{-[l]}\iota_{l+1,*}q_{l+1}^*F$, which is 
which proves ${}^{l-2}\CUP\; {}_k\CAP={}_k\CAP\;{}^l\CUP$.
\end{proof}
The relations here plus the Kauffman triangle of the next section would be enough to prove all the necessary moves that guarantee that we get an isotopy invariant of the tangles.

\subsection{The Kauffman Distinguished Triangle}
The usual Kauffman bracket \cite{LHKAUFFMAN1987395} relates a crossing with the two ways of resolving it; it is a highly efficient algorithm to compute a knot polynomial. Here, the Kauffman bracket gets promoted to distinguished triangles in $D^b(Y_n)$.
\begin{proposition}\label{prop_Kauffman}
  There are two distinguished triangles that categorify the Kauffman brackets
  \bea 
  &&F \to t^{-\frac14}{}_k\UC F\to t^{-\frac12}\;{}^k\CUP\;{}_k\CAP F\stackrel{+1}{\to},\nn\\
  &&t^{\frac12}{}^k\CUP\;{}_k\CAP F\to t^{\frac14}{}_k\OC F\to F\stackrel{+1}{\to}.\nn\eea
  In pictures they read
  \bea 
  \begin{tikzpicture}[scale=1]
\draw [blue,-] (2,0) -- (3,1);
\draw [white,-,line width=2mm] (3,0) -- (2,1);
\draw [blue,-] (3,0) -- (2,1);

\node at (1.6,0.5) {$t^{-\frac14}$};
\node at (1,0.5) {\scriptsize$\longrightarrow~~$};

\centerarc[blue](-0.5,0.5)(65:-65:0.5);
\centerarc[blue](0.8,0.5)(115:245:0.5);

\node at (3.8,0.5) {{\scriptsize$\longrightarrow$}  $t^{-\frac12}$};
\centerarc[blue](5,1.2)(-35:-145:0.5);
\centerarc[blue](5,-0.2)(35:145:0.5);
\node at (5.9,0.5) {\scriptsize{$\stackrel{+1}{\longrightarrow}$}};
\end{tikzpicture}\nn\\
\begin{tikzpicture}[scale=1]
\draw [blue,-] (3.5,0) -- (2.5,1);
\draw [white,-,line width=2mm] (2.5,0) -- (3.5,1);
\draw [blue,-] (2.5,0) -- (3.5,1);

\node at (-0.4,0.5) {$t^{\frac12}$};

\centerarc[blue](4.7,0.5)(65:-65:0.5);
\centerarc[blue](6,0.5)(115:245:0.5);

\node at (2,0.5) {$t^{\frac14}$};
\node at (1.3,0.5) {\scriptsize$\longrightarrow~~$};
\node at (4,0.5) {\scriptsize$\longrightarrow$};
\node at (6.2,0.5) {\scriptsize$\stackrel{+1}{\longrightarrow}$};

\centerarc[blue](.5,1.2)(-35:-145:0.5);
\centerarc[blue](.5,-0.2)(35:145:0.5);

\end{tikzpicture}\nn
\eea
\end{proposition}
\begin{proof}
    From prop.\ref{prop_twists}, the $k^{th}$-twists are defined by the following DT's (as usual, we set aside the fractional powers of $t$)
    \bea &&F\oplus \iota_*q^*q_{*}\iota^*F\stackrel{({\rm res},\ep)}{\longrightarrow} F|_{D_{k+1}}\to {}_k\OC F[1]\stackrel{+1}{\to},\nn\\
    &&F(D_{k+1})\oplus L_{-[k]}\otimes\iota_*q^*q_{*}\iota^*F_{[k+1]}\stackrel{({\rm res},\ep)}{\longrightarrow} F(D_{k+1})|_{D_{k+1}}\to {}_k\UC F[1]\stackrel{+1}{\to},\nn  \eea
    where $F|_{D_{k+1}}=\iota_*\iota^*F$ and $F_{[k+1]}=F\otimes L_{[k+1]}$. 
    We deal with ${}_k\UC F$ first, which is easier. Using the octahedral axiom, we can rewrite the DT as
    \bea &&L_{-[k]}\otimes\iota_*q^*q_{*}\iota^*F_{[k+1]}\to\opn{cone}(F(D_{k+1})\to \iota_*\iota^*F(D_{k+1}))\to {}_k\UC F[1]\stackrel{+1}{\to}\nn\\
    \To&&L_{-[k]}\otimes\iota_*q^*q_{*}\iota^*F_{[k+1]}\to F[1]\to {}_k\UC F[1]\stackrel{+1}{\to}.\nn\eea
    In computing $\opn{cone}(F(D_{k+1})\to \iota_*\iota^*F(D_{k+1}))$ we used the fact that $\iota_*\iota^*F(D_{k+1})$ fits into a short exact sequence $0\to F\to F(D_{k+1})\to \iota_*\iota^*F(D_{k+1})\to0$.

    Now we get the following DT after a rotation
    \bea
    F \to {}_k\UC F\to L_{-[k]}\otimes\iota_*q^*q_{*}\iota^*F_{[k+1]}\stackrel{+1}{\to}.\nn\eea
    We recognise the rightmost term to be $t^{-1/2}\;{}^k\CUP\;{}_k\CAP F$, restoring the $t$-factors, we get the first Kauffman relation.

    Now for ${}_k\OC F$. Using again the octahedral axiom and the fact that $\iota_*$ is exact and requires no derivation  
    \bea F\to \iota_*\opn{cone}(q^*q_{*}\iota^*F\to \iota^*F) \to {}_k\OC F[1]\stackrel{+1}{\to}.\label{Kauffman_DT}\eea
    We compute the middle cone next. We have that $q_*q^*=1$, which follows from the projection formula and that the fibrewise pushdown of the structure sheaf is still the structure sheaf when the fibre is projective. We conclude that $\opn{cone}(q^*q_{*}\iota^*F\to \iota^*F)$ is killed by $q_*$ and has to be `proportional to'
    a line bundle of fibre degree $-1$, in the sense that we have an isomorphism
    \bea \opn{cone}(q^*q_{*}\iota^*F\to \iota^*F)\cong q_*\iota^*(F\otimes L_{[k+1]})\otimes_{{\cal O}_{Y_{n-2}}}t\iota^*L_{-[k]}[1].\label{proportional}\eea
    where $\iota^*L_{-[k]}$ has fibre degree $-1$ by prop.\ref{prop_pushdown_q}.
    Assuming this for the moment, \eqref{Kauffman_DT} becomes 
    \bea F\to t^{\frac12}\,{}^k\CUP\,{}_k\CAP F[1]\to {}_k\OC F[1]\stackrel{+1}{\to}.\nn\eea
    Rotate this triangle and we are done.

    Now for the proof \eqref{proportional}.  
    As $D_{k+1}$ is a $\BB{P}^1$-bundle over $Y_{n-2}$, we can 
    apply the fibrewise tilting complex. For this, we use the structure sheaf ${\cal L}_0:={\cal O}_{D_{k+1}}$ and ${\cal L}_1:=\iota^*L_{-[k+1]}$ as the two line bundles that restrict to degree 0 and 1 on the $\BB{P}^1$-fibre. 
    Fibrewise tilting says that a sheaf $G$ on $D_{k+1}$ is isomorphic to the cone of
    \bea \begin{tikzpicture}
    \matrix (m) [matrix of math nodes, row sep=1em, column sep=2em]
    { q_*(G\otimes {\cal L}_1^{\vee})\otimes_{{\cal O}_{Y_{n-2}}}q_*{\cal L}_1\otimes_{{\cal O}_{Y_{n-2}}} {\cal L}_0 & q_*G\otimes_{{\cal O}_{Y_{n-2}}}{\cal L}_0 \\
        & q_*(G\otimes{\cal L}_1^{\vee})\otimes_{{\cal O}_{Y_{n-2}}}{\cal L}_1 \\   };
    \path[->,font =\scriptsize]
    (m-1-1) edge (m-1-2)
    (m-1-1) edge (m-2-2);
    \node at (3,0) {$\oplus$};
    \end{tikzpicture}\label{replace_G}\eea
    Now apply \eqref{replace_G} to $G=\opn{cone}(q^*q_{*}\iota^*F\to \iota^*F)$, and noting $q_*G=0$ in this case
    \bea \opn{cone}(q^*q_{*}\iota^*F\to \iota^*F)\simeq q_*(\iota^*F\otimes{\cal L}_1^{\vee})\otimes_{{\cal O}_{Y_{n-2}}}\opn{cone}\big(q_*{\cal L}_1\otimes_{{\cal O}_{Y_{n-2}}} {\cal L}_0 \rightarrow {\cal L}_1\big).
    \label{replace_F}\eea
    One recognises from the right hand cone the fibrewise Euler-sequence, indeed, this cone is isomorphic to $t\iota^*L_{-[k]}[1]$. 
    This can be seen by simply replacing $G$ with $\iota^*L_{-[k]}$ (which has fibre degree $-1$) in the very same complex \eqref{replace_G}.
    As $q_*\iota^*L_{-[k]}=0$, and so we get that $\iota^*L_{-[k]}$ is isomorphic to
    \bea q_*\iota^*({\cal O}(D_{k+1}))\otimes_{{\cal O}_{Y_{n-2}}}\opn{cone}\big(q_*{\cal L}_1\otimes_{{\cal O}_{Y_{n-2}}} {\cal L}_0
    \to {\cal L}_1\big).\nn\eea
    The pushdown $q_*\iota^*({\cal O}(D_{k+1}))$ is a trivial line bundle $t^{-1}{\cal O}_{Y_{n-2}}[-1]$ over $Y_{n-2}$, up to a shift (see ex.\ref{ex_for_pullback}). 
    Now \eqref{replace_F} can be written as
    \bea \opn{cone}(q^*q_{*}\iota^*F\to \iota^*F)\simeq q_*\iota^*(F\otimes L_{[k+1]})\otimes_{{\cal O}_{Y_{n-2}}}t\iota^*L_{-[k]}[1], \nn\eea
    which is \eqref{proportional}. This proves the claim.
\end{proof}

\subsection{Rouquier Construction}\label{sec_RC}
Here we take a detour to relate our treatment to Rouquier's formulation using adjoint functors \cite{rouquier2006categorification}.
The maps in the DTs of prop.\ref{prop_Kauffman} follow implicitly from the axioms of derived categories. They can also be viewed as units and counits. From cor.\ref{cor_2_side_adj} we have the adjoints 
\bea t^{1/2}\CAP_k[1] \dashv \CUP_k \dashv t^{-1/2}\CAP_k[-1],\nn\eea  
so we rewrite the Kauffman DTs as
\bea %
\begin{tikzpicture}[scale=1]
\draw [blue,-] (2,0) -- (3,1);
\draw [white,-,line width=2mm] (3,0) -- (2,1);
\draw [blue,-] (3,0) -- (2,1);

\node at (1.3,0.5) {\scriptsize$\longrightarrow~~$};

\draw[blue](-0.2,0) -- (-0.2,1);
\draw[blue](0.6,0) -- (0.6,1);

\node at (-0.5,0.6) {\scriptsize$\stackrel{\ep}{\longrightarrow}~~$};
\node at (-1.3,0.1) {\scriptsize$[-1]$};

\node at (-2.8,0.1) {\scriptsize$t^{-1/2}$};
\centerarc[blue](-2,1.2)(-35:-145:0.5);
\centerarc[blue](-2,-0.2)(35:145:0.5);
\node at (3.5,0.6) {\scriptsize{$\stackrel{+1}{\longrightarrow}$}};
\node at (2,0.5) {\scriptsize$t^{-1/4}$};
\end{tikzpicture}\nn\\
\begin{tikzpicture}[scale=1]
\draw [blue,-] (3,0) -- (2,1);
\draw [white,-,line width=2mm] (2,0) -- (3,1);
\draw [blue,-] (2,0) -- (3,1);
\node at (3,0.5) {\scriptsize$[1]$};

\node at (-0.5,0.1) {\scriptsize$t^{1/2}$};
\node at (0.8,0.1) {\scriptsize$[1]$};

\draw[blue](-2.3,0) -- (-2.3,1);
\draw[blue](-1.5,0) -- (-1.5,1);

\node at (-1,0.6) {\scriptsize$\stackrel{\eta}{\longrightarrow}~~$};
\node at (1.4,0.5) {\scriptsize$\longrightarrow$};
\node at (3.5,0.5) {\scriptsize$\stackrel{+1}{\longrightarrow}$};

\centerarc[blue](.2,1.2)(-35:-145:0.5);
\centerarc[blue](.2,-0.2)(35:145:0.5);

\node at (2,0.5) {\scriptsize$t^{1/4}$};
\end{tikzpicture}\nn\eea
where the lower cap with $t^{-1/2}[-1]$ (resp. $t^{1/2}[1]$) is the right (resp. left) adjoint to the cup.
With this view, one has immediately:
\begin{proposition}
  The two crossings defined above as cones of units/counits are inverse isomorphisms. 
\end{proposition}
\begin{proof}
  We apply prop.8.1 of \cite{rouquier2006categorification}, for which we set
  \bea &&{\cal C}=D^b(Y_{n-2})~\mathop{\rightleftharpoons}^F_G~{\cal D}=D^b(Y_n),\nn\\
  &&F=\CUP,~~G=t^{-1/2}\CAP[-1],~~\Phi=t[2],\nn\\
  &&\Upsilon=\opn{cone}(FG\to 1_{\cal D})[-1],~~\Upsilon'=\opn{cone}(1_{\cal D}\to F\Phi G).\nn\eea
  This fulfills the assumption that there is a sequence of natural transformations
  \bea 1_{\cal C}\to GF\to \Phi^{-1} \stackrel{+1}{\to}, \nn \eea 
  the first map being the unit of the adjunction, which gives a `DT of functors' in the sense that we get a distinguished triangle upon applying it to any individual object in ${\cal C}$.
  
  In fact, we have $GF=t^{-1/2}\bigcirc[-1]=t^{-1/2}(t^{1/2}[1]\oplus t^{-1/2}[-1])[-1]=1\oplus t^{-1}[-2]$, so the triangle above is
  \bea 1\to 1\oplus t^{-1}[-2] \to t^{-1}[-2].\nn\eea
  Then prop.8.1 of loc.cit. says that $\Upsilon$ and $\Upsilon'$ are inverse equivalences. In fact, $\Upsilon$ is the undercrossing shifted by $[-1]$ and $\Upsilon'$ is the overcrossing shifted by [1].
\end{proof}

\smallskip

In practice, this description of the maps in the DTs in terms of counit and unit transformation is still too abstract to be effective; nevertheless, when we apply the DT to a concrete object in $D^b(Y_n)$, we can often determine the maps in an ad hoc manner, via some simple tests relying on equivariance and degree considerations. This point will be illustrated in the calculations given in sec.\ref{sec_CS}. 

We have yet to show that the crossings, cap/cups above respect a list of isotopy moves so that one obtains a tangle invariant. They are usually easy consequences of the Kauffman DT and prop.\ref{prop_cap_cup}. We will prove the key isotopy moves in secs.\ref{sec_RIM}, \ref{sec_PosIIS} once we have seen some sample calculations in sec.\ref{sec_CS} to thaw the fingers.

\section{Some Calculation Examples}\label{sec_CS}
In this section, we illustrate our calculation methods with some examples. 
The bag of tricks is actually quite limited: apart from the standard projection formula, and flat base change, we will repeatedly use
\begin{enumerate}
\item[$\spadesuit$]
For a sheaf $F$ on $Y_n$, $F|_D=\iota_*\iota^*F=\opn{cone}(F(-D)\to F)$, where $\iota:\,D\hookrightarrow Y_n$ is the inclusion. In contrast, if $G$ is a sheaf on $D$, then $\iota^*\iota_*G=(G(-D)[1]\oplus G)$. \label{item_spade}
\item[$\heartsuit$] Rather than all sheaves, we will only need certain line bundles $L_{i_1,\cdots,i_n}$ given in def.\ref{def_exp_collection}. Their sheaf cohomology can be computed quickly as shown in ex.\ref{ex_L_sheaf_coh}.\label{item_heart}
  \item[$\diamondsuit$] The formula for pushdown, pullback along the projection $q_k:\, D_k\to Y_{n-2}$, as well as the fibre degree for the line bundles $L_{i_1,\cdots,i_n}$, given in prop.\ref{prop_pushdown_q}.\label{item_diamond}
  \item[$\clubsuit$] Restricting $L_{i_1,\cdots,i_n}$ to $D_k$, we can adjust its indices to $(s/t^{k-1})L_{i_1-2,\cdots,i_{k-1}-2,i_k-1,i_{k+1},\cdots,i_n}$, shown in lem.\ref{lem_trivial_lb}.\label{item_club}
\end{enumerate}
\begin{remark}\label{rem_BWB}
We would like to make a final remark about the organisation of sheaf cohomologies.
    The varieties $Y_n$ are iterated $\mathbb{P}^1$-fibrations, we can better organise the cohomology by means of the Borel--Weil--Bott theorem: using $\BB{P}^1=SU(2)/U(1)=GL_2/B$, we identify the vector space $H^0(\BB{P}^1,{\cal O}(n))$ with the spin $n/2$ representation of $SU(2)$ (of dimension $n+1$). Remembering the equivariant degrees we write
    \[ H^0(\BB{P}^1,{\cal O}(1))=V^*=1+s, \quad H^1(\BB{P}^1,s{\cal O}(-3))=H^1(\BB{P}^1,K(-1))=V=1+s^{-1}, \]  
    which is consistent with Serre duality. Note the convention of writing $V^*=1+s$ was explained in the beginning of sec.\ref{sec_BoY}.
    
    We also have $sV \cong V^*$ as the fundamental representation of $SU(2)$ is pseudo real, but the isomorphism is non-trivial: if $\vec u=u_1+u_2s^{-1}$, then $s\vec u$ is sent to $(u_2,-u_1s)\in V^*$. This will be responsible for making certain differentials square to zero in our calculations below.
\end{remark}

\subsection{Maps of the Kauffman Triangle}
In this section, we address the issuing of pinning down the maps in the Kauffman triangles in prop. \ref{prop_Kauffman}.
As we just mentioned above, there is often only one possible candidate for these maps. Thus we need only decide whether the map is zero or not, which can be done by looking at the associated long exact sequences of sheaf cohomology.
We illustrate our point by investigating whether the map $t^{1/4}\OC F\to F$ (of 0 degree and 0 equivariant weight).

To do this, we observe that in the original formula \eqref{braid} for the crossings acting on an object $F$, all maps involved are canonical, so we can fairly easily compute the equivariant sheaf cohomology $H^{\sbullet}(Y_n,L\otimes \OC F)$ of $F$ tensored with any line bundle. The DT in prop.\ref{prop_Kauffman} implies a long exact sequence
\bea H^{\sbullet-1}(Y_n,L\otimes F)\to t^{\frac12}H^{\sbullet}(Y_n,L\otimes\CUP\CAP F) \to  t^{\frac14}H^{\sbullet}(Y_n,L\otimes \OC F)\to H^{\sbullet}(Y_n,L\otimes F)\to t^{\frac12}H^{\sbullet+1}(Y_n,L\otimes\CUP\CAP F).\nn\eea
As the example below shows, the third arrow is nonzero for some $L$, then the map $t^{1/4}\OC F\to F$ cannot be zero, thereby fixing it up to a scalar.
\begin{example}\label{ex_DT_map}
    We investigate the two-strand case (that is, $n=2$). In appendix sec. \ref{sec_ECoDY}, an exceptional collection made up of four line bundles $L_{00},L_{10},L_{21},L_{31}$ is given for the category $D^b(Y_2)$, where we denote $L_{pq}={\cal O}(pD_1+qD_2)$. We now look at $\OC$ on each of these line bundles. First, for $L_{00}$, the formula for $\OC$ gives
    \bea t^{1/4}\OC L_{00}=\opn{cone}(L_{00}\oplus \iota_*q^*q_*\iota^*L_{00}\to \iota_*\iota^*L_{00})[-1],\nn\eea
    where $q=q_2$ and $\iota=\iota_2$ and the two maps are restriction and evaluation, so their effect is easily determined. It is easy to compute each term and it turns out that $\OC L_{00}=L_{00}$. In the same way we calculate that $\OC L_{21}=L_{21}$. But for $L_{10},L_{31}$, the story is different. Prop.\ref{prop_Kauffman} gives the DT
    \bea L_{10}\to t^{1/2}\CUP \CAP L_{10}[1]\to t^{1/4}\OC L_{10}[1]\stackrel{+1}{\to}.\label{DT_L10}\eea
    To understand the first map, we note that we have
    \bea \hom^0(L_{10},t^{1/2}\CUP \CAP L_{10}[1])=\hom^0(L_{10},s\iota^*L_{-10}[1])=\ext^1(L_{10},s\iota^*L_{-10})=1,\nn\eea
    where $\CUP\CAP L_{10}$ is computed as (use the formula for $q_*$ and $q^*$ provided in prop.\ref{prop_pushdown_q}) 
    \bea \CUP\CAP L_{10}=t\CUP q_*\iota^*L_{21}=\CUP \frac{s}{t}=st^{-1/2}\iota_*\iota^*L_{-10}.\nn\eea

    Thus the map $L_{10}\to t^{1/2}\CUP \CAP L_{10}[1]$ is either zero or given by the above 1-dimensional $\ext^1$. In order to determine which case it is, we investigate the sheaf cohomology long exact sequence
    \bea H^{\sbullet}(Y_2,L_{10})\to t^{1/2}H^{\sbullet}(Y_2,\CUP \CAP L_{10}[1])\to t^{1/4}H^{\sbullet}(Y_2,\OC L_{10}[1])\to H^{\sbullet+1}(Y_2,L_{10})\to \cdots\nn\eea
    Here the second term is zero, this sequence is exact regardless whether $L_{10}\to t^{1/2}\CUP\CAP L_{10}[1]$ is zero or not. 
    But since tensoring with a locally free sheaf preserves the DT, and if by inspecting the new long exact sequence
    \bea H^{\sbullet}(Y_2,E\otimes L_{10})\to t^{1/2}H^{\sbullet}(Y_2,E\otimes \CUP \CAP L_{10}[1])\to H^{\sbullet}(Y_2,E\otimes \OC L_{10}[1])\to H^{\sbullet+1}(Y_2,E\otimes L_{10})\to \cdots\nn\eea
    we get that the first map is nonzero, then the map $L_{10}\to t^{1/2}\cup \cap L_{10}[1]$ cannot be zero.
    We pick a line bundle $E=L_{-1,0}$, and now 
    \bea &&H^{\sbullet}(Y_2,E\otimes \OC L_{10}[1])=H^{\sbullet}(Y_2,\opn{cone}(L_{00}\oplus L_{-1,0}\otimes \iota_*q^*q_*\iota^*L_{10} \to \iota_*\iota^*L_{00}))=0\nn\\
    &&t^{1/2}H^{\sbullet}(Y_2,E\otimes \CUP \CAP L_{10}[1])=H^{\sbullet}(Y_2,s\iota_*\iota^*L_{-20}[1])=1,\nn\eea
    where $H^{\sbullet}(Y_2,L_{-1,0}\otimes \iota_*q^*q_*\iota^*L_{10})=H^{\sbullet}(Y_2,(1+s)L_{-1,0}\otimes \iota_*\iota^*L_{00})=(1+s)H^{\sbullet}(Y_2,\iota_*\iota^*L_{-1,0})=0$, while the map $H^{\sbullet}(Y_2,L_{00})\to H^{\sbullet}(Y_2,\iota_*\iota^*L_{00})$ is an isomorphism, causing the first line to vanish.
    The long exact sequence becomes
    \bea 1\to 1\to 0 \to [1]\to [1] \cdots\nn\eea
    which is inconsistent with the first map being zero. 
    One can draw similar conclusions for $\OC L_{31}$.

    To better present the action of $\OC$ on $L_{\sbullet\sbullet}$, we express the result using the tilting complex (we shall explain this more thoroughly in sec.\ref{sec_LthHC}). For example we have determined the maps in \eqref{DT_L10} and $t^{1/2}\cup\cap L_{10}=s\iota_*\iota^*L_{-10}$ is quasi-isomorphic to the complex of exceptional objects
    \bea s\iota_*\iota^*L_{-10} \stackrel{q.i.}{\sim} (tL_{10}\to t(1+s^{-1})L_{21}\to ts^{-1}L_{31})\nn\eea
    with the middle term at degree 0. So the map in $\opn{cone}(L_{10}[-1]\to t^{1/2}\cup \cap L_{10})$ is induced by the unique map $L_{10}\to st^{-1}L_{31}$ of weight zero.
    We summarise the computation for all the exceptional objects in the following picture:
    \bea  \begin{tikzpicture}
    \matrix (m) [matrix of math nodes, row sep=1em, column sep=2.5em]
    { L_{00} & & L_{00} & \\ 
    L_{10} & tL_{10} & 
    \begin{matrix}
        \frac{t}{s}(1+s)L_{21} \\ L_{10}
    \end{matrix} &
    \frac{t}{s}L_{31} \\ 
    L_{21} & & L_{21} & \\
    L_{31} & sL_{10} & \begin{matrix} (1+s)L_{21} \\ L_{31} \end{matrix} & L_{31} \\ };
    \node[above left=-0.6 and -0.5 of m-2-3] (x1) {\phantom{x}};
    \node[above right=-0.6 and -0.5 of m-2-3] (x2) {\phantom{x}};
    \node[below=-0.5 of m-2-3] (y) {\phantom{y}};
    \node[above right =-0.6 and -0.7 of m-2-3] (y1) {\phantom{x}};
    \node[above left=-0.6 and -0.5 of m-4-3] (z1) {\phantom{x}};
    \node[above right=-0.6 and -0.5 of m-4-3] (z2) {\phantom{x}};
    \node[below=-0.5 of m-4-3] (w) {\phantom{w}};
    \node[above right =-0.6 and -0.7 of m-4-3] (w1) {\phantom{w}};
    \path[->,font =\scriptsize]
    (x2) edge  node[above]{\textcircled{b}} (m-2-4)
    (y) edge node[below]{\textcircled{c}} (m-2-4)
    (m-2-2) edge node[above]{\textcircled{a}} (x1)
    (m-4-2) edge (z1)
    (z2) edge (m-4-4);
    \path[->,font =\scriptsize]
    (m-2-2) edge[red] node[right]{$s/t$} (m-4-2)
    (m-2-4) edge[red] node[right]{$s/t$} (m-4-4)
    (x2) edge[red] node[right]{$s/t$} (z2)
    (m-1-3) edge[blue] (y)
    (m-1-3) edge[blue] (y1)
    (m-1-3) edge[red,bend right=15] (m-3-3.west) 
    (y) edge[blue] (m-3-3)
    (m-3-3) edge[blue] (w)
    (m-3-3) edge[blue] (w1)
    (y.west) edge[red, bend right =15] (w.west)
    (w) edge[black] node[below] {id} (m-4-4);
    \path[->,font =\scriptsize]
    (m-1-1) edge[blue] node {$a_i$} (m-2-1)
    (m-2-1) edge[blue] node {$b_i$} (m-3-1)
    (m-3-1) edge[blue] node {$a_i$} (m-4-1)
    (m-1-1) edge[bend right =25,red] node[left] {$u_1$} (m-3-1)
    (m-2-1) edge[bend left =25,red] node[right] {$u_2$} (m-4-1);
    \node at (-2.5,0) {$\longrightarrow$};
    \draw[blue] (-2.6,0.1) -- (-2.4,0.3);
    \draw[blue] (-2.4,0.1) -- (-2.45,0.15);
    \draw[blue] (-2.55,0.25) -- (-2.6,0.3);
    \end{tikzpicture}\label{T_1_Y_2}\eea
    To read this table, the left column is the four exceptional objects, and the complex representing $\OC L_{\sbullet\sbullet}$ is on the right. The differentials (horizontal maps) are explained below
    \bea \text{\textcircled{a}}:&&L_{10} \to V\otimes V^*L_{10} \stackrel{1\otimes ev}{\to} V\otimes L_{21},~~~(\textrm{since }\hom(L_{10},L_{21})=V^*), \nn\\
    \text{\textcircled{b}}:&&V^*L_{21} \stackrel{ev}{\to} L_{31},~~~(\textrm{since }\hom(L_{21},L_{31})=V^*),\nn\\
    \text{\textcircled{c}}:&&L_{10} \stackrel{t/s\otimes u_1}{\to} \frac{t}{s} L_{31},~~~(\textrm{since $u_1$ has weight $s/t$}).\nn\eea
    We note that it is important to use the appropriate identification between $sV$ and $V^*$ (see rmk \ref{rem_BWB}) for the maps \textcircled{a}, \textcircled{b} to compose to zero. 
        
    Now for the vertical arrows of the left column, these are the irreducible maps between the objects of the exceptional collection, explained in sec.\ref{sec_CwKaB}, especially ex.\ref{ex_quiver_Y2}. For example, denoting by $\vec a\in \hom(L_{00},L_{10})=V^*$ (some arbitrary linear combination of $a_1$ and $a_2$), then $\OC(\vec a)$ gives the map of complexes on the right
    \bea L_{00} \stackrel{-\frac{t}{s}\vec a\otimes u_1\oplus \vec a}{\longrightarrow}\frac{t}{s}V^*L_{21}\oplus L_{10}.\label{vertical_map_a}\eea
    
\end{example}

\subsection{Reidemeister I Moves}\label{sec_RIM}
The Reidemeister I moves are among the simplest calculations, but they are important to get right the framing dependence of the tangle invariants. We add a cap to a crossing to make a Reidemeister I move. We have: $\CAP \OC L_{00}=tq_*\iota^*L_{11}=0$ and $L_{21}=tq_*\iota^*L_{32}=0$, but
\bea \CAP \OC L_{10}&=&tq_*\iota^*(tL_{21}\to t(1+s^{-1})L_{32}\oplus L_{21}\to \frac{t}{s}L_{42})
=t(s\to 0\oplus \frac{s}{t} \to \frac{s}{t})=st[1],\nn\\
\CAP \OC L_{31}&=&tq_*\iota^*(sL_{21}\to (1+s)L_{32}\to 0)=s^2[1].\nn\eea
Comparing this with the action of $\CAP$ alone, we conclude 
\bea \begin{tikzpicture}[scale=.6]
\draw[blue] plot[smooth, tension=.7] coordinates {(-0.4,-0.5) (0.4,0.4) (0,0.7) (-0.4,0.4) (-0.05,0.05)};
\draw[blue] plot[smooth, tension=.7] coordinates {(0.15,-0.15) (0.4,-0.5)};
\draw[blue] plot[smooth, tension=.7] coordinates {(1.8,-0.5) (1.9,0.4) (2.3,0.7) (2.7,0.4) (2.8,-0.5)};
\node at (1,0) {\scriptsize$=t[1]$};
\node at (-1,0) {\scriptsize$t^{1/4}$};
\end{tikzpicture}.\nn\eea
From the point of view of the Kauffman triangle
\bea 
\begin{tikzpicture}[scale=.6]
\draw[blue] plot[smooth, tension=.7] coordinates {(-0.4,-0.5) (0.4,0.4) (0,0.7) (-0.4,0.4) (-0.05,0.05)};
\draw[blue] plot[smooth, tension=.7] coordinates {(0.15,-0.15) (0.4,-0.5)};
\draw[blue] plot[smooth, tension=.7] coordinates {(1.8-8,-0.5) (1.9-8,0.4) (2.3-8,0.7) (2.7-8,0.4) (2.8-8,-0.5)};
\draw[blue] plot[smooth, tension=.7] coordinates {(2-5,-0.5) (2.4-5,0.1) (2.8-5,-0.5)};
\draw[blue] (2.4-5,.5) circle (0.3);
\node at (-4.7,0) {\scriptsize$[-1]$};
\node at (-0.8,0.1) {\scriptsize$t^{1/4}$};
\node at (1-5,0) {\scriptsize$\to$};
\node at (1.7-5,0.1) {\scriptsize$t^{1/2}$};
\node at (-1.8,0) {\scriptsize$\to$};
\end{tikzpicture}.\label{split_DT}\eea
The three terms are $\CAP[-1] \to t^{1/2}(t^{1/2}[1]+t^{-1/2}[-1])\CAP\to t[1]\CAP$, so in this case, the Kauffman triangle is actually split.

\begin{example}\label{ex_Reide_I}
In this example, we illustrate the same calculations, but \emph{without} using the Kauffman triangle as in the last example. We will check
\bea \begin{tikzpicture}[scale=.6]
\begin{scope}[xscale=-1]
\draw[blue] plot[smooth, tension=.7] coordinates {(-0.4,0.5) (0.4,-0.4) (0,-0.7) (-0.4,-0.4) (-0.05,-0.05)};
\draw[blue] plot[smooth, tension=.7] coordinates {(0.15,0.15) (0.4,0.5)};
\end{scope}
\draw[blue] plot[smooth, tension=.7] coordinates {(2.7,0.5) (2.8,-0.4) (3.2,-0.7) (3.6,-0.4) (3.6,0.5)};
\node at (1.7,0) {\scriptsize$=~t[1]$};
\node at (-1,0) {\scriptsize$t^{1/4}$};
\end{tikzpicture}.\nn\eea
We will consider these tangles as a map from $D^b(Y_{\sbullet})$ to $D^b(Y_{\sbullet+2})$, i.e. they are \emph{part of a bigger tangle}. Because of that, this calculation will establish the desired invariance in general.  

Let $F$ be a sheaf (resolved as a finite complex of locally free sheaves), then by definition 
\bea t^{\frac{n}{2}}\,{}^n\CUP F=L_{-[n]}\iota_*q^*F,~~~\iota=\iota_{n+1},~~q=q_{n+1}.\nn\eea
Now apply the crossing
\bea t^{\frac14+\frac{n}{2}}\,{}_n\OC\,{}^n\CUP F=\opn{cone}(L_{-[n]}\iota_*q^*F\oplus \iota_*q^*q_*\iota^*L_{-[n]}\iota_*q^*F \to \iota_*\iota^*L_{-[n]}\iota_*q^*F)[-1].\nn\eea
Note $\iota^*\iota_*q^*F={\cal O}_D(-D)q^*F[1]\oplus q^*F$, where $D=D_{n+1}$.
\bea t^{\frac14+\frac{n}{2}}\,{}_n\OC\,{}^n\CUP F=\opn{cone}(\hcancel[red]{L_{-[n]}\iota_*q^*F}\oplus \iota_*q^*q_*(\iota^*L_{-[n+1]}[1]+\hcancel[blue]{\iota^*L_{-[n]}})q^*F \to (L_{-[n+1]}[1]+\hcancel[red]{L_{-[n]}})\iota_*q^*F)[-1],\nn\eea
where the blue term vanishes by the projection formula: $q_*(\iota^*L_{-[n]})q^*F=(q_*\iota^*L_{-[n]})F$, then as $\iota^*L_{-[n]}$ has fibre degree $-1$, it has zero pushdown. In contrast, the same projection formula gives for $q_*\iota^*L_{-[n+1]}q^*F=(q_*\iota^*L_{-[n+1]})F=\opn{cone}(tL_{-[n-1]}[-1]\to t^ns^{-1}L_{[n-1]})F$.
We get
\bea t^{\frac14+\frac{n}{2}}\,{}_n\OC\,{}^n\CUP F
&=&\opn{cone}(\opn{cone}(tq^*L_{-[n-1]}[-1]\to \frac{t^n}{s}q^*L_{[n-1]})\iota_*q^*F\to L_{-[n+1]}\iota_*q^*F)\nn\\
&=&\iota_*\opn{cone}(\opn{cone}(t\iota^*L_{-[n-1]}[-1]\to \iota^*L_{-[n-1],-2,-1})q^*F\to(\iota^*L_{-[n+1]})q^*F)\nn\\
&=&\iota_*\opn{cone}(t\iota^*L_{-[n-1]}q^*F\to (\iota^*\iota_{n*}\iota^*_nL_{-[n+1]})q^*F)\nn\\
&=&\iota_*\opn{cone}(t\iota^*L_{-[n-1]}\to t(\iota^*\iota_{n*}\iota^*_nL_{-[n-1]}))q^*F\nn\\
&=&t\iota_*(\iota^*L_{-[n]})q^*F[1]=tL_{-[n]}\iota_*q^*F[1]=t(t^{\frac{n}{2}}{}^n\CUP F[1]).\nn\eea
In the calculation, we used $\spadesuit,\heartsuit,\diamondsuit$ listed in the beginning of the section.
Also in the second last step, the map $\iota^*L_{-[n-1]}\to \iota^*\iota_{n*}\iota^*_nL_{-[n-1]}$ comes out as a result of using the octahedral axiom to cancel the two red terms above, so one needs to fix this map. But the method is the same as in ex.\ref{ex_DT_map}.
We encourage the reader to take $F=\BB{C}\in D^b(Y_0)$, so as to see the manipulation more clearly.
We will leave out the proof of the other Reidemeister moves.

  Now we can add a cap and compute $\CAP \OC^{\pm n}\CUP$, that is, the tangle
  \bea \begin{tikzpicture}[scale=.6]
\draw [blue,-] (1,0) -- (0,1);
\draw [white,-,line width=2mm] (0,0) -- (1,1);
\draw [blue,-] (0,0) -- (1,1);
\draw [blue,-] (1,1.5) -- (0,2.5);
\draw [white,-,line width=2mm] (0,1.5) -- (1,2.5);
\draw [blue,-] (0,1.5) -- (1,2.5);

\centerarc[blue](0.5,2.5)(180:0:0.5);
\centerarc[blue](0.5,0)(180:360:0.5);

\node at (0.5,1.5) {$\vdots$};
\end{tikzpicture}\nn\eea
We get
\bea \CAP \OC^n\CUP=t^n[n](t^{1/2}[1]\oplus t^{-1/2}[-1]),~~~n\in\BB{Z}.\nn\eea
Taking the Poincar\'e polynomial gives
\bea \begin{tikzpicture}[scale=0.8]
    \node at (0,-0.20) {$\bigcup$};
    \node at (0,0.67) {$\bigcap$};
    \node at (0,0.25) {\scriptsize{$\boxed{n}$}};
    \node at (-0.5,.35) {$\bra$};
    \node at (0.5,.35) {$\ket$};
    \node at (3.34,0.35) {$=t^nq^{-n}(t^{1/2}q^{-1}+t^{-1/2}q),$};
    \end{tikzpicture}\nn\eea
    where $q$ keeps track of the cohomological degree (i.e. $[1]\sim q^{-1}$).
    \end{example}
    We summarise all the Reidemeister I moves here:
    \bea \begin{tikzpicture}[scale=.6]
    \draw[blue] plot[smooth, tension=.7] coordinates {(-0.4,-0.5) (0.4,0.4) (0,0.7) (-0.4,0.4) (-0.05,0.05)};
    \draw[blue] plot[smooth, tension=.7] coordinates {(0.15,-0.15) (0.4,-0.5)};
    \draw[blue] plot[smooth, tension=.7] coordinates {(1.8,-0.5) (1.9,0.4) (2.3,0.7) (2.7,0.4) (2.8,-0.5)};
    \node at (1,0) {\scriptsize$=t[1]$};
    \node at (-1,0) {\scriptsize$t^{1/4}$};
    \end{tikzpicture}~~~\begin{tikzpicture}[scale=.6]
    \draw[blue] plot[smooth, tension=.7] coordinates {(-0.4,0.5) (0.4,-0.4) (0,-0.7) (-0.4,-0.4) (-0.05,-0.05)};
    \draw[blue] plot[smooth, tension=.7] coordinates {(0.15,0.15) (0.4,0.5)};
    \draw[blue] plot[smooth, tension=.7] coordinates {(2.7,0.5) (2.8,-0.4) (3.2,-0.7) (3.6,-0.4) (3.6,0.5)};
    \node at (1.5,0) {\scriptsize$=t^{-1}[-1]$};
    \node at (-1,0) {\scriptsize$t^{-1/4}$};
\end{tikzpicture}~~~
\begin{tikzpicture}[scale=.8]
    \draw[blue] plot[smooth, tension=.7] coordinates {(0.05,0.05) (0.1,0.5)};
    \draw[blue] plot[smooth, tension=.7] coordinates {(0.1,-0.5) (0,0) (-0.3,0.3) (-0.5,0) (-0.3,-0.3) (-0.05,-0.05)};
    \node at (0.8,0) {\scriptsize$=t^{1/2}[1]$};
    \node at (-1,0) {\scriptsize$t^{-1/4}$};
    \draw[blue] (1.7,-0.5) -- (1.7,0.5);
\end{tikzpicture}~~~
\begin{tikzpicture}[scale=.8]
    \draw[blue] plot[smooth, tension=.7] coordinates {(-0.05,0.05) (-0.1,0.5)};
    \draw[blue] plot[smooth, tension=.7] coordinates {(-0.1,-0.5) (0,0) (0.3,0.3) (0.5,0) (0.3,-0.3) (0.05,-0.05)};
    \node at (1.5,0) {\scriptsize$=t^{-1/2}[-1]$};
    \node at (-.6,0) {\scriptsize$t^{-1/4}$};
    \draw[blue] (2.7,-0.5) -- (2.7,0.5);
\end{tikzpicture}\label{Reidemeister_I}\eea
We see that the normalising factors $t^{\mp1/4}$ in \eqref{braid} are necessary to make the pictures invariant under $90^{\circ}$ rotations.

\subsection{Some Tangle Invariants}
\begin{example}\label{ex_Hopf}
After the unknot there is the inescapable Hopf link. It can be quickly computed using the Kauffman DT of prop.\ref{prop_Kauffman}:
\bea \begin{tikzpicture}[scale=.6]
\draw[blue] plot[smooth, tension=.7] coordinates {(0.1,0.1) (0.5,0.5) (0.5,-1.5) (0.1,-0.1) (-0.5,0.5) (-0.5,-1.5) (-0.1,-0.1)};
\node at (0,-1.5) {\phantom{x}};
\node at (1.5,-0.4) {$\longrightarrow$};
\draw[dotted] (-1,-.4) -- (1,-0.4);
\node at (-1.3,-0.2) {\scriptsize$t^{-1/4}$};
\end{tikzpicture} 
\begin{tikzpicture}[scale=.6]
\draw[blue] plot[smooth, tension=.7] coordinates {(0.1,0.1) (0.5,0.5) (0.5,-1.5) (-0.3,-0.6) (-0.1,-0.1)};
\draw[blue] plot[smooth, tension=.7] coordinates {(0.05,-0.95) (0.3,-0.4) (-0.5,0.5) (-0.5,-1.5) (-0.1,-1.1)};
\node at (-1.3,-0.2) {\scriptsize$t^{-1/2}$};
\node at (1.3,-0.4) {$\longrightarrow$};
\draw[dotted] (-1,-.4) -- (1,-0.4);
\end{tikzpicture} 
\begin{tikzpicture}[scale=.6]
\draw[blue] plot[smooth, tension=.7] coordinates {(0.1,0.1) (0.5,0.5) (0.5,-1.5) (0,-0.9) (-0.5,-1.5) (-0.5,0.5) (-0.1,0.1) (0.3,-0.5) (0,-0.7) (-0.3,-0.5) (-0.1,-0.1)};
\node at (-1.2,-0.2) {\scriptsize$t^{-3/4}$};
\node at (1.5,-0.2) {$\stackrel{+1}{\longrightarrow}$};
\draw[dotted] (-1,-.4) -- (1,-0.4);
\end{tikzpicture} \nn\eea
The left term is the unknot with a twist which gives $t^{1/2}[1](t^{1/2}[1]\oplus t^{-1/2}[-1])$ according to \eqref{Reidemeister_I}, while the right term gives $t^{-1/2}t^{-1}[-1](t^{1/2}[1]\oplus t^{-1/2}[-1])$. Due to degree reasons, the middle term can only be the direct sum of the two side terms
\bea \Phi({\rm Hopf})=(t^{1/2}[1]\oplus t^{-3/2}[-1])(t^{1/2}[1]\oplus t^{-1/2}[-1]),\nn\eea
with Poincar\'e polynomial 
\bea t^{-1/2}(t^{3/2}q^{-2}+t^{1/2}+t^{-1/2}+t^{-3/2}q^2).\nn\eea
Setting $q=-1$ and correcting for the writhe, this gives the Jones polynomial. 

We do another exercise to understand the maps in the above DT better, and to give more examples of working with the line bundles $L_{\vec a}$. We start with the empty link and compute upward to the dotted line, which would involve $Y_4$ at intermediate stages. First set $F={}^3\CUP {}^1\CUP~\BB{C}=t^{-2}L_{-2,-1,-1,0}|_{D_2\cap D_4}$. 
From now on we use the notation in \eqref{L_with_line}, using an underline for the restriction to a divisor, so $L_{-2,\udl{-1},-1,\udl{0}}:=L_{-2,-1,-1,0}|_{D_2\cap D_4}$.

For the twist, we have the DT
\bea &F\to {}_2\UC F\to t^{-\frac12}{}^2\CUP\,{}_2\CAP F\stackrel{+1}{\to}.\nn\eea
We can compute 
\bea {}^2\CUP\,{}_2\CAP F=t^{1/2}t^{-2}\opn{cone}(t^{-1}L_{-2,-2,-1,-1}|_{D_3}\to L_{-2,-1,0,0}|_{D_3})
={}^2\CUP\,{}^1\CUP\BB{C},\nn\eea
and then look at what the connecting map can be. The relevant cohomology is 
\bea \hom^0(t^{-\frac12}{}^2\CUP\,{}_2\CAP F,F[1])=\ext^1(\opn{cone}(t^{-1}L_{-2,-2,\udl{-1},-1}\to L_{-2,-1,\udl{0},0}),L_{-2,\udl{-1},-1,\udl{0}})=1.\nn\eea
Therefore the connecting map is a (possibly zero) multiple of a generator of the $\ext^1$ above. An argument similar to that of ex.\ref{ex_DT_map} shows that it is nonzero.

This map is induced from the restriction $\hom(L_{-2,-1,-1,0}\to L_{-2,\udl{-1},-1,\udl{0}})$. We keep applying the functors above the dotted line and show that this map gets nullified, agreeing with what we argued above that the DT splits. Acting with ${}_2\UC$ on $L_{-2,-1,-1,0}$ we get
\bea {}_2\UC L_{-2,-1,-1,0}&=&L_{-1,-1,0,0}\otimes\opn{cone}(L_{-1,0,0,0}\oplus \iota_{3*}q_3^*q_{3*}\iota_3^*L_{-1,0,0,0} \to L_{-1,0,\udl{0},0})[-1]\nn\\
&=&L_{-1,-1,0,0}\otimes L_{-1,0,0,0}=L_{-2,-1,0,0}, \nn\eea
which is annihilated by both ${}_1\CAP$ and ${}_3\CAP$, hence the agreement. 
\end{example}

\begin{example}
We also compute a tangle with open ends, with presentation $K_{{\rm fig\,8}}={}_1\CAP\;{}_2\UC\;{}_1\OC\;{}_2\UC\;{}^1\CUP$. It is depicted in fig.\ref{fig_8_stopper}.
\begin{figure}
\begin{center}
\begin{tikzpicture}[scale=.6]
\draw [blue,-] (0.8,0.2) -- (0,1);
\draw [white,-,line width=2mm] (0,0) -- (1,1);
\draw [blue,-] (0,0) -- (1,1);
\draw [blue,-] (0.5,-0.5) -- (1.5,.5);
\draw [white,-,line width=2mm] (0.8,0.2) -- (1.2,-0.2);
\draw [blue,-] (0.8,0.2) -- (1.4,-0.4);
\draw [blue,-] (0,0) -- (0,-0.5);
\draw[blue] plot[smooth, tension=.7] coordinates {(0,-0.5) (0.1,-0.9) (0.5,-0.5)};
\draw[blue] plot[smooth, tension=.7] coordinates {(1,1) (1.5,1.5) (1.5,2.2)};
\draw [white,-,line width=2mm] (1.2,0.9) -- (0.8,1.3);
\draw[blue] plot[smooth, tension=.7] coordinates {(0,1) (0,1.6) (0.3,1.9) (1.5,0.5)};
\draw[blue] plot[smooth, tension=.7] coordinates {(1.4,-0.4) (1.5,-0.6) (1.5,-1.2)};
\end{tikzpicture}  
\caption{The figure 8 stopper. If one closes the tangle, one gets the right-handed trefoil.}\label{fig_8_stopper}
\end{center}
\end{figure}
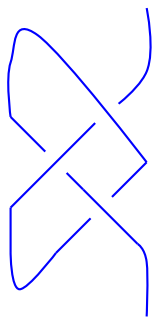

By computing this with the Kauffman DT, we expect that the answer will be proportional to the identity morphism. 
Indeed, we can always reduce the tangle into a straight line with some disconnected unknots, each giving $t^{1/2}[1]+t^{-1/2}[-1]$. From prop.\ref{prop_Kauffman}, we get
\bea 
\begin{tikzpicture}[scale=.6]
\draw [blue,-] (0.9,0.1) -- (0,1);
\draw [white,-,line width=2mm] (0,0) -- (1,1);
\draw [blue,-] (0,0) -- (1,1);
\draw [blue,-] (0,0) -- (0,-0.5);
\draw[blue] plot[smooth, tension=.7] coordinates {(0,-0.5) (0.1,-0.9) (0.5,-0.5) (0.9,0) (0.9,0.1)};
\draw[blue] plot[smooth, tension=.7] coordinates {(1,1) (1.5,1.5) (1.5,2.2)};
\draw [white,-,line width=2mm] (1.2,0.9) -- (0.8,1.3);
\draw[blue] plot[smooth, tension=.7] coordinates {(0,1) (0,1.6) (0.3,1.9) (1.5,0.5)};
\draw[blue] plot[smooth, tension=.7] coordinates {(1.5,0.5) (1.5,0.4) (1.2,0) (1.5,-0.6) (1.5,-1.2)};
\node at (2.2,0.4) {$\longrightarrow$};
\end{tikzpicture}
\begin{tikzpicture}[scale=.6]
\draw [blue,-] (0.8,0.2) -- (0,1);
\draw [white,-,line width=2mm] (0,0) -- (1,1);
\draw [blue,-] (0,0) -- (1,1);
\draw [blue,-] (0.5,-0.5) -- (1.5,.5);
\draw [white,-,line width=2mm] (0.8,0.2) -- (1.2,-0.2);
\draw [blue,-] (0.8,0.2) -- (1.4,-0.4);
\draw [blue,-] (0,0) -- (0,-0.5);
\draw[blue] plot[smooth, tension=.7] coordinates {(0,-0.5) (0.1,-0.9) (0.5,-0.5)};
\draw[blue] plot[smooth, tension=.7] coordinates {(1,1) (1.5,1.5) (1.5,2.2)};
\draw [white,-,line width=2mm] (1.2,0.9) -- (0.8,1.3);
\draw[blue] plot[smooth, tension=.7] coordinates {(0,1) (0,1.6) (0.3,1.9) (1.5,0.5)};
\draw[blue] plot[smooth, tension=.7] coordinates {(1.4,-0.4) (1.5,-0.6) (1.5,-1.2)};
\node at (-1,0.5) {$t^{-1/4}$};
\node at (2.2,0.4) {$\longrightarrow$};
\end{tikzpicture}
\begin{tikzpicture}[scale=.6]
\draw [blue,-] (0.9,0.1) -- (0,1);
\draw [white,-,line width=2mm] (0,0) -- (1,1);
\draw [blue,-] (0,0) -- (1,1);
\draw [blue,-] (0,0) -- (0,-0.5);
\draw[blue] plot[smooth, tension=.7] coordinates {(0,-0.5) (0.1,-0.9) (0.8,-0.2) (1.1,-0.1) (1.4,-0.2) (1.5,-0.5) (1.5,-1.2)};
\draw[blue] plot[smooth, tension=.7] coordinates {(1,1) (1.5,1.5) (1.5,2.2)};
\draw [white,-,line width=2mm] (1.2,0.9) -- (0.8,1.3);
\draw[blue] plot[smooth, tension=.7] coordinates {(0,1) (0,1.6) (0.3,1.9) (1.5,0.5)};
\draw[blue] plot[smooth, tension=.7] coordinates {(1.5,0.5) (1.5,0.4) (1.2,0.1) (0.9,0.1)};
\node at (-1,0.5) {\small$t^{-1/2}$};
\node at (2.2,0.6) {$\stackrel{+1}{\longrightarrow}$};
\end{tikzpicture} \nn\eea
The left term is just the straight line with two Reidemeister I moves \eqref{Reidemeister_I}, which we have worked out in ex.\ref{ex_Reide_I}. The simplification of the right term is similar to ex.\ref{ex_Hopf}. We get
\bea 
\begin{tikzpicture}[scale=.6]
\draw [blue,-] (0.3,1.5) -- (0.3,-1.5);
\node at (-0.9,0.2) {\small$t^{3/2}[2]$};
\node at (1.1,.1) {$~\longrightarrow$};
\end{tikzpicture}
\begin{tikzpicture}[scale=.6]
\draw [blue,-] (0.8,0.2) -- (0,1);
\draw [white,-,line width=2mm] (0,0) -- (1,1);
\draw [blue,-] (0,0) -- (1,1);
\draw [blue,-] (0.5,-0.5) -- (1.5,.5);
\draw [white,-,line width=2mm] (0.8,0.2) -- (1.2,-0.2);
\draw [blue,-] (0.8,0.2) -- (1.4,-0.4);
\draw [blue,-] (0,0) -- (0,-0.5);
\draw[blue] plot[smooth, tension=.7] coordinates {(0,-0.5) (0.1,-0.9) (0.5,-0.5)};
\draw[blue] plot[smooth, tension=.7] coordinates {(1,1) (1.5,1.5) (1.5,2.2)};
\draw [white,-,line width=2mm] (1.2,0.9) -- (0.8,1.3);
\draw[blue] plot[smooth, tension=.7] coordinates {(0,1) (0,1.6) (0.3,1.9) (1.5,0.5)};
\draw[blue] plot[smooth, tension=.7] coordinates {(1.4,-0.4) (1.5,-0.6) (1.5,-1.2)};
\node at (-1,0.5) {$t^{-1/4}$};
\node at (2.4,0.4) {$\longrightarrow$};
\end{tikzpicture}
\begin{tikzpicture}[scale=.6]
\draw [blue,-] (5.6,1.5) -- (5.6,-1.5);
\node at (3.5,.2) {\small$t^{-1/2}(t^{-1}[-1]+t^1[1])\quad \stackrel{+1}{\longrightarrow}$};
\end{tikzpicture}\nn\eea
Here again we need to understand the maps in between. But the straight lines correspond to the kernel ${\cal O}_{\Gd}\in D^b(Y_1\times Y_1)$. It will turn out that the latter has a resolution in terms of the exceptional objects concentrated at degrees $-1,0$. So the middle term can only be a direct sum. We conclude that
\bea \Phi(K_{\rm fig\,8})=t^{\frac14}(t^{\frac32}[2]\oplus t^{-\frac32}[-1]\oplus t^{\frac12}[1]){\cal O}_{\Gd}.\nn\eea
\end{example}
\begin{example}
 Another commonly used knot in sailing is the bowline knot fig.\ref{fig_palstek}. It has presentation \bea {}_1\CAP\;{}_2\OC\;{}_1\UC\;{}_3\CAP\;{}_2\UC\;{}_4\OC\;{}_1\OC\;{}_2\UC\;{}_4\OC\;{}^3\CUP\;{}^1\CUP\label{palstek}\eea
\begin{figure}[h]
\begin{center}
\includegraphics[width=3cm]{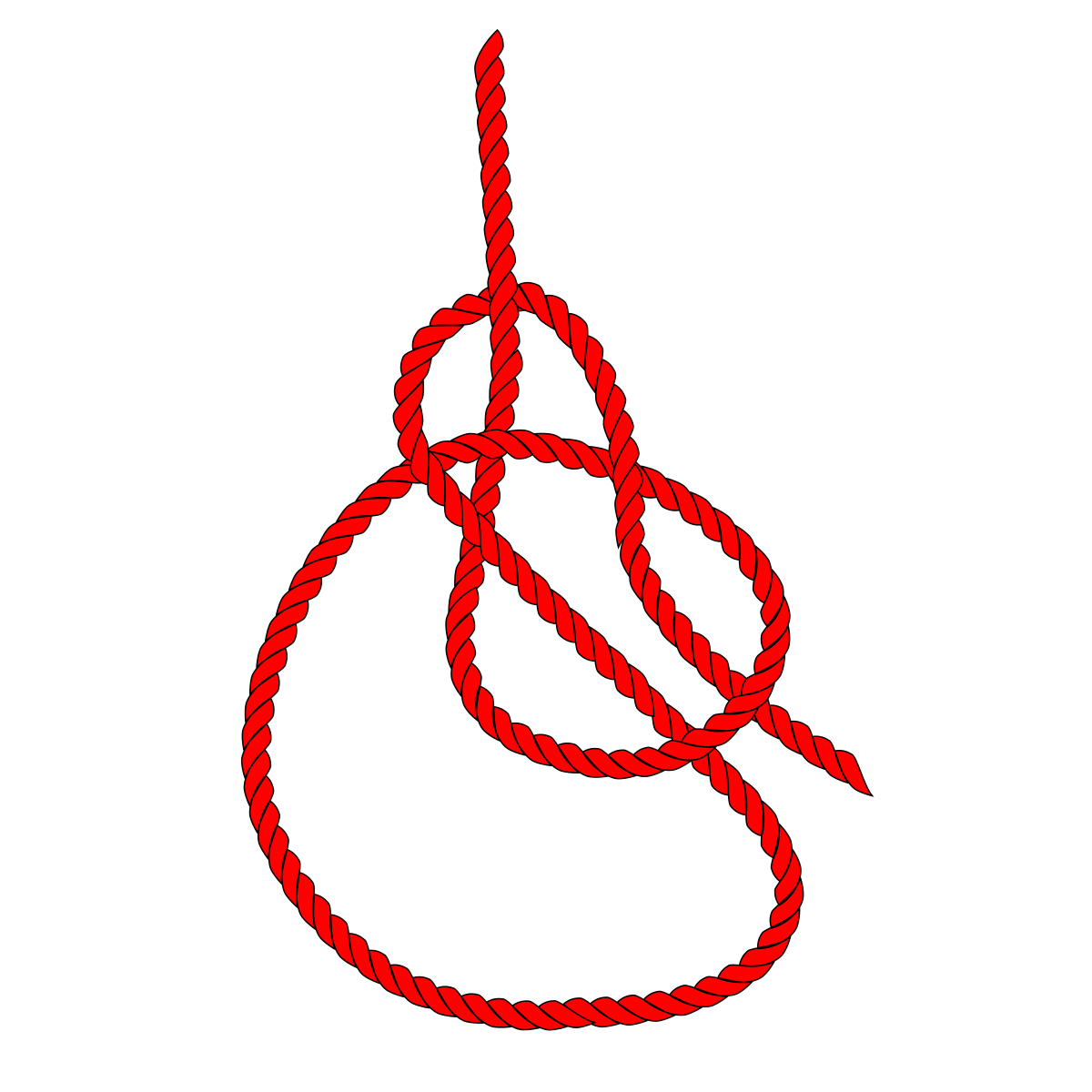}~~~
\begin{tikzpicture}[scale=.8]
\draw[blue] (0.6,1.6) -- (0,1.1)-- (0,0) -- (0.5,0.5);
\draw[blue] (1.6,1.5) -- (2,1.1) -- (1,0.1) -- (0.2,1.1); \draw[blue] (.11,1.25) -- (0,1.4) -- (0,1.9) -- (0.15,2.05);
\draw[blue] (0.65,0.65) -- (1,1) -- (1,1.4) -- (0,2.3);
\draw[blue] (1.5,0) -- (1.5,0.5); 
\draw[blue] (1.5,0.68) -- (1.5,2.5) -- (1.2,2.75); \draw[blue] (1.1,2.85) -- (.9,3.0) -- (0,2.3);
\draw[blue] (1.4,1.72) -- (1.05,2.05) -- (0.72,1.72); 
\draw[blue] (0.3,2.15) -- (1.5,3.05);
\end{tikzpicture}\caption{The bowline or p{\aa}lstek knot.}\label{fig_palstek}
\end{center}
\end{figure}
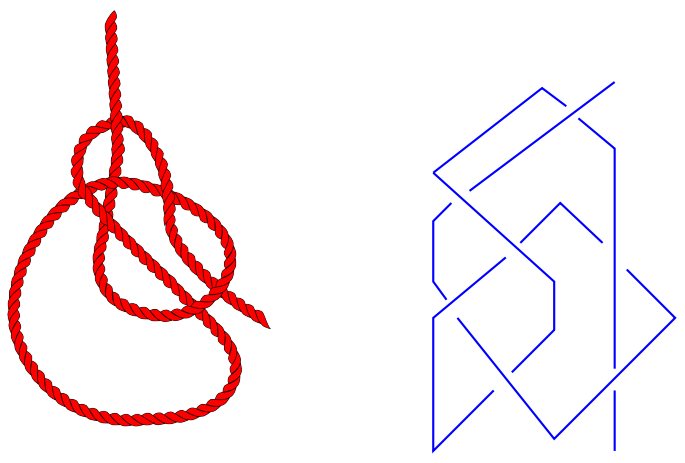
This is still a 1 strand to 1 strand tangle, so is `proportional' to the identity map $D^b(Y_1)\to D^b(Y_1)$, and its kernel is proportional to ${\cal O}_{\Gd}$. 
In this calculation, we did not use the Kauffman triangles but simply stacked up the maps in \eqref{palstek}, and applied them to the two exceptional objections $L_0,L_1$. We obtain
\bea \Phi(\textrm{bowline})=t^{\frac34}(t^2[3]+2t[2]+2[1]+2t^{-1}[0]+2t^{-2}[-1]+t^{-3}[-2]+t^{-4}[-3]){\cal O}_{\Gd}.\nn\eea
\end{example}

\begin{example}
We consider a tangle with two open ends
\begin{center}
\begin{tikzpicture}[scale=1]
\draw[blue] plot[smooth, tension=.7] coordinates {(0,0) (0.40,0.75) (0.6,0.9) (0.8,0.7)};
\draw[blue] plot[smooth, tension=.7] coordinates {(0.88,0.5) (1.25,0)};
\draw[blue] plot[smooth, tension=.7] coordinates {(0,1.4) (0.35,0.8)};
\draw[blue] plot[smooth, tension=.7] coordinates {(0.42,0.6) (0.6,0.43) (0.82,0.55) (1.25,1.4)};
\node[black] at (-.5,0.7) {\scriptsize$B=$};
\end{tikzpicture}
\end{center}
This gives an autofunctor of $D^b(Y_2)$, so as in example \ref{ex_DT_map} we can understand it by testing it on the exceptional collection of that category. Applying the Kauffman DT we get a DT of functors or kernels
\bea 
\begin{tikzpicture}[scale=1]
\draw[blue] plot[smooth, tension=.7] coordinates {(0,1.2) (0.4,.8) (0.6,0.95) (0.8,0.65)};
\draw[blue] plot[smooth, tension=.7] coordinates {(0.88,0.5) (1.25,0)};
\draw[blue] plot[smooth, tension=.7] coordinates {(0,0) (0.4,.65) (0.82,0.55) (1.25,1.4)};
\node[blue] at (1.6,0.65) {$\to$};
\node[blue] at (-0.3,0.65) {\scriptsize$t^{3/4}$};
\end{tikzpicture}~\begin{tikzpicture}[scale=1]
\draw[blue] plot[smooth, tension=.7] coordinates {(0,0) (0.40,0.75) (0.6,0.9) (0.8,0.7)};
\draw[blue] plot[smooth, tension=.7] coordinates {(0.88,0.5) (1.25,0)};
\draw[blue] plot[smooth, tension=.7] coordinates {(0,1.4) (0.35,0.8)};
\draw[blue] plot[smooth, tension=.7] coordinates {(0.42,0.6) (0.6,0.43) (0.82,0.55) (1.25,1.4)};
\node[blue] at (-0.3,0.65) {\scriptsize$t^{1/2}$};
\node[blue] at (1.6,0.65) {$\to$};
\end{tikzpicture}~
\begin{tikzpicture}[scale=1]
\draw[blue] plot[smooth, tension=.7] coordinates {(0,0) (0.3,0.6) (0,1.4)};
\draw[blue] plot[smooth, tension=.7] coordinates {(0.8,0.65) (0.6,0.85) (0.4,0.5) (0.6,0.3) (0.82,0.56)};
\draw[blue] plot[smooth, tension=.7] coordinates {(0.88,0.5) (1.25,0)};
\draw[blue] plot[smooth, tension=.7] coordinates {(0.82,0.55) (1.25,1.4)};
\node[blue] at (-0.3,0.65) {\scriptsize$t^{1/4}$};
\end{tikzpicture}.\nn\eea
The left term is just the overcrossing $t^{3/4}\OC$, while the right one is $t^{-1/2}[-1]$. 
We rotate the rightmost term to the left, giving $t^{-1/2}[-2]$, so the braid $B$ induces the map
\bea B=\opn{cone}(t^{-1}[-2]\to t^{1/4}\OC).\nn\eea
In particular, if we test this map on the exceptional objects $\{L_{00},L_{10},L_{21},L_{31}\}$, then $\OC L_{\sbullet\sbullet}$ only has terms  
at degrees $-1,0,1$ (see \eqref{T_1_Y_2}), thus the $-2$ shift excludes any possible maps. Note this argument is only applicable when the objects are resolved using the strong exceptional collection, as maps within the collection are always of degree 0. We conclude that
\bea BF=t^{1/4}\OC F\oplus t^{-1}F[-1],~~~F=L_{\sbullet\sbullet}.\label{consistent_with}\eea

As a verification, we look into how $B$ acts on $L_{00}$, without resorting to the Kauffman triangles. Using the notation from \eqref{L_with_line}, we have:
\bea L_{00}&\stackrel{{}^2\CUP}{\to}& t^{-1}L_{-1,-1,\udl{0},0}\nn\\
&\stackrel{t^{1/4}{}_1\OC}{\to}&t^{-1}\opn{cone}
(\iota_{2*}q_2^*q_{2*}\iota_2^*L_{-1,-1,\udl{0},0}[-1] \to L_{-1,-2,\udl{0},0})\nn\\
&=&t^{-1}\opn{cone}
(\frac{t^2}{s}L_{1,\udl{1},1,0}[-1]\to \frac{t^2}{s}L_{1,0,\udl{1},0})
=\frac{t}{s}\opn{cone}(L_{1000}\to L_{1110})\nn\\
&=&\frac{t}{s}\opn{cone}(L_{11\udl{1}0}[-1]\to L_{1\udl{1}00})
.\nn\eea
We compute ${}_3\OC$. It acts trivially on $L_{1\udl{1}00}$, but on $L_{11\udl{1}0}$
\bea t^{\frac14}{}_3\OC L_{11\udl{1}0}&=&\opn{cone}(\iota_{4*}q_4^*q_{4*}\iota_4^*L_{11\udl{1}0}[-1]\to L_{11\udl{1},-1})
=\opn{cone}(\frac{s}{t^2}L_{-1,-1,0,\udl{0}}[-1]\to L_{11\udl{1},-1})\nn\\
{}_2\CAP\, t^{\frac14}{}_3\OC L_{11\udl{1}0}&=&t^{3/2}\opn{cone}(\frac{s^2}{t^4}L_{-2,-2}[-1]\to \opn{cone}(\frac{s}{t^3}L_{00}[-2]\to \frac{s^2}{t^4}L_{-2,-2}[-1]))
=t^{3/2}\frac{s}{t^3}L_{00}[-1].\nn\eea
Putting it all together
\bea {}_2\CAP\,(t^{\frac14}{}_3\OC)(t^{\frac14}{}_1\OC)\,{}^2\CUP L_{00}=\frac{t}{s}t^{3/2}\opn{cone}(\frac{s}{t^3}L_{00}[-2]\to \frac{s}{t^2}L_{00})
=t^{3/2}(\frac{1}{t^2}L_{00}[-1]\oplus \frac{1}{t}L_{00}).\nn\eea
We conclude that the whole tangle acts on $L_{00}$ as
\bea BL_{00}=(t^{-1}[-1]+1)L_{00}.\nn\eea
This is consistent with the result \eqref{consistent_with}, as we have shown already that $t^{1/4}\OC L_{00}=L_{00}$. We can check the action of $B$ on $L_{10},L_{21},L_{31}$ the same way.
\end{example}

\subsection{Proofs of some Isotopy Invariance Statements}\label{sec_PosIIS}
That the tangle invariant constructed above satisfies the requisite isotopy invariance was already proved in \cite{Cautis2007KnotHV}, sec.5 (see their lem.4.1 for the list of moves to be checked). Our proofs are slightly different, and rely almost exclusively on the Kauffman triangle. We will illustrate a few cases.
\begin{lemma}
  We have the pitchfork moves
  \bea t^{\frac14}{}_n\CAP\;{}_{n-1}\OC=t^{-\frac14}{}_{n-1}\CAP\;{}_n\UC,~~~t^{-\frac14}{}_n\CAP\;{}_{n-1}\UC=t^{\frac14}{}_{n-1}\CAP\;{}_n\OC\nn\\
  t^{\frac14}{}_{n-1}\OC\;{}^n\CUP=t^{-\frac14}{}_n\UC\;{}_{n-1}\CUP,~~~t^{-\frac14}{}_{n-1}\UC\;{}^n\CUP=t^{\frac14}{}_n\OC\;{}^{n-1}\CUP.\nn\eea
\end{lemma}
\begin{proof}
We prove the first one. In pictures, it looks like
\begin{center}
\begin{tikzpicture}[scale=.8]
\draw[blue] plot[smooth, tension=.7] coordinates {(0,0) (0.3,.8) (0.5,1) (0.7,0.8) (1,0)};
\draw[blue] plot[smooth, tension=.7] coordinates {(0.5,0) (0.25,.4)};
\draw[blue] plot[smooth, tension=.7] coordinates {(0.15,0.6) (0,1)};
\node at (-0.6,0.5) {\scriptsize$t^{1/4}$};
\node at (1.5,0.5) {\scriptsize$=$};
\end{tikzpicture}
\begin{tikzpicture}[scale=.8]
\draw[blue] plot[smooth, tension=.7] coordinates {(0,0) (0.3,.8) (0.5,1) (0.7,0.8) (1,0)};
\draw[blue] plot[smooth, tension=.7] coordinates {(0.5,0) (0.75,.4)};
\draw[blue] plot[smooth, tension=.7] coordinates {(0.85,0.6) (1,1)};
\node at (-0.6,0.55) {\scriptsize$t^{-1/4}$};
\end{tikzpicture}
\end{center}
Before jumping into the details, we attempt a graphic proof using the Kauffman DT:
\begin{center}
\begin{tikzpicture}[scale=.8]
\draw[blue] plot[smooth, tension=.7] coordinates {(0.5,0) (0.25,.4) (0.3,.8) (0.5,1) (0.7,0.8) (1,0)};
\draw[blue] plot[smooth, tension=.7] coordinates {(0,0) (0.15,0.6) (0,1)};
\node at (1.3,0.5) {\scriptsize$[-1]$};
\node at (2.2,0.5) {\scriptsize$\to$};
\begin{scope}[xshift=3.5cm] 
\draw[blue] plot[smooth, tension=.7] coordinates {(0,0) (0.15,.4) (0.25,.4)};
\draw[blue] plot[smooth, tension=.7] coordinates {(0.5,0) (0.25,.4)};
\draw[blue] plot[smooth, tension=.7] coordinates {(0,1) (0.15,0.6) (0.25,0.6) (0.5,1) (0.7,0.8) (1,0)};
\node at (-0.5,0.5) {\scriptsize$t^{1/2}$};
\end{scope}
\begin{scope}[yshift=-1.5cm]
  \draw[blue] plot[smooth, tension=.7] coordinates {(1,0) (0.85,.4) (0.75,.4)};
\draw[blue] plot[smooth, tension=.7] coordinates {(0.5,0) (0.75,.4)};
\draw[blue] plot[smooth, tension=.7] coordinates {(1,1) (0.85,0.6) (0.75,0.6) (0.5,1) (0.3,0.8) (0,0)};
\node at (1.3,0.5) {\scriptsize$[-1]$};
\node at (2.2,0.5) {\scriptsize$\to$};
\node at (.5,1.2) {\scriptsize\rotatebox{90}{$\simeq$}};
\end{scope}
\begin{scope}[xshift=3.5cm,yshift=-1.5cm]
\draw[blue] plot[smooth, tension=.7] coordinates {(0.5,0) (0.75,.4) (0.7,.8) (0.5,1) (0.3,0.8) (0,0)};
\draw[blue] plot[smooth, tension=.7] coordinates {(1,0) (0.85,0.6) (1,1)};
\node at (-0.5,0.5) {\scriptsize$t^{1/2}$};
\node at (.5,1.2) {\scriptsize\rotatebox{90}{$\simeq$}};
\end{scope}
\end{tikzpicture}
\end{center}
The cone of the top row is the left pitchfork. We use prop.\ref{prop_cap_cup} to reach the lower row. If we could prove that the induced map of the lower row is part of the Kauffman DT, then its cone would give the right pitchfork. This is slightly tricky as maps in the Kauffman DT's only exist in the derived category, and are not easy to describe. So we return to the formula \eqref{braid} where all maps were still canonical.

Add a cap on either sides to \eqref{braid}
\bea &&t^{\frac14-\frac{n+1}{2}}{}_n\!\CAP\; {}_{n-1}\OC F[1]=q_{n+1,*}\iota_{n+1}^*L_{[n+1]}\opn{cone}(F\oplus \iota_{n*}q_n^*q_{n*}\iota_n^*F\to\iota_{n*}\iota_n^*F)\nn\\
&&t^{-\frac14-\frac{n}{2}}{}_{n-1}\CAP\;{}_n\UC F[1]=q_{n*}\iota_n^*\opn{cone}(F_{[n+1]}\oplus \iota_{n+1,*}q_{n+1}^*q_{n+1,*}\iota_{n+1}^*F_{[n+1]}\to\iota_{n+1,*}\iota_{n+1}^*F_{[n+1]}).\nn\eea
We simplify each term.
For the top middle term
\bea &&q_{n+1,*}\iota_{n+1}^*L_{[n+1]}\iota_{n*}q_n^*q_{n*}\iota_n^*F
  \stackrel{1}{=}q_{n+1,*}\iota_{n+1}^*\iota_{n*}(\iota_n^*L_{[n+1]})q_n^*q_{n*}\iota_n^*F
  \stackrel{2}{=}t^{-1}q_{n+1,*}\iota_{n+1}^*\iota_{n*}q_n^*L_{[n-1]}q_{n*}\iota_n^*F\nn\\
  &\stackrel{3}{=}&q_{n+1,*}\iota_{n+1}^*\iota_{n*}q_n^*q_{n*}(L_{[n+1]}\iota_n^*F)
  \stackrel{4}{=}q_{n+1,*}\iota_{n+1}^*\iota_{n*}q_n^*q_{n*}\iota_n^*F_{[n+1]},\nn\eea
  where we used $\iota_n^*L_{[n+1]}=t^{-1}q_n^*L_{[n-1]}$ (see lem.\ref{prop_pushdown_q}) in steps 2,4 and projection formula in steps 1, 3. 
  Now recall from lem.\ref{lem_trivial_lb} that the two projections $Y_{\sbullet-2}\stackrel{q_{n+1}}{\ot}D_n\cap D_{n+1}\stackrel{q_n}{\to} Y_{\sbullet-2}$ are both isomorphisms, so $q_{n+1,*}\iota_{n+1}^*\iota_{n*}q_n^*=1$ and the rhs is just $q_{n*}\iota_n^*F_{[n+1]}$. 
  The simplification of other terms are similar if not easier, we get
  \bea &&t^{\frac14-\frac{n+1}{2}}{}_n\!\CAP\; {}_{n-1}\OC F[1]=\opn{cone}(q_{n+1,*}\iota_{n+1}^*F_{[n+1]}\oplus q_{n*}\iota_n^*F_{[n+1]}\to q_*\iota^*F_{[n+1]})\nn\\
&&t^{-\frac14-\frac{n}{2}}{}_{n-1}\CAP\;{}_n\UC F[1]=\opn{cone}(q_{n*}\iota_n^*F_{[n+1]}\oplus q_{n+1,*}\iota_{n+1}^*F_{[n+1]}\to q_*\iota^*F_{[n+1]}),\nn\eea
  where in the rightmost terms, $\iota$ is the inclusion of $D_n\cap D_{n+1}$ into $Y_{\sbullet}$, and $q$ is the projection from $D_n\cap D_{n+1}$ to $Y_{\sbullet-2}$ (using either $q_n$ or $q_{n+1}$). Noting that the horizontal maps are just restrictions from $D_n$ or $D_{n+1}$ to $D_n\cap D_{n+1}$, it is clear that the two cones are isomorphic.
\end{proof}
Next we prove a key isotopy move, after which the other proofs are mere repetitions of the same.
\begin{lemma}(Reidemeister III)\label{lem_Rd_III}
  We have ${}_{n-1}\OC\;{}_n\OC\;{}_{n-1}\OC={}_n\OC\;{}_{n-1}\OC\;{}_n\OC$, or in picture
  \begin{center}
\begin{tikzpicture}
\begin{scope}[scale=.4]
\draw[blue] plot[smooth, tension=.7] coordinates {(0,0) (2,3)};
\draw[blue] plot[smooth, tension=.7] coordinates {(1,0) (.6,0.6)};
\draw[blue] plot[smooth, tension=.7] coordinates {(0.4,0.9) (0.1,1.65) (1,3)};
\draw[blue] plot[smooth, tension=.7] coordinates {(2,0) (1.1,1.35)};
\draw[blue] plot[smooth, tension=.7] coordinates {(0.9,1.65) (0.6,2.10)};
\draw[blue] plot[smooth, tension=.7] coordinates {(0.4,2.4) (0,3)};
\end{scope}
\node at (1.4,0.6) {$=$};
\begin{scope}[xshift=2cm,scale=.4]
\draw[blue] plot[smooth, tension=.7] coordinates {(0,0) (2,3)};
\draw[blue] plot[smooth, tension=.7] coordinates {(1,0) (1.9,1.3) (1.6,2.1)};
\draw[blue] plot[smooth, tension=.7] coordinates {(1.4,2.4) (1,3)};
\draw[blue] plot[smooth, tension=.7] coordinates {(2,0) (1.6,0.6)};
\draw[blue] plot[smooth, tension=.7] coordinates {(1.4,0.9) (1.1,1.35)};
\draw[blue] plot[smooth, tension=.7] coordinates {(0.9,1.65) (0,3)};
\end{scope}
\end{tikzpicture}
\end{center}
and same for the undercrossing.
\end{lemma}
\begin{proof}
  Start from the lhs, apply the Kauffman DT to the top crossing, we get that it is the cone of
  \begin{center}
    \begin{tikzpicture}[scale=.45]
\draw[blue] plot[smooth, tension=.7] coordinates {(0,0) (2,3)};
\draw[blue] plot[smooth, tension=.7] coordinates {(2,0) (0,3)};
\draw[blue] plot[smooth, tension=.2] coordinates {(1,0) (0.1,1.35) (0.1,1.65) (1,3)};
\begin{scope}[xshift=0.4cm,yshift=2.11cm]%
  \filldraw[white](0,0) rectangle (0.2,0.3);
\draw[blue] plot[smooth, tension=.7] coordinates {(0,0) (0.06,0.15) (0,0.3)};
\draw[blue] plot[smooth, tension=.7] coordinates {(0.2,0) (0.14,0.15) (0.2,0.3)};
\end{scope}
\begin{scope}[xshift=0.9cm,yshift=1.35cm]%
  \filldraw[white](0,0) rectangle (0.2,0.3);
\draw[blue] plot[smooth, tension=.7] coordinates {(0,0) (0.2,0.3)};
\draw[blue] plot[smooth, tension=.7] coordinates {(0.2,0) (0.15,0.075)};
\draw[blue] plot[smooth, tension=.7] coordinates {(0.05,0.225) (0,0.3)};
\end{scope}
\begin{scope}[xshift=0.4cm,yshift=0.6cm]%
  \filldraw[white](0,0) rectangle (0.2,0.3);
\draw[blue] plot[smooth, tension=.7] coordinates {(0,0) (0.2,0.3)};
\draw[blue] plot[smooth, tension=.7] coordinates {(0.2,0) (0.15,0.075)};
\draw[blue] plot[smooth, tension=.7] coordinates {(0.05,0.225) (0,0.3)};
\end{scope}
\node at (1.7,1.5) {\scriptsize$[-1]$};
\draw[blue,->] (2.2,1.5) -- (2.8,1.5);
\begin{scope}[xshift=3cm]
\draw[blue] plot[smooth, tension=.7] coordinates {(0,0) (2,3)};
\draw[blue] plot[smooth, tension=.7] coordinates {(2,0) (0,3)};
\draw[blue] plot[smooth, tension=.2] coordinates {(1,0) (0.1,1.35) (0.1,1.65) (1,3)};
\end{scope}
\begin{scope}[xshift=3.4cm,yshift=2.11cm]%
  \filldraw[white](0,0) rectangle (0.2,0.3);
\draw[blue] plot[smooth, tension=.7] coordinates {(0,0) (0.1,0.1) (0.2,0)};
\draw[blue] plot[smooth, tension=.7] coordinates {(0,0.3) (0.1,0.2) (0.2,0.3)};
\end{scope}
\begin{scope}[xshift=3.9cm,yshift=1.35cm]%
  \filldraw[white](0,0) rectangle (0.2,0.3);
\draw[blue] plot[smooth, tension=.7] coordinates {(0,0) (0.2,0.3)};
\draw[blue] plot[smooth, tension=.7] coordinates {(0.2,0) (0.15,0.075)};
\draw[blue] plot[smooth, tension=.7] coordinates {(0.05,0.225) (0,0.3)};
\end{scope}
\begin{scope}[xshift=3.4cm,yshift=0.6cm]%
  \filldraw[white](0,0) rectangle (0.2,0.3);
\draw[blue] plot[smooth, tension=.7] coordinates {(0,0) (0.2,0.3)};
\draw[blue] plot[smooth, tension=.7] coordinates {(0.2,0) (0.15,0.075)};
\draw[blue] plot[smooth, tension=.7] coordinates {(0.05,0.225) (0,0.3)};
\end{scope}
\node at (4.7,1.6) {\scriptsize$t^{1/2}$};
\end{tikzpicture}
\end{center}
Apply DT to the middle crossing, we get the cone of cone of the rows of
\begin{center}
\begin{tikzpicture}[scale=.45]
\draw[blue] plot[smooth, tension=.7] coordinates {(0,0) (2,3)};
\draw[blue] plot[smooth, tension=.7] coordinates {(2,0) (0,3)};
\draw[blue] plot[smooth, tension=.2] coordinates {(1,0) (0.1,1.35) (0.1,1.65) (1,3)};
\begin{scope}[xshift=0.4cm,yshift=2.11cm]%
  \filldraw[white](0,0) rectangle (0.2,0.3);
\draw[blue] plot[smooth, tension=.7] coordinates {(0,0) (0.06,0.15) (0,0.3)};
\draw[blue] plot[smooth, tension=.7] coordinates {(0.2,0) (0.14,0.15) (0.2,0.3)};
\end{scope}
\begin{scope}[xshift=0.9cm,yshift=1.35cm]%
\filldraw[white](0,0) rectangle (0.2,0.3);
\draw[blue] plot[smooth, tension=.7] coordinates {(0,0) (0.06,0.15) (0,0.3)};
\draw[blue] plot[smooth, tension=.7] coordinates {(0.2,0) (0.14,0.15) (0.2,0.3)};
\end{scope}
\begin{scope}[xshift=0.4cm,yshift=0.6cm]%
  \filldraw[white](0,0) rectangle (0.2,0.3);
\draw[blue] plot[smooth, tension=.7] coordinates {(0,0) (0.2,0.3)};
\draw[blue] plot[smooth, tension=.7] coordinates {(0.2,0) (0.15,0.075)};
\draw[blue] plot[smooth, tension=.7] coordinates {(0.05,0.225) (0,0.3)};
\end{scope}
\node at (1.7,1.5) {\scriptsize$[-2]$};
\draw[blue,->] (2.2,1.5) -- (2.8,1.5);
\begin{scope}[xshift=3cm]
\draw[blue] plot[smooth, tension=.7] coordinates {(0,0) (2,3)};
\draw[blue] plot[smooth, tension=.7] coordinates {(2,0) (0,3)};
\draw[blue] plot[smooth, tension=.2] coordinates {(1,0) (0.1,1.35) (0.1,1.65) (1,3)};
\end{scope}
\begin{scope}[xshift=3.4cm,yshift=2.11cm]%
  \filldraw[white](0,0) rectangle (0.2,0.3);
\draw[blue] plot[smooth, tension=.7] coordinates {(0,0) (0.1,0.1) (0.2,0)};
\draw[blue] plot[smooth, tension=.7] coordinates {(0,0.3) (0.1,0.2) (0.2,0.3)};
\end{scope}
\begin{scope}[xshift=3.9cm,yshift=1.35cm]%
\filldraw[white](0,0) rectangle (0.2,0.3);
\draw[blue] plot[smooth, tension=.7] coordinates {(0,0) (0.06,0.15) (0,0.3)};
\draw[blue] plot[smooth, tension=.7] coordinates {(0.2,0) (0.14,0.15) (0.2,0.3)};
\end{scope}
\begin{scope}[xshift=3.4cm,yshift=0.6cm]%
  \filldraw[white](0,0) rectangle (0.2,0.3);
\draw[blue] plot[smooth, tension=.7] coordinates {(0,0) (0.2,0.3)};
\draw[blue] plot[smooth, tension=.7] coordinates {(0.2,0) (0.15,0.075)};
\draw[blue] plot[smooth, tension=.7] coordinates {(0.05,0.225) (0,0.3)};
\end{scope}
\node at (5.2,1.6) {\scriptsize$t^{1/2}[-1]$};
\begin{scope}[yshift=-4cm]
\draw[blue] plot[smooth, tension=.7] coordinates {(0,0) (2,3)};
\draw[blue] plot[smooth, tension=.7] coordinates {(2,0) (0,3)};
\draw[blue] plot[smooth, tension=.2] coordinates {(1,0) (0.1,1.35) (0.1,1.65) (1,3)};
\begin{scope}[xshift=0.4cm,yshift=2.11cm]%
  \filldraw[white](0,0) rectangle (0.2,0.3);
\draw[blue] plot[smooth, tension=.7] coordinates {(0,0) (0.06,0.15) (0,0.3)};
\draw[blue] plot[smooth, tension=.7] coordinates {(0.2,0) (0.14,0.15) (0.2,0.3)};
\end{scope}
\begin{scope}[xshift=0.9cm,yshift=1.35cm]%
\filldraw[white](0,0) rectangle (0.2,0.3);
\draw[blue] plot[smooth, tension=.7] coordinates {(0,0) (0.1,0.1) (0.2,0)};
\draw[blue] plot[smooth, tension=.7] coordinates {(0,0.3) (0.1,0.2) (0.2,0.3)};
\end{scope}
\begin{scope}[xshift=0.4cm,yshift=0.6cm]%
  \filldraw[white](0,0) rectangle (0.2,0.3);
\draw[blue] plot[smooth, tension=.7] coordinates {(0,0) (0.2,0.3)};
\draw[blue] plot[smooth, tension=.7] coordinates {(0.2,0) (0.15,0.075)};
\draw[blue] plot[smooth, tension=.7] coordinates {(0.05,0.225) (0,0.3)};
\end{scope}
\node at (1.7,1.5) {\scriptsize$[-1]$};
\node at (-0.4,1.5) {\scriptsize$t^{1/2}$};
\draw[blue,->] (2.2,1.5) -- (2.8,1.5);
\begin{scope}[xshift=3cm]
\draw[blue] plot[smooth, tension=.7] coordinates {(0,0) (2,3)};
\draw[blue] plot[smooth, tension=.7] coordinates {(2,0) (0,3)};
\draw[blue] plot[smooth, tension=.2] coordinates {(1,0) (0.1,1.35) (0.1,1.65) (1,3)};
\end{scope}
\begin{scope}[xshift=3.4cm,yshift=2.11cm]%
\filldraw[white](0,0) rectangle (0.2,0.3);
\draw[blue] plot[smooth, tension=.7] coordinates {(0,0) (0.1,0.1) (0.2,0)};
\draw[blue] plot[smooth, tension=.7] coordinates {(0,0.3) (0.1,0.2) (0.2,0.3)};
\end{scope}
\begin{scope}[xshift=3.9cm,yshift=1.35cm]%
\filldraw[white](0,0) rectangle (0.2,0.3);
\draw[blue] plot[smooth, tension=.7] coordinates {(0,0) (0.1,0.1) (0.2,0)};
\draw[blue] plot[smooth, tension=.7] coordinates {(0,0.3) (0.1,0.2) (0.2,0.3)};
\end{scope}
\begin{scope}[xshift=3.4cm,yshift=0.6cm]%
  \filldraw[white](0,0) rectangle (0.2,0.3);
\draw[blue] plot[smooth, tension=.7] coordinates {(0,0) (0.2,0.3)};
\draw[blue] plot[smooth, tension=.7] coordinates {(0.2,0) (0.15,0.075)};
\draw[blue] plot[smooth, tension=.7] coordinates {(0.05,0.225) (0,0.3)};
\end{scope}
\node at (4.7,1.5) {\scriptsize$t$};
\end{scope}
\node at (1,-0.5) {{\color{blue}\scriptsize$\downarrow$}};
\node at (4,-0.5) {{\color{blue}\scriptsize$\downarrow$}};
\end{tikzpicture}
\end{center}
Finally applying DT to the bottom crossing, we get a triply iterated cone, drawn as a cube
\begin{center}
\begin{tikzpicture}[scale=0.6]
\begin{scope}[scale=0.7,xshift=-1cm,yshift=4cm]
\draw[blue] plot[smooth, tension=.7] coordinates {(0,0) (2,3)};
\draw[blue] plot[smooth, tension=.7] coordinates {(2,0) (0,3)};
\draw[blue] plot[smooth, tension=.2] coordinates {(1,0) (0.1,1.35) (0.1,1.65) (1,3)};
\begin{scope}[xshift=0.4cm,yshift=2.11cm]%
  \filldraw[white](0,0) rectangle (0.2,0.3);
\draw[blue] plot[smooth, tension=.7] coordinates {(0,0) (0.06,0.15) (0,0.3)};
\draw[blue] plot[smooth, tension=.7] coordinates {(0.2,0) (0.14,0.15) (0.2,0.3)};
\end{scope}
\begin{scope}[xshift=0.9cm,yshift=1.35cm]%
\filldraw[white](0,0) rectangle (0.2,0.3);
\draw[blue] plot[smooth, tension=.7] coordinates {(0,0) (0.06,0.15) (0,0.3)};
\draw[blue] plot[smooth, tension=.7] coordinates {(0.2,0) (0.14,0.15) (0.2,0.3)};
\end{scope}
\begin{scope}[xshift=0.4cm,yshift=0.6cm] %
\filldraw[white](0,0) rectangle (0.2,0.3);
\draw[blue] plot[smooth, tension=.7] coordinates {(0,0) (0.06,0.15) (0,0.3)};
\draw[blue] plot[smooth, tension=.7] coordinates {(0.2,0) (0.14,0.15) (0.2,0.3)};
\end{scope}
\node at (1.5,0.8) {\scriptsize$[-3]$};
\end{scope}
\begin{scope}[scale=0.7,xshift=4.2cm,yshift=4cm]
  \draw[blue] plot[smooth, tension=.7] coordinates {(0,0) (2,3)};
\draw[blue] plot[smooth, tension=.7] coordinates {(2,0) (0,3)};
\draw[blue] plot[smooth, tension=.2] coordinates {(1,0) (0.1,1.35) (0.1,1.65) (1,3)};
\begin{scope}[xshift=.4cm,yshift=2.11cm]%
  \filldraw[white](0,0) rectangle (0.2,0.3);
\draw[blue] plot[smooth, tension=.7] coordinates {(0,0) (0.1,0.1) (0.2,0)};
\draw[blue] plot[smooth, tension=.7] coordinates {(0,0.3) (0.1,0.2) (0.2,0.3)};
\end{scope}
\begin{scope}[xshift=.9cm,yshift=1.35cm]%
\filldraw[white](0,0) rectangle (0.2,0.3);
\draw[blue] plot[smooth, tension=.7] coordinates {(0,0) (0.06,0.15) (0,0.3)};
\draw[blue] plot[smooth, tension=.7] coordinates {(0.2,0) (0.14,0.15) (0.2,0.3)};
\end{scope}
\begin{scope}[xshift=.4cm,yshift=0.6cm]%
\filldraw[white](0,0) rectangle (0.2,0.3);
\draw[blue] plot[smooth, tension=.7] coordinates {(0,0) (0.06,0.15) (0,0.3)};
\draw[blue] plot[smooth, tension=.7] coordinates {(0.2,0) (0.14,0.15) (0.2,0.3)};
\end{scope}
\node at (2.2,1) {\scriptsize$t^{1/2}[-2]$};
\end{scope}%
 \begin{scope}[scale=0.7,xshift=-1cm,yshift=-1.5cm]
\draw[blue] plot[smooth, tension=.7] coordinates {(0,0) (2,3)};
\draw[blue] plot[smooth, tension=.7] coordinates {(2,0) (0,3)};
\draw[blue] plot[smooth, tension=.2] coordinates {(1,0) (0.1,1.35) (0.1,1.65) (1,3)};
\begin{scope}[xshift=0.4cm,yshift=2.11cm]%
  \filldraw[white](0,0) rectangle (0.2,0.3);
\draw[blue] plot[smooth, tension=.7] coordinates {(0,0) (0.06,0.15) (0,0.3)};
\draw[blue] plot[smooth, tension=.7] coordinates {(0.2,0) (0.14,0.15) (0.2,0.3)};
\end{scope}
\begin{scope}[xshift=0.9cm,yshift=1.35cm]%
\filldraw[white](0,0) rectangle (0.2,0.3);
\draw[blue] plot[smooth, tension=.7] coordinates {(0,0) (0.1,0.1) (0.2,0)};
\draw[blue] plot[smooth, tension=.7] coordinates {(0,0.3) (0.1,0.2) (0.2,0.3)};
\end{scope}
\begin{scope}[xshift=0.4cm,yshift=0.6cm]%
\filldraw[white](0,0) rectangle (0.2,0.3);
\draw[blue] plot[smooth, tension=.7] coordinates {(0,0) (0.06,0.15) (0,0.3)};
\draw[blue] plot[smooth, tension=.7] coordinates {(0.2,0) (0.14,0.15) (0.2,0.3)};
\end{scope}
\node at (0.7,0) {\scriptsize$t^{1/2}[-2]$};
\end{scope}
\begin{scope}[scale=0.7,xshift=4.3cm,yshift=-1.5]
\draw[blue] plot[smooth, tension=.7] coordinates {(0,-1.5) (2,1.5)};
\draw[blue] plot[smooth, tension=.7] coordinates {(2,-1.5) (0,1.5)};
\draw[blue] plot[smooth, tension=.2] coordinates {(1,-1.5) (0.1,-0.15) (0.1,0.15) (1,1.5)};

\begin{scope}[xshift=.4cm,yshift=.61cm]%
\filldraw[white](0,0) rectangle (0.2,0.3);
\draw[blue] plot[smooth, tension=.7] coordinates {(0,0) (0.1,0.1) (0.2,0)};
\draw[blue] plot[smooth, tension=.7] coordinates {(0,0.3) (0.1,0.2) (0.2,0.3)};
\end{scope}
\begin{scope}[xshift=.9cm,yshift=-0.15cm]%
\filldraw[white](0,0) rectangle (0.2,0.3);
\draw[blue] plot[smooth, tension=.7] coordinates {(0,0) (0.1,0.1) (0.2,0)};
\draw[blue] plot[smooth, tension=.7] coordinates {(0,0.3) (0.1,0.2) (0.2,0.3)};
\end{scope}
\begin{scope}[xshift=.4cm,yshift=-0.9cm]%
\filldraw[white](0,0) rectangle (0.2,0.3);
\draw[blue] plot[smooth, tension=.7] coordinates {(0,0) (0.06,0.15) (0,0.3)};
\draw[blue] plot[smooth, tension=.7] coordinates {(0.2,0) (0.14,0.15) (0.2,0.3)};
\end{scope}
\node at (1.5,-0.8) {\scriptsize$t[-1]$};
\end{scope}
\begin{scope}[scale=0.8,xshift=-3.5cm,yshift=0.2cm]
\draw[red] plot[smooth, tension=.7] coordinates {(0,0) (2,3)};
\draw[red] plot[smooth, tension=.7] coordinates {(2,0) (0,3)};
\draw[red] plot[smooth, tension=.2] coordinates {(1,0) (0.1,1.35) (0.1,1.65) (1,3)};
\begin{scope}[xshift=0.4cm,yshift=2.11cm]%
  \filldraw[white](0,0) rectangle (0.2,0.3);
\draw[red] plot[smooth, tension=.7] coordinates {(0,0) (0.06,0.15) (0,0.3)};
\draw[red] plot[smooth, tension=.7] coordinates {(0.2,0) (0.14,0.15) (0.2,0.3)};
\end{scope}
\begin{scope}[xshift=0.9cm,yshift=1.35cm]%
\filldraw[white](0,0) rectangle (0.2,0.3);
\draw[red] plot[smooth, tension=.7] coordinates {(0,0) (0.06,0.15) (0,0.3)};
\draw[red] plot[smooth, tension=.7] coordinates {(0.2,0) (0.14,0.15) (0.2,0.3)};
\end{scope}
\begin{scope}[xshift=0.4cm,yshift=0.6cm] %
\filldraw[white](0,0) rectangle (0.2,0.3);
\draw[red] plot[smooth, tension=.7] coordinates {(0,0) (0.1,0.1) (0.2,0)};
\draw[red] plot[smooth, tension=.7] coordinates {(0,0.3) (0.1,0.2) (0.2,0.3)};
\end{scope}
\node at (-0.8,1.2) {\scriptsize$t^{1/2}[-2]$};
\end{scope}
\begin{scope}[scale=0.8,xshift=1cm,yshift=0.2cm]
  \draw[red] plot[smooth, tension=.7] coordinates {(0,0) (2,3)};
\draw[red] plot[smooth, tension=.7] coordinates {(2,0) (0,3)};
\draw[red] plot[smooth, tension=.2] coordinates {(1,0) (0.1,1.35) (0.1,1.65) (1,3)};
\begin{scope}[xshift=.4cm,yshift=2.11cm]%
\filldraw[white](0,0) rectangle (0.2,0.3);
\draw[red] plot[smooth, tension=.7] coordinates {(0,0) (0.1,0.1) (0.2,0)};
\draw[red] plot[smooth, tension=.7] coordinates {(0,0.3) (0.1,0.2) (0.2,0.3)};
\end{scope}
\begin{scope}[xshift=.9cm,yshift=1.35cm]%
\filldraw[white](0,0) rectangle (0.2,0.3);
\draw[red] plot[smooth, tension=.7] coordinates {(0,0) (0.06,0.15) (0,0.3)};
\draw[red] plot[smooth, tension=.7] coordinates {(0.2,0) (0.14,0.15) (0.2,0.3)};
\end{scope}
\begin{scope}[xshift=.4cm,yshift=0.6cm]%
\filldraw[white](0,0) rectangle (0.2,0.3);
\draw[red] plot[smooth, tension=.7] coordinates {(0,0) (0.1,0.1) (0.2,0)};
\draw[red] plot[smooth, tension=.7] coordinates {(0,0.3) (0.1,0.2) (0.2,0.3)};
\end{scope}
\node at (1.5,0.9) {\scriptsize$t[-1]$};
\end{scope}%
 \begin{scope}[scale=0.8,xshift=-3.2cm,yshift=-4.7cm]
\draw[red] plot[smooth, tension=.7] coordinates {(0,0) (2,3)};
\draw[red] plot[smooth, tension=.7] coordinates {(2,0) (0,3)};
\draw[red] plot[smooth, tension=.2] coordinates {(1,0) (0.1,1.35) (0.1,1.65) (1,3)};
\begin{scope}[xshift=0.4cm,yshift=2.11cm]%
  \filldraw[white](0,0) rectangle (0.2,0.3);
\draw[red] plot[smooth, tension=.7] coordinates {(0,0) (0.06,0.15) (0,0.3)};
\draw[red] plot[smooth, tension=.7] coordinates {(0.2,0) (0.14,0.15) (0.2,0.3)};
\end{scope}
\begin{scope}[xshift=0.9cm,yshift=1.35cm]%
\filldraw[white](0,0) rectangle (0.2,0.3);
\draw[red] plot[smooth, tension=.7] coordinates {(0,0) (0.1,0.1) (0.2,0)};
\draw[red] plot[smooth, tension=.7] coordinates {(0,0.3) (0.1,0.2) (0.2,0.3)};
\end{scope}
\begin{scope}[xshift=0.4cm,yshift=0.6cm]%
\filldraw[white](0,0) rectangle (0.2,0.3);
\draw[red] plot[smooth, tension=.7] coordinates {(0,0) (0.1,0.1) (0.2,0)};
\draw[red] plot[smooth, tension=.7] coordinates {(0,0.3) (0.1,0.2) (0.2,0.3)};\end{scope}
\node at (1.5,0.9) {\scriptsize$t[-1]$};
\end{scope}
\begin{scope}[scale=0.8,xshift=1.4cm,yshift=-3.2cm]
\draw[red] plot[smooth, tension=.7] coordinates {(0,-1.5) (2,1.5)};
\draw[red] plot[smooth, tension=.7] coordinates {(2,-1.5) (0,1.5)};
\draw[red] plot[smooth, tension=.2] coordinates {(1,-1.5) (0.1,-0.15) (0.1,0.15) (1,1.5)};

\begin{scope}[xshift=.4cm,yshift=.61cm]%
\filldraw[white](0,0) rectangle (0.2,0.3);
\draw[red] plot[smooth, tension=.7] coordinates {(0,0) (0.1,0.1) (0.2,0)};
\draw[red] plot[smooth, tension=.7] coordinates {(0,0.3) (0.1,0.2) (0.2,0.3)};
\end{scope}
\begin{scope}[xshift=.9cm,yshift=-0.15cm]%
\filldraw[white](0,0) rectangle (0.2,0.3);
\draw[red] plot[smooth, tension=.7] coordinates {(0,0) (0.1,0.1) (0.2,0)};
\draw[red] plot[smooth, tension=.7] coordinates {(0,0.3) (0.1,0.2) (0.2,0.3)};
\end{scope}
\begin{scope}[xshift=.4cm,yshift=-0.9cm]%
\filldraw[white](0,0) rectangle (0.2,0.3);
\draw[red] plot[smooth, tension=.7] coordinates {(0,0) (0.1,0.1) (0.2,0)};
\draw[red] plot[smooth, tension=.7] coordinates {(0,0.3) (0.1,0.2) (0.2,0.3)};\end{scope}
\node at (1.6,-0.3) {\scriptsize$t^{3/2}$};
\end{scope}
\draw[black,->] (-0.2,2.7) -- (-0.2,1); \filldraw[white] (-0.3,1.3) rectangle (-0.1,1.5);
\draw[black,->] (0.3,0) -- (2.6,0); \draw[black,->] (3,-0.6) -- (1.8,-2.3);
\draw[black,->] (-0.75,-0.3) -- (-2,-2.25); 
\draw[black,->] (-1.5,-2.6) -- (0.8,-2.6);
\begin{scope}[xshift=-0.2cm,yshift=4cm]
\draw[black,->] (0.3,0) -- (2.6,0); \draw[green,->] (3,-0.6) -- (1.8,-2.3);
\draw[black,->] (-0.5,-0.3) -- (-1.75,-2.25); 
\draw[black,->] (-1.5,-2.6) -- (0.8,-2.6);
\end{scope}
\filldraw[white] (1.5,-0.05) rectangle (1.7,0.05);
\begin{scope}[xshift=4cm]
\draw[purple,->] (-0.2,2.7) -- (-0.2,1);
\end{scope}
\draw[green,->] (1.6,0.2) -- (1.6,-1.5);
\begin{scope}[xshift=-3.5cm,yshift=0cm]
  \draw[black,->] (1.6,0.2) -- (1.6,-1.5);
\end{scope}
\end{tikzpicture}
\end{center}
The three terms connected by the green arrows can be cancelled via a chain homotopy, since the middle one is actually the direct sum of the two other. After cancellation, we are left with the diagram
\begin{center}
\begin{tikzpicture}
  \draw[blue,-] (0,-0.5) -- (0,0.5); \draw[blue,-] (0.3,-0.5) -- (0.3,0.5); \draw[blue,-] (0.6,-0.5) -- (0.6,0.5);
  \node at (.9,0) {\scriptsize[-3]};
  
  \begin{scope}[xshift=1.5cm,yshift=.6cm]
    \draw[blue] plot[smooth, tension=.5] coordinates {(0,0.5) (0.25,0.1) (0.5,0.5)};
  \draw[blue] plot[smooth, tension=.5] coordinates {(0,-0.5) (0.25,-0.1) (0.5,-0.5)};
  \draw[blue,-] (0.6,-0.5) -- (0.6,0.5);
  \node at (1.2,0) {\scriptsize{$t^{1/2}[-2]$}};
  \end{scope}
  \begin{scope}[xshift=1.5cm,yshift=-0.6cm]
  \draw[blue] plot[smooth, tension=.5] coordinates {(0.6,0.5) (0.35,0.1) (0.1,0.5)};
  \draw[blue] plot[smooth, tension=.5] coordinates {(0.6,-0.5) (0.35,-0.1) (0.1,-0.5)};
  \draw[blue,-] (0,-0.5) -- (0,0.5);
  \node at (1.2,0) {\scriptsize{$t^{1/2}[-2]$}};
  \end{scope}
  \begin{scope}[xshift=4cm,yshift=-.6cm]
    \draw[blue] plot[smooth, tension=.5] coordinates {(0,0.5) (0.25,0.1) (0.5,0.5)};
  \draw[blue] plot[smooth, tension=.5] coordinates {(0.6,-0.5) (0.35,-0.1) (0.1,-0.5)};
  \draw[blue] plot[smooth, tension=.5] coordinates {(-0.1,-0.5) (0.25,0.05) (0.35,-0.05) (0.7,0.5)};
  \node at (1,0) {\scriptsize{$t[-1]$}};
  \end{scope}
  \begin{scope}[xshift=4cm,yshift=0.6cm]
  \draw[blue] plot[smooth, tension=.5] coordinates {(0.6,0.5) (0.35,0.1) (0.1,0.5)};
  \draw[blue] plot[smooth, tension=.5] coordinates {(0,-0.5) (0.25,-0.1) (0.5,-0.5)};
  \draw[blue] plot[smooth, tension=.5] coordinates {(0.7,-0.5) (0.35,0.05) (0.25,-0.05) (-0.1,0.5)};
  \node at (1,0) {\scriptsize{$t[-1]$}};
  \end{scope}
  \draw[blue,->] (0.8,0.25) -- (1.5,0.6);
  \draw[blue,->] (0.8,-0.25) -- (1.5,-0.6);
  \draw[blue,->] (3.3,0.65) -- (4,0.65);
  \draw[blue,->] (3.3,-0.65) -- (4,-0.65);
  \draw[blue,->,dotted] (3.3,0.6) -- (4,-0.6);
  \draw[blue,->] (3.3,-0.6) -- (4,0.6);
\end{tikzpicture}
\end{center}
Note that the dashed arrow did not exist in the original cube, it came about as a result of rerouting the purple map in the cube during the chain homotopy. If we carry out the same exercise for the rhs of the Reidemeister III move, we will get exactly the same diagram with the same maps, thereby finishing the proof.
\end{proof}
\begin{remark}
  The reader will immediately recognise that the proof runs along the logic of the Khovanov homology, by resolving each crossing in two ways. See also sec.\ref{sec_RtKh} for an actual Khovanov homology calculation.
\end{remark}

\subsection{Relation to Khovanov Homology}\label{sec_RtKh}
In this example, we perform the computation in a way that makes the relation to Khovanov homology \cite{khovanov2000categorification} \cite{bar2002khovanov} more evident. Here we follow Bar-Natan's exposition \cite{bar2002khovanov} of Khovanov's categorification of the Jones polynomial as a categorification of the Kauffman bracket. One resolves each crossing in two ways, and the result of the resolution is a number of disjoint circles, each of which is assigned a dim 2 vector space $C$, whose base vectors have equivariant weights $t^{1/2}$ and $t^{-1/2}$, respectively.\footnote{Our use of $t$ and $q$ is opposite that of \cite{bar2002khovanov}.} Each vertex of the cube corresponds to one of the $2^{\#{\rm crossing}}$ ways of resolution. An edge would then connect two different resolutions that differ only in one crossing, and to this edge is assigned a map of $C^{\otimes m}\to C^{\otimes n}$ where $m,n$ are the numbers of disjoint circles in the two resolutions. These maps together give the differential of the chain complex coming from flattening the cube. The equivariant Poincar\'e polynomial then is a categorification of Jones polynomial.

For us, we have already categorified the Kauffman bracket in prop.\ref{prop_Kauffman}. But we rotate the DT's
\bea 
\begin{tikzpicture}[scale=.8]
\draw [blue,-] (2,0) -- (3,1);
\draw [white,-,line width=2mm] (3,0) -- (2,1);
\draw [blue,-] (3,0) -- (2,1);

\node at (1.3,0.5) {\scriptsize$\longrightarrow~~$};

\centerarc[blue](-0.5,0.5)(65:-65:0.5);
\centerarc[blue](0.8,0.5)(115:245:0.5);

\node at (-0.5,0.5) {\scriptsize$\longrightarrow~~$};
\node at (-1.3,0.5) {\scriptsize$[-1]$};

\node at (-2.7,0.5) {\scriptsize$t^{-1/2}$};
\centerarc[blue](-2,1.2)(-35:-145:0.5);
\centerarc[blue](-2,-0.2)(35:145:0.5);
\node at (3.5,0.6) {\scriptsize{$\stackrel{+1}{\longrightarrow}$}};
\node at (2,0.5) {\scriptsize$t^{-1/4}$};
\end{tikzpicture}\nn\\
\begin{tikzpicture}[scale=.8]
\draw [blue,-] (3,0) -- (2,1);
\draw [white,-,line width=2mm] (2,0) -- (3,1);
\draw [blue,-] (2,0) -- (3,1);

\node at (-0.1,0.5) {\scriptsize$t^{1/2}$};

\centerarc[blue](-2.8,0.5)(65:-65:0.5);
\centerarc[blue](-1.5,0.5)(115:245:0.5);
\node at (-1.5,0.5) {\scriptsize$[-1]$};

\node at (-.7,0.5) {\scriptsize$\longrightarrow~~$};
\node at (1.4,0.5) {\scriptsize$\longrightarrow$};
\node at (3.5,0.5) {\scriptsize$\stackrel{+1}{\longrightarrow}$};

\centerarc[blue](.5,1.2)(-35:-145:0.5);
\centerarc[blue](.5,-0.2)(35:145:0.5);
\node at (2,0.5) {\scriptsize$t^{1/4}$};
\end{tikzpicture}\label{eq_used}\eea
so as to suggest that each crossing implies a map between its two resolutions.

For the remainder of the calculation, we take
the trefoil knot as example.
\bea 
\begin{tikzpicture}[scale=0.8]
\draw[blue] (0,0) -- (0.8,-0.8) -- (0.8,-1.8) -- (0.1,-1.1);
\draw[blue] (-0.1,-0.9) -- (-0.4,-0.6) -- (-0.4,0.4) -- (0.4,0.4) -- (0.4,-0.3);
\draw[blue] (0.4,-0.5) -- (0.4,-0.6) -- (-0.8,-1.8) -- (-0.8,-0.8) -- (-0.5,-0.5);
\draw[blue] (-0.3,-0.3) -- (0,0);
\node at (-1.5,-.5) {\scriptsize$t^{-1/4}$};
\end{tikzpicture}\nn\eea
and by resolving  each crossing in two ways, we get a cube
\bea 
\begin{tikzpicture}[scale=1]
  \matrix (m) [matrix of math nodes, row sep=2em,
    column sep=2em]{
    & 1 & & 2 \\
    5 & & 6 & \\
    & 3 & & 4 \\
    7 & & 8 & \\};
  \path[-stealth]
    (m-1-2) edge (m-1-4) edge (m-2-1)
            edge [densely dotted] (m-3-2)
    (m-1-4) edge (m-3-4) edge (m-2-3)
    (m-2-1) edge [-,line width=6pt,draw=white] (m-2-3)
            edge (m-2-3) edge (m-4-1)
    (m-3-2) edge [densely dotted] (m-3-4)
            edge [densely dotted] (m-4-1)
    (m-4-1) edge (m-4-3)
    (m-3-4) edge (m-4-3)
    (m-2-3) edge [-,line width=6pt,draw=white] (m-4-3)
            edge (m-4-3);
            
\begin{scope}[scale=0.6,xshift=-2cm,yshift=4cm]
 \draw[blue] (0,0) -- (0.8,-0.8) -- (0.8,-1.8) -- (0.1,-1.1);
\draw[blue] (-0.1,-0.9) -- (-0.4,-0.6) -- (-0.4,0.4) -- (0.4,0.4) -- (0.4,-0.3);
\draw[blue] (0.4,-0.5) -- (0.4,-0.6) -- (-0.8,-1.8) -- (-0.8,-0.8) -- (-0.5,-0.5);
\draw[blue] (-0.3,-0.3) -- (0,0);
\draw[fill,white] (0.3,-0.3) rectangle (0.5,-0.5);
\draw[fill,white] (-0.3,-0.3) rectangle (-0.5,-0.5);
\draw[fill,white] (-0.1,-0.9) rectangle (0.1,-1.1);
\draw[blue] (-0.1,-0.9) -- (-0.1,-1.1);
\draw[blue] (0.1,-0.9) -- (0.1,-1.1);
\draw[blue] (-0.3,-0.3) -- (-0.4,-0.3);
\draw[blue] (-0.4,-0.5) -- (-0.5,-0.5);
\draw[blue] (0.3,-0.3) -- (0.4,-0.3);
\draw[blue] (0.4,-0.5) -- (0.5,-0.5);
\node at (1.8,-0.6) {\scriptsize{$t^{-1}[-3]$}};
\end{scope}
\begin{scope}[scale=0.6,xshift=4cm,yshift=4cm]
  \draw[blue] (0,0) -- (0.8,-0.8) -- (0.8,-1.8) -- (0.1,-1.1);
\draw[blue] (-0.1,-0.9) -- (-0.4,-0.6) -- (-0.4,0.4) -- (0.4,0.4) -- (0.4,-0.3);
\draw[blue] (0.4,-0.5) -- (0.4,-0.6) -- (-0.8,-1.8) -- (-0.8,-0.8) -- (-0.5,-0.5);
\draw[blue] (-0.3,-0.3) -- (0,0);
\draw[fill,white] (0.3,-0.3) rectangle (0.5,-0.5);
\draw[fill,white] (-0.3,-0.3) rectangle (-0.5,-0.5);
\draw[fill,white] (-0.1,-0.9) rectangle (0.1,-1.1);
\draw[blue] (-0.1,-0.9) -- (-0.1,-1.1);
\draw[blue] (0.1,-0.9) -- (0.1,-1.1);
\draw[blue] (-0.3,-0.3) -- (-0.4,-0.3);
\draw[blue] (-0.4,-0.5) -- (-0.5,-0.5);
\draw[blue] (0.3,-0.3) -- (0.4,-0.5);
\draw[blue] (0.4,-0.3) -- (0.5,-0.5);
\node at (1.9,-0.6) {\scriptsize{$t^{-1/2}[-2]$}};
\end{scope}
\begin{scope}[scale=0.6,xshift=-2cm,yshift=0.5cm]
  \draw[blue] (0,0) -- (0.8,-0.8) -- (0.8,-1.8) -- (0.1,-1.1);
\draw[blue] (-0.1,-0.9) -- (-0.4,-0.6) -- (-0.4,0.4) -- (0.4,0.4) -- (0.4,-0.3);
\draw[blue] (0.4,-0.5) -- (0.4,-0.6) -- (-0.8,-1.8) -- (-0.8,-0.8) -- (-0.5,-0.5);
\draw[blue] (-0.3,-0.3) -- (0,0);
\draw[fill,white] (0.3,-0.3) rectangle (0.5,-0.5);
\draw[fill,white] (-0.3,-0.3) rectangle (-0.5,-0.5);
\draw[fill,white] (-0.1,-0.9) rectangle (0.1,-1.1);
\draw[blue] (-0.1,-0.9) -- (-0.1,-1.1);
\draw[blue] (0.1,-0.9) -- (0.1,-1.1);
\draw[blue] (-0.3,-0.3) -- (-0.4,-0.5);
\draw[blue] (-0.4,-0.3) -- (-0.5,-0.5);
\draw[blue] (0.3,-0.3) -- (0.4,-0.3);
\draw[blue] (0.4,-0.5) -- (0.5,-0.5);
\node at (1.9,-0.6) {\scriptsize{$t^{-1/2}[-2]$}};
\end{scope}
\begin{scope}[scale=0.6,xshift=4cm,yshift=0.5cm]
\draw[blue] (0,0) -- (0.8,-0.8) -- (0.8,-1.8) -- (0.1,-1.1);
\draw[blue] (-0.1,-0.9) -- (-0.4,-0.6) -- (-0.4,0.4) -- (0.4,0.4) -- (0.4,-0.3);
\draw[blue] (0.4,-0.5) -- (0.4,-0.6) -- (-0.8,-1.8) -- (-0.8,-0.8) -- (-0.5,-0.5);
\draw[blue] (-0.3,-0.3) -- (0,0);
\draw[fill,white] (0.3,-0.3) rectangle (0.5,-0.5);
\draw[fill,white] (-0.3,-0.3) rectangle (-0.5,-0.5);
\draw[fill,white] (-0.1,-0.9) rectangle (0.1,-1.1);
\draw[blue] (-0.1,-0.9) -- (-0.1,-1.1);
\draw[blue] (0.1,-0.9) -- (0.1,-1.1);
\draw[blue] (-0.3,-0.3) -- (-0.4,-0.5);
\draw[blue] (-0.4,-0.3) -- (-0.5,-0.5);
\draw[blue] (0.3,-0.3) -- (0.4,-0.5);
\draw[blue] (0.4,-0.3) -- (0.5,-0.5);
\node at (1.3,-0.6) {\scriptsize{$[-1]$}};
\end{scope}
\begin{scope}[red,scale=0.6,xshift=-4cm,yshift=2cm]
  \draw[red] (0,0) -- (0.8,-0.8) -- (0.8,-1.8) -- (0.1,-1.1);
\draw[red] (-0.1,-0.9) -- (-0.4,-0.6) -- (-0.4,0.4) -- (0.4,0.4) -- (0.4,-0.3);
\draw[red] (0.4,-0.5) -- (0.4,-0.6) -- (-0.8,-1.8) -- (-0.8,-0.8) -- (-0.5,-0.5);
\draw[red] (-0.3,-0.3) -- (0,0);
\draw[fill,white] (0.3,-0.3) rectangle (0.5,-0.5);
\draw[fill,white] (-0.3,-0.3) rectangle (-0.5,-0.5);
\draw[fill,white] (-0.1,-0.9) rectangle (0.1,-1.1);
\draw[red] (-0.1,-0.9) -- (0.1,-0.9);
\draw[red] (-0.1,-1.1) -- (0.1,-1.1);
\draw[red] (-0.3,-0.3) -- (-0.4,-0.3);
\draw[red] (-0.4,-0.5) -- (-0.5,-0.5);
\draw[red] (0.3,-0.3) -- (0.4,-0.3);
\draw[red] (0.4,-0.5) -- (0.5,-0.5);
\node at (1.8,-0.6) {\scriptsize{$t^{-1/2}[-2]$}};
\end{scope}
\begin{scope}[red,scale=0.6,xshift=2cm,yshift=2cm]
  \draw[red] (0,0) -- (0.8,-0.8) -- (0.8,-1.8) -- (0.1,-1.1);
\draw[red] (-0.1,-0.9) -- (-0.4,-0.6) -- (-0.4,0.4) -- (0.4,0.4) -- (0.4,-0.3);
\draw[red] (0.4,-0.5) -- (0.4,-0.6) -- (-0.8,-1.8) -- (-0.8,-0.8) -- (-0.5,-0.5);
\draw[red] (-0.3,-0.3) -- (0,0);
\draw[fill,white] (0.3,-0.3) rectangle (0.5,-0.5);
\draw[fill,white] (-0.3,-0.3) rectangle (-0.5,-0.5);
\draw[fill,white] (-0.1,-0.9) rectangle (0.1,-1.1);
\draw[red] (-0.1,-0.9) -- (0.1,-0.9);
\draw[red] (-0.1,-1.1) -- (0.1,-1.1);
\draw[red] (-0.3,-0.3) -- (-0.4,-0.3);
\draw[red] (-0.4,-0.5) -- (-0.5,-0.5);
\draw[red] (0.3,-0.3) -- (0.4,-0.5);
\draw[red] (0.4,-0.3) -- (0.5,-0.5);
\node at (1.4,-0.6) {\scriptsize{$[-1]$}};
\end{scope}
\begin{scope}[red,scale=0.6,xshift=-4cm,yshift=-2cm]
  \draw[red] (0,0) -- (0.8,-0.8) -- (0.8,-1.8) -- (0.1,-1.1);
\draw[red] (-0.1,-0.9) -- (-0.4,-0.6) -- (-0.4,0.4) -- (0.4,0.4) -- (0.4,-0.3);
\draw[red] (0.4,-0.5) -- (0.4,-0.6) -- (-0.8,-1.8) -- (-0.8,-0.8) -- (-0.5,-0.5);
\draw[red] (-0.3,-0.3) -- (0,0);
\draw[fill,white] (0.3,-0.3) rectangle (0.5,-0.5);
\draw[fill,white] (-0.3,-0.3) rectangle (-0.5,-0.5);
\draw[fill,white] (-0.1,-0.9) rectangle (0.1,-1.1);
\draw[red] (-0.1,-0.9) -- (0.1,-0.9);
\draw[red] (-0.1,-1.1) -- (0.1,-1.1);
\draw[red] (-0.3,-0.3) -- (-0.4,-0.5);
\draw[red] (-0.4,-0.3) -- (-0.5,-0.5);
\draw[red] (0.3,-0.3) -- (0.4,-0.3);
\draw[red] (0.4,-0.5) -- (0.5,-0.5);
\node at (1.4,-0.6) {\scriptsize{$[-1]$}};
\end{scope}
\begin{scope}[red,scale=0.6,xshift=2cm,yshift=-2cm]
\draw[red] (0,0) -- (0.8,-0.8) -- (0.8,-1.8) -- (0.1,-1.1);
\draw[red] (-0.1,-0.9) -- (-0.4,-0.6) -- (-0.4,0.4) -- (0.4,0.4) -- (0.4,-0.3);
\draw[red] (0.4,-0.5) -- (0.4,-0.6) -- (-0.8,-1.8) -- (-0.8,-0.8) -- (-0.5,-0.5);
\draw[red] (-0.3,-0.3) -- (0,0);
\draw[fill,white] (0.3,-0.3) rectangle (0.5,-0.5);
\draw[fill,white] (-0.3,-0.3) rectangle (-0.5,-0.5);
\draw[fill,white] (-0.1,-0.9) rectangle (0.1,-1.1);
\draw[red] (-0.1,-0.9) -- (0.1,-0.9);
\draw[red] (-0.1,-1.1) -- (0.1,-1.1);
\draw[red] (-0.3,-0.3) -- (-0.4,-0.5);
\draw[red] (-0.4,-0.3) -- (-0.5,-0.5);
\draw[red] (0.3,-0.3) -- (0.4,-0.5);
\draw[red] (0.4,-0.3) -- (0.5,-0.5);
\node at (1.3,-0.6) {\scriptsize{$t^{1/2}$}};
\end{scope}
\end{tikzpicture}\nn\eea
All the maps are derived from adding a cup and cap to the graphs \eqref{eq_used}, and we get 
\bea& t^{-1/2}C^{\otimes2}[-1]\stackrel{m}{\to} C \to t^{-1}C[-1]\stackrel{+1}{\to}\nn\\
&C[-1]\stackrel{\Gd}{\to} t^{1/2}C^{\otimes2}\to tC[1]\stackrel{+1}{\to}\nn\eea
where $C=t^{1/2}[1]\oplus t^{-1/2}[-1]$ is the complex for a circle. We have also named the two maps $m,\Gd$ as in \cite{bar2002khovanov}.
That the DT's give rise to long exact sequence of cohomology pins $m,\Gd$ completely. The two DT's above are enough for computing the cohomology of the cube: we replace the back and front face with the chain homotopic complex
\bea {\rm back~face}:&&t^{-2}[-3]C,\nn\\
{\rm front~face}:&&(t^{-1}[-1]+t[1])C.\nn\eea
The cone between the back and front face is
\bea t^{3/2}[2]+t^{1/2}+t^{-1/2}+t^{-5/2}[-3],\nn\eea
Removing the $t^{-1/4}$ prefactor, we get
\bea \Phi(\textrm{trefoil}_L)=t^{7/4}[2]+t^{3/4}+t^{-1/4}+t^{-9/4}[-3].\nn\eea
Similarly, for the right-handed trefoil, one gets $t^2/q^3 + 1/q^2 + 1/t + 1/t^2$
\bea \Phi(\textrm{trefoil}_R)=t^{-7/4}[-2]+t^{-3/4}+t^{1/4}+t^{9/4}[3].\nn\eea

If we descend to the Grothendieck group by setting $[n]=(-1)^n$, and normalise by the result of the unknot, we get $t^{1/4}(t^{-2}-t-t^{-1})$.
Further normalise by the writhe of the knot $(-t^{-3/4})^{-w(L)}=-t^{-9/4}$, we get the Jones polynomial of the trefoil (up to $t\to t^{-1}$ to match the convention)
\bea \textrm{Jones}(\textrm{trefoil}_L)=-t^{-4}+t^{-1}+t^{-3}.\nn\eea

\section{Lifting to the Homotopy Category}\label{sec_LthHC}
\subsection{Motivation}
Our tangle invariants assign to a tangle $K$ a kernel $\Phi(K)$ in $D^b(Y_m\times Y_n)$. Isotopic tangles lead to quasi-isomorphic kernels. However, checking quasi-isomorphisms can be non-trivial. 
What is more, the composition of the kernels is done by a convolution, i.e. involve the pushforward or pullback functors, making practical computations beyond reach. 

We will lift these kernels from the derived category to the homotopy category as we described around \eqref{tilting_equ}, \eqref{proj/htpy}. The first lift is to $K^b(Y_m\times Y_n)$ and then to an object in
\bea \Phi(K)\in K^b(_{\Gl_m}\opn{proj}_{\Gl_n}),\nn\eea
where $D^b(\opn{mod}_{\Gl_n})$ is the category derived equivalent to $D^b(Y_n)$ via tilting. Now quasi-isomorphisms will correspond to chain homotopies, which are much easier to check. Moreover, the convolution operation, which can be hard to understand in the derived category, becomes the usual tensor product of (complexes of) bi-modules.

These motivations here are similar to those of Rouquier in \cite{rouquier2006categorification}, namely, to categorify the braid group action by means of \emph{Rouquier complexes}. Let us recall this idea: one considers $\FR{sl}(n)$, with Weyl group $S_{n+1}$ and natural action on the polynomial algebra $R=k[x_0,\cdots,x_n]$. This also induces a natural action on the category of $R$-modules, and to each Weyl reflection $s_i$ (permutation of $x_i,x_{i+1}$), Rouquier assigns a chain complex of $R$-bi-modules that serves as the FM kernel of $s_i$.

\subsection{A matrix presentation of kernels}
To apply the same procedure to our setting, we proceed as follows. In sec.\ref{sec_ECoDY} we find the exceptional collection $\{L_I\}$ for $Y_n$ and the tilting bundle $T_n=\oplus_IL_I$. The endomorphism algebra of $T_n$ is a quiver algebra $\Gl_n$. Tilting then turns $D^b(Y_n)$ into the equivalent module category $D^b(\opn{mod}_{\Gl_n})$: the exceptional objects $\{L_I\}$ pass over to the projective right modules $\{e_i\Gl_n\}$ with $e_i$ being the idempotents of $\Gl_n$. Likewise, the dual bundles $\{L_I^{\vee}\}$ pass over to the left modules $\Gl_n e_i$; the hom space $\hom_{D^b(Y_n)}(L_I,L_J)$ is isomorphic to the vector space $e_j\Gl_n e_i$, the space of paths from $e_i$ to $e_j$ (and there is no higher hom).

Recall that $\Gl_n$ is the kernel of the identity FM transform $D^b(\opn{mod}_{\Gl_n})\to D^b(\opn{mod}_{\Gl_n})$.
Assuming that we found the resolution of $\Gl_n$ as a $\Gl_n$-bi-module in terms of the left/right projectives $\Gl_n e_i$ and $e_i \Gl_n$, then this resolution serves as the Rouquier complex for the identity FM transform. Furthermore, from this kernel, we can easily find the Rouquier complex for the other kernels. To clarify the idea, we suggestively write the resolution of $\Gl$ (we will stop writing the index of $\Gl$) as
\bea \Gl \sim \Gl e_i\Gd_{ij} e_j\Gl\nn\eea
as if it were a matrix.Here the entries $\Gd_{ij}$ are (cohomologically and equivariantly) graded vector spaces. But between the entries there are also maps encoding the differentials in the resolution.\footnote{A proper treatment would be to use the twisted complex of \cite{BondalKapranov}. But at this stage, this would be an unnecessary complication.} 
For example, for $Y_1$, the quiver algebra $\Gl$ is reviewed in sec.\ref{sec_TB}, and the bi-module resolution for $\Gl$ is
\bea
\begin{tikzpicture}
  \matrix (m) [matrix of math nodes, row sep=1em, column sep=1.5em]
    { \Gl e_1\otimes V^*\otimes e_0\Gl & \Gl e_0\otimes e_0\Gl \\
    & \Gl e_1\otimes e_1\Gl \\
    };
    \path[->] (m-1-1) edge (m-1-2)
    (m-1-1) edge (m-2-2);
    \node at (1.5,0) {$\oplus$};
\end{tikzpicture}\nn\eea
where $V^*=\hom(L_0,L_1)=1\oplus s$, as in remark \ref{rem_BWB}. So the matrix $\Gd_{ij}$ is
\bea
\begin{tikzpicture}
\matrix (m) [left delimiter={[},right delimiter={]},matrix of math nodes, row sep=1.5em, column sep=.8em]
    { \BB{C}  & 0 \\
       V^*[1]  & \BB{C} \\ };
      \path[->,font =\scriptsize]
      (m-2-1) edge[red] (m-1-1)
      (m-2-1) edge (m-2-2);
\end{tikzpicture}\nn\eea

Let now $F_K$ be the FM transform assigned to a tangle, say $K_m^n$. Then the kernel for $F$ is
\bea \Phi(F)= \Gl_m e_i\Gd_{ij} F(e_j\Gl_m)=\Gl_m e_i\Gd_{ij} \hom(\Gl_ne_k,F e_j\Gl_m) \Gd_{kl} e_l\Gl_n,\label{scheme_kernel_F}\eea
where we applied $\hom(\Gl_ne_k,-) \Gd_{kl} e_l\Gl_n$ to $F(e_j\Gl_m)$ to resolve it as a complex of projectives $e_l\Gl_n$.
The rhs is essentially a matrix multiplication, though one is allowed to perform chain homotopy.

Rather than doing this with all tangles, we need only the building blocks $\cap,\cup$, as all the others are derived from these.
\begin{example}
Starting with examples with two strands.
For $Y_2$ the matrix $\Gd_2$ is given in \eqref{res_Y2}, along with the description of how to read off the maps among the matrix entries.

The cup sends $\BB{C}$ to the complex
\bea \CUP\BB{C}=t^{-1/2}(\frac{t}{s}e_1\Gl \to {\color{blue}\frac{t}{s}Ve_2\Gl} \to \frac{t}{s^2}e_3\Gl),~~~\Gl:=\Gl_2,\nn\eea
with the blue term at degree 0. 
The maps in this complex are as follows. The first map goes as $ts^{-1}\to ts^{-1}V\otimes V^*$, as $V^*$ is the paths space $e_2\Gl e_1$, we multiply it to the right. The second map goes as $ts^{-1}V\cong ts^{-2}V^*$ and again multiply $V^*=e_3\Gl e_2$ to the right. This gives the kernel for $\CUP$
\bea \Phi(\CUP)=t^{-1/2}
\begin{bmatrix}
  0 & \frac{t}{s}[1] & \frac{t}{s}V & \frac{t}{s^2}[-1]
\end{bmatrix}
\begin{bmatrix}
  e_0\Gl \\ e_1\Gl \\ e_2\Gl \\ e_3\Gl
\end{bmatrix}\in K(\opn{proj}_{\Gl_2}).\nn\eea
according to the scheme of \eqref{scheme_kernel_F}.

Now for the cap, it sends $e_i\Gl_2$ to
\bea e_0\Gl\to 0,~~e_1\Gl\to s,~~e_2\Gl\to 0,~~e_3\Gl \to \frac{s^2}{t}.\nn\eea
Using the procedure introduced above, we get the kernel also for $\CAP$
\bea \Phi(\CAP)=
[\Gl_2 e_{\sbullet}]\Gd_2 [\CAP e_{\sbullet}\Gl_2]&=&
\Big[
  \Gl e_0 ~ \Gl e_1 ~ \Gl e_2 ~ \Gl e_3
\Big]
\begin{bmatrix}
  1 & 0 & 0 & 0 \\ V^*[1] & 1 & 0 & 0 \\
  s[2]\oplus\frac{s}{t}[1] & V^*[1] & 1 & 0\\
  \frac{s}{t}V^*[2] & s[2]\oplus\frac{s}{t}[1] & V^*[1] & 1
\end{bmatrix} \begin{bmatrix}
    0 \\ s \\ 0 \\ \frac{s^2}{t} 
\end{bmatrix}\nn\\
&=&\Big[
  \Gl e_0 ~ \Gl e_1 ~ \Gl e_2 ~ \Gl e_3
\Big]\begin{bmatrix}
    0 \\ s \\ sV^*[1] \\ s^2[2] 
\end{bmatrix}\in K({}_{\Gl_2}\!\opn{proj}) \nn\eea 

As we mentioned, once the FM kernel is resolved in terms of projectives, composition of kernels is just tensoring of bi-modules, and in fact, matrix multiplications
\bea \Phi(\CAP\CUP)&=&\Phi(\CUP)\otimes_{\Gl}\Phi(\CAP)\nn\\
&=&\begin{bmatrix}
  0 & \frac{t}{s}[1] & \frac{t}{s}V & \frac{t}{s^2}[-1]
\end{bmatrix}
\begin{bmatrix}
  e_0\Gl \\ e_1\Gl \\ e_2\Gl \\ e_3\Gl
\end{bmatrix}\otimes_{\Gl}\Big[
  \Gl e_0 ~ \Gl e_1 ~ \Gl e_2 ~ \Gl e_3
\Big]\begin{bmatrix}
    0 \\ s \\ sV^*[1] \\ s^2[2] 
\end{bmatrix}\nn\\
&=&\begin{bmatrix}
  0 & \frac{t}{s}[1] & \frac{t}{s}V & \frac{t}{s^2}[-1]
\end{bmatrix}
\begin{bmatrix}
  1 & \\
  V^* & 1 & \\
  S^2V^*+\frac{s}{t} & V^* & 1 &\\
  S^3V^*+\frac{s}{t}V^* & S^2V^*+\frac{s}{t} & V^* & 1
\end{bmatrix}\begin{bmatrix}
    0 \\ s \\ sV^*[1] \\ s^2[2] 
\end{bmatrix}.\nn\eea
Multiply the matrices out, and perform chain homotopies along the (unmarked, alas) maps, we get
\bea \Phi(\CAP\CUP)
&=&t^{-1/2}\begin{bmatrix}
  0 & \frac{t}{s}[1] & \frac{t}{s}V & \frac{t}{s^2}[-1]
\end{bmatrix}
\begin{bmatrix}
    0 \\ s \\ sV^*[1]+sV^* \\ s^2[2]+sV^*\otimes V^*[1]+sS^2V^*+\frac{s^2}{t}
\end{bmatrix}\nn\\
&=&t^{-1/2}\begin{bmatrix}
  0 & \frac{t}{s}[1] & \frac{t}{s}V & \frac{t}{s^2}[-1]
\end{bmatrix}
\begin{bmatrix}
    0 \\ s \\ 0 \\ \frac{s^2}{t}
\end{bmatrix}=t^{-1/2}(t[1]+[-1]).\nn\eea
So we have recovered the same result in ex.\ref{ex_Reide_I}.

We worked out in sec.\ref{sec_CS} the action of $\OC$ on the exceptional objects $\{L_{\sbullet\sbullet}\}$ of $Y_2$, see \eqref{T_1_Y_2}, as well as how it acts on maps between these objects. Applying these results, we get the kernel for $\OC$.
\bea \Phi(\OC)&=&
[\Gl_2 e_{\sbullet}]\Gd_2 [\OC e_{\sbullet}\Gl_2]\nn\\
&=&\Big[
  \Gl_2e_{\sbullet}
\Big]
\underbrace{\begin{bmatrix}
  1 & 0 & 0 & 0 \\ V^*[1] & 1 & 0 & 0 \\
  s[2]\oplus\frac{s}{t}[1] & V^*[1] & 1 & 0\\
  \frac{s}{t}V^*[2] & s[2]\oplus\frac{s}{t}[1] & V^*[1] & 1
\end{bmatrix}}_{\Gd_2}
\underbrace{\begin{bmatrix}
  1 & 0 & 0 & 0 \\
  0 & 1+t[1] & tV & \frac{t}{s}[-1] \\
  0 & 0 & 1 & 0 \\
  0 & s[1] & V^* & 1+[-1]
\end{bmatrix}}_{M(\OC)}
\begin{bmatrix}
  e_{\sbullet}\Gl_2
\end{bmatrix}.\label{see_eg}\eea
Carrying out the matrix multiplication and chain homotopy 
\bea \Phi(\OC)=\Big[
  \Gl e_0 ~ \Gl e_1 ~ \Gl e_2 ~ \Gl e_3
\Big]\!
\begin{bmatrix}
  1 & 0 & 0 & 0 \\
  V^*[1] & 1+t[1] & tV & \frac{t}{s}[-1] \\
  \wedge^2V^*[2]+\frac{s}{t}[1] & V^*([1]+t[2]) & 1+tV\otimes V^*[1] & tV \\
  \frac{s}{t}V^*[2] & \wedge^2V^*(t[3]+[2])+\frac{s}{t}[1] & V^*([1]+t[2]) & 1+t[1]
\end{bmatrix}\!
\begin{bmatrix}
  e_0\Gl \\ e_1\Gl \\ e_2\Gl \\ e_3\Gl
\end{bmatrix}\label{M_OC}\eea
In the above calculation, we have used the isomorphism $sV\cong V^*$ freely, see remark \ref{rem_BWB}.

\end{example}

In view of the above example and the clarity of using matrices to present the kernel, we:
\begin{definition}\label{def_Phi_V_M}
Define for any tangle $K_m^n$ the matrix $M(K)$, so that 
\bea \Phi(K)(e_a\Gl_m)=M(K)_{ai}e_i\Gl_n.\nn\eea
In this notation, the FM kernel for $K$ is $[\Gl_m e_{\sbullet}]\Gd_m M(K)[e_{\sbullet}\Gl_n]$ (compare with \eqref{see_eg}).
\end{definition}
We record here the matrix for the cups and caps for $Y_2,Y_3$.
\begin{proposition}
  For the cup $D^b(Y_0)\to D^b(Y_2)$ and cap $D^b(Y_2)\to D^b(Y_0)$
  \bea  M(\CUP)=t^{-1/2}\left[\begin{array}{cccc} 0 & \frac{t}{s}[1] & \frac{t}{s}V & \frac{t}{s^2}[-1] \end{array}\right],~~~
  M(\CAP)=t
  \begin{bmatrix}
    0 \\ \frac{s}{t} \\ 0 \\ \frac{s^2}{t^2} 
  \end{bmatrix}\nn\eea
  For the $\CUP$'s: $D^b(Y_1)\to D^b(Y_3)$ we have
  \bea && M({}^1\CUP)=t^{-1/2}\left[\begin{array}{cccccccc} 0 & \frac{t}{s}[1] & \frac{t}{s}V & \frac{t}{s^2}[-1] & 0 & 0 & 0 & 0\\ 0 & 0 & 0 & 0 & \frac{t^3}{s^2}[1] & 0 & \frac{t^3}{s^2}S^2V & \frac{t^3}{s^3}V[-1] \end{array}\right],\nn\\
  && M({}^2\CUP)=t^{-1}
\begin{bmatrix}
  0 & 0 & \frac{t^2}{s}V[1] & \frac{t^2}{s^2} & \frac{t^2}{s}S^3V + \frac{t^3}{s^2}V & (\frac{t^2}{s^2}S^2V+\frac{t^3}{s^3})[-1] & \frac{t^3}{s^2}S^3V[-1] & \frac{t^3}{s^3}S^2V[-2] \\
  0 & 0 & \frac{t^2}{s}[1] & 0 & \frac{t^2}{s}S^2V+\frac{t^3}{s^2} & \frac{t^2}{s^2}V[-1] & \frac{t^3}{s^2}S^2V[-1]
  & \frac{t^3}{s^3}V[-2] \end{bmatrix}.\nn\eea
 For the $\CAP$'s: $D^b(Y_3)\to D^b(Y_1)$, we have
 \bea M({}_1\CAP)=t
\begin{bmatrix}
  0 & 0 \\ \frac{s}{t} & 0 \\ 0 & 0 \\ \frac{s^2}{t^2} & 0 \\ 0 & \frac{s^2}{t^3} \\ 0 & \frac{s^2}{t^3}V^* \\
  0 & \frac{s^3}{t^4} \\ 0 & \frac{s^3}{t^4}V^*.\end{bmatrix},~~~
  M({}_2\CAP)=t^{3/2}
  \begin{bmatrix}
    0 & 0 \\ 
    0 & 0 \\ 
    0 & \frac{s}{t^2} \\
    \frac{s^2}{t^2}[1] & \frac{s^2}{t^2}V \\
    0 & 0 \\ 
    0 & 0 \\ 
    \frac{s^4}{t^4}V[1] & \frac{s^4}{t^4}S^2V \\
    \frac{s^5}{t^4}S^2V[1] & \frac{s^5}{t^4}S^3V \\
  \end{bmatrix}\nn\eea
\end{proposition}
\noindent The proof is via direct computation.
\begin{remark}
Once again, we have to omit the chain maps between entries of the matrices, but they all
go to the right whenever the entry degrees differ by 1 and the equivariance allows such a map.
For example, take the second row of $M_{{}^1\cup}$, the maps are
\bea &&L_{421}\to (S^2V)\otimes(S^2V^*)L_{421}\to (S^2V)L_{631},\nn\\
&&S^2VL_{631}\to\frac{1}{s}V\otimes V^*L_{631}\to \frac{1}{s}VL_{731}.\nn\eea
Note rather than writing $e_4\Gl_3$ or $e_6\Gl_3$, we wrote $L_{421}$ and $L_{631}$ for easier reference, since it is easier to compute $\hom_{\Gl}(e_4\Gl_3,e_6\Gl_3)$ via $\hom(L_{421},L_{631})=S^2V^*+st^{-1}$.
\end{remark}
Effectively, we reduce the composition of FM kernels as multiplying the matrices named $M$:
\bea \Phi(K_2K_1)=\Phi(K_1)\otimes_{\Gl}\Phi(K_2)=[\Gl e_{\sbullet}]\Gd M(K_1)\underbrace{[e_{\sbullet}\Gl]\otimes_{\Gl}
[\Gl e_{\sbullet}]\Gd}_{{\rm cancel}} M(K_2)[e_{\sbullet}\Gl]
=[\Gl e_{\sbullet}]\Gd M(K_1)M(K_2)[e_{\sbullet}\Gl].\nn\eea
Do notice that we multiply the matrices from the right, so $M(K_1)M(K_2)$ corresponds to $M(K_2\,K_1)$. 
\begin{example}
We provide some examples, and verify some statements of prop.\ref{prop_cap_cup} now in the homotopy category. First
\begin{tikzpicture}[scale=0.5]
\centerarc[blue](0.3,0)(-180:0:0.3);
\centerarc[blue](0.3,0)(180:0:0.3);
\end{tikzpicture} as a cup followed by a cap.
We multiply $M({}^2\CUP)M({}_2\CAP)$ to get $M(\bigcirc)$
\bea M({}^2\CUP)M({}_2\CAP)=(t^{1/2}[1]+t^{-1/2}[-1]),\nn\eea
we leave the chain homotopy to the reader.

Next, we verify straightening out a U-bend
\begin{tikzpicture}[scale=0.5]
\draw[blue] (-0.6,.6) -- (-0.6,0);
\centerarc[blue](-0.3,0)(-180:0:0.3);
\centerarc[blue](0.3,0)(180:0:0.3);
\draw[blue] (0.6,0) -- (0.6,-0.6);

\node[blue] at (1.1,-0.1) {\small{=}};

\draw[blue] (1.6,0.6) -- (1.6,-0.6);
\end{tikzpicture}~~~~~\begin{tikzpicture}[scale=0.5]
\draw[blue] (0.6,.6) -- (0.6,0);
\centerarc[blue](0.3,0)(-180:0:0.3);
\centerarc[blue](-0.3,0)(180:0:0.3);
\draw[blue] (-0.6,0) -- (-0.6,-0.6);

\node[blue] at (1.1,-0.1) {\small{=}};

\draw[blue] (1.6,0.6) -- (1.6,-0.6);
\end{tikzpicture}
\bea &&M({}^2\CUP)M({}_1\CAP)=
\begin{bmatrix}
  1 & 0 \\
  0 & \frac{s}{t}S^2V+1 + \frac{1}{t}V\otimes V^*[-1] + \frac{s}{t}S^2V[-1]+ \frac{1}{t}V\otimes V^*[-2]
\end{bmatrix}\sim 
\begin{bmatrix}
  1 & 0 \\ 0 & 1
\end{bmatrix},\nn\\
&&M({}^1\CUP)M({}_2\CAP)=
\begin{bmatrix}
  1 & (V+V[-1]) \\
  s^2S^2V\otimes (V[1]+V) & s^2(S^2V\otimes S^2V+V\otimes S^3V[-1])
\end{bmatrix}\sim 
\begin{bmatrix}
  1 & 0 \\ 0 & 1
\end{bmatrix}.\nn\eea
We hope that, even though we did not make all the maps in the matrices explicit, it is clear that a chain homotopy would simplify both to the identity matrix.
\end{example}

\begin{example}\label{ex_reef}
The final example is the reef knot widely used in sailing. 
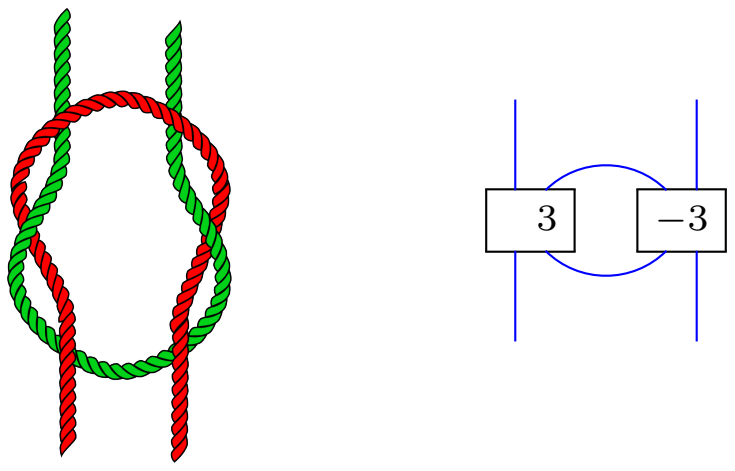
\begin{figure}[h]
\begin{center}
\includegraphics[width=3.5cm]{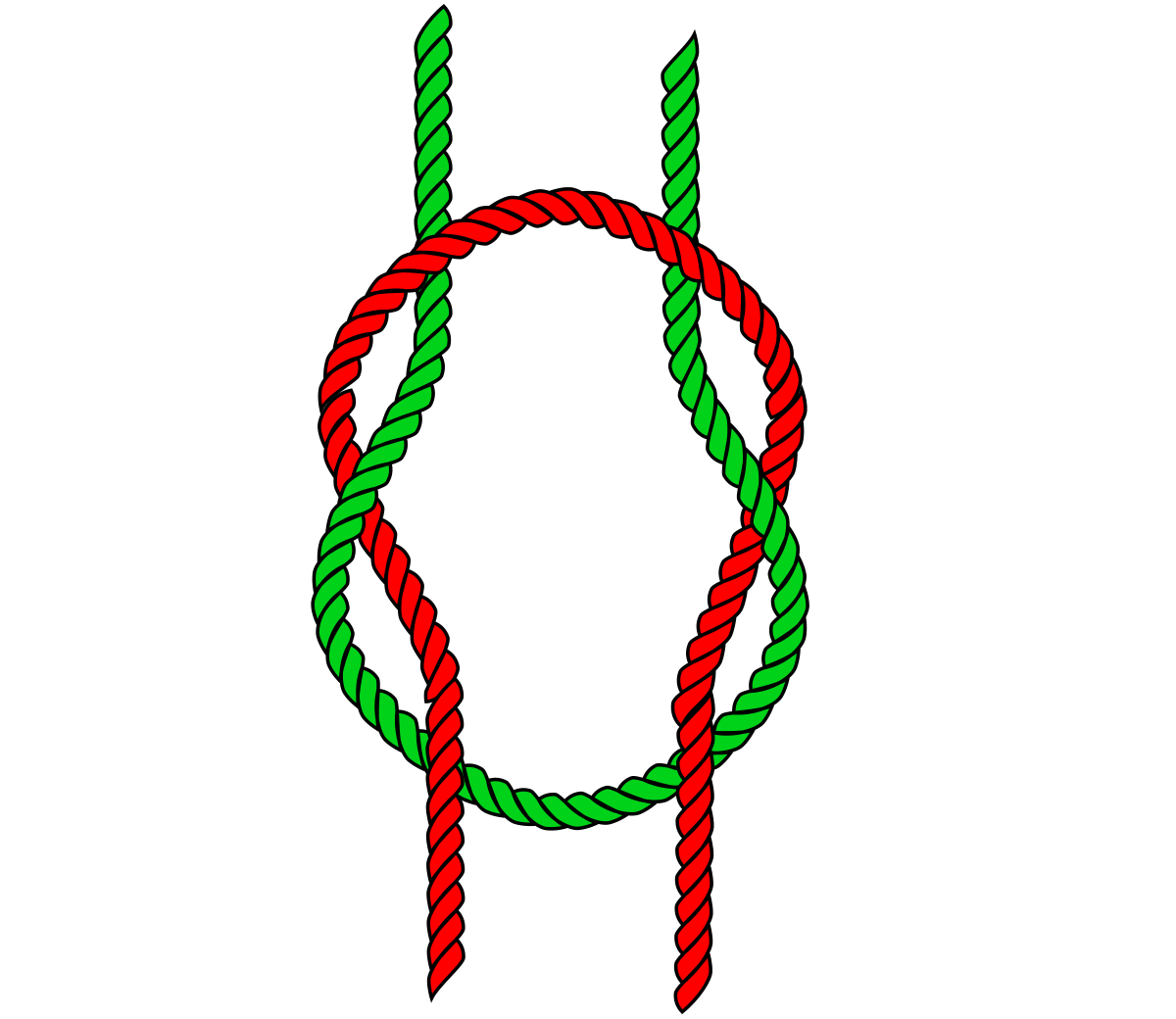}~~~
\begin{tikzpicture}[scale=1]
  \node at (0,0.5) {\scriptsize{$\boxed{\phantom{-}3}$}};
  \node at (1,0.5) {\scriptsize{$\boxed{-3}$}};
  \draw[blue] (-0.1,-0.3) -- (-0.1,0.3); \draw[blue] (1.1,-0.3) -- (1.1,0.3);
  \path[blue] (0.1,0.3) edge[bend right =45] (0.9,0.3);
  \path[blue] (0.1,0.7) edge[bend left =45] (0.9,0.7);
  \draw[blue] (-0.1,0.7) -- (-0.1,1.3); \draw[blue] (1.1,0.7) -- (1.1,1.3);
  \node at (0,-1) {};
\end{tikzpicture}\caption{The reef (also known as the square knot or Råbandsknop).}\label{fig_raband}
\end{center}
\end{figure}
Stacking up the maps as indicated in fig.\ref{fig_raband}, and testing the map on the exceptional objects, 
We can write the action of the tangle as 
\bea  
\begin{tikzpicture}
    \matrix (m) [matrix of math nodes, row sep=1em, column sep=2.5em]
    { L_{00} & & BL_{00} & \\ 
    L_{10} & A[-1]L_{10} & 
    \begin{matrix}
        BL_{10} \\ \\ A[-1]VL_{21}  
    \end{matrix} &
    \frac{1}{s}A[-1]L_{31} \\ 
    L_{21} & & BL_{21} & \\
    L_{31} & \frac{s}{t}A[-1]L_{10} & 
    \begin{matrix} 
        \frac{s}{t}A[-1]VL_{21} \\ \\ BL_{31} 
    \end{matrix} & \frac{1}{t}A[-1]L_{31} \\ };
    \node[below left=-0.6 and -0.5 of m-2-3] (x1) {\phantom{x}};
    \node[below right=-0.6 and -0.5 of m-2-3] (x2) {\phantom{x}};
    \node[above left =-0.6 and -0.85 of m-2-3] (y1) {\phantom{y}};
    \node[above=-0.6 of m-2-3] (y) {\phantom{y}};
    \node[above right =-0.6 and -0.85 of m-2-3] (y2) {\phantom{y}};
    \node[below =0.4 of y] (y3) {\phantom{$BL_{10}$}};
    \node[above left=-0.6 and -0.5 of m-4-3] (z1) {\phantom{x}};
    \node[above right=-0.6 and -0.5 of m-4-3] (z2) {\phantom{x}};
    \node[below left=-0.5 and -1 of m-4-3] (w1) {\phantom{w}};
    \node[below=-0.6 of m-4-3] (w) {\phantom{w}};
    \node[below right=-0.5 and -1 of m-4-3] (w2) {\phantom{w}};
    \node[above=0.5 of w] (w3) {\phantom{w}};
    \path[->,font =\scriptsize]
    (m-2-2) edge (y1)
    (m-2-2) edge (x1)
    (x2) edge  (m-2-4)
    (y2) edge (m-2-4)
    (y) edge (y3)
    (m-2-2) edge[bend right =30] (m-2-4);
    \path[->,font =\scriptsize]
    (m-4-2) edge (z1)
    (m-4-2) edge (w1)
    (z2) edge (m-4-4)
    (w2) edge (m-4-4)
    (w3) edge (w)
    (m-4-2) edge[bend left =30] (m-4-4);
  \node at (-3.5,0) {$\stackrel{\textrm{reef}}{\longrightarrow}$};
\end{tikzpicture}\label{Reef_Y_2}\eea
where $A,B$ are themselves complexes
\bea &&A=t^3[3]+2t^2[2]+3t[1]+2[0]+\frac{1}{t}[-1],\nn\\
&&B=t^3[2]+t^2[1]+\frac{1}{t}[-1]+\frac{1}{t^2}[-2].\nn\eea
We choose to make the similarity with \eqref{T_1_Y_2} clearer, so the left (right) column has a shift of $[1]$ ($[-1]$).
But as the entries themselves are complexes, there can be differentials from the left to the right column, e.g.  $A[-1]L_{10}\to s^{-1}A[-1]L_{31}$, which should be read as a map of degree 1 from $AL_{10}$ to $s^{-1}A[-2]L_{31}$.
The same remark applies to chain maps, say, between $\textrm{reef}(L_{00})$ and $\textrm{reef}(L_{10})$.

We can present \eqref{Reef_Y_2} in matrix form, from which we can obtain the kernel for the Reef knot (def.\ref{def_Phi_V_M})
\bea t^{1/2}M({\rm Reef})=\begin{bmatrix}
    B & 0 & 0 & 0  \\ 
      0 & A+B & A[-1]V & \frac{1}{s}A[-2] \\
    0 & 0 & B & 0 \\
   0 & \frac{s}{t}A & \frac{s}{t}A[-1]V & B+t^{-1}A[-2]  
\end{bmatrix}\nn\eea
\end{example}

\appendix
\section{More details on \texorpdfstring{$Y_n$}{Yn}}\label{Sec_DCoY_n}
In this appendix we describe the varieties $Y_n$ in some more detail, with some explicit coordinates that will be useful for describing the homological algebra of coherent sheaves on them.
\subsection{Geometry of \texorpdfstring{$Y_n$}{Yn}}\label{sec_GoY}
Given the definition \ref{def_Yn} of $Y_n$, we look at some small examples.  
\begin{example}\label{ex_Y2} We have an isomorphism $Y_2 \cong \BB{F}_2$, the Hirzebruch surface. Indeed, since we need $z{\cal L}_1=0$, so ${\cal L}_1=\textrm{span}\bra e_1,f_1\ket$, and its projectivisation $\BB{P}{\cal L}_1$ is $\BB{P}^1$. Going up to ${\cal L}_2$, since ${\cal L}_2\subset z^{-1}({\cal L}_1)$ with the latter of dim 3 containing ${\cal L}_1$. Therefore ${\cal L}_2/{\cal L}_1$ is spanned by a vector in $z^{-1}({\cal L}_1)/{\cal L}_1$ which has dim 2, this gives another $\BB{P}^1$ under projectivisation. In conclusion, we see $Y_2$ is a $\BB{P}^1$ fibred over the base $\BB{P}^1=Y_1$.

We can further get the degree of this fibration. 
We introduce another shift $\gs$
\bea \gs(e_i)=e_{i+1},~~\gs(f_i)=f_{i+1}\label{sigma}.\eea
Note $\gs$ is the right inverse of $z$: $z\gs=1$, but $\gs z\neq 1$. We reserve the notation $z^{-1}$ for the pre-image of $z$.
Now write $\ell=x_0e_1+x_1f_1$ for the parametrisation of the base $\BB{P}^1$ and ${\cal L}_2$ is included in
\bea z^{-1}({\cal L}_1)=\opn{span}\bra\gs(\ell_1),e_1,f_1\ket,\nn\eea
In the open sets $x_0\neq 0$ or $x_1\neq 0$, we can parametrise ${\cal L}_2/{\cal L}_1$ as
\bea U_0=\{x_0\neq 0\}:&&{\cal L}_2/{\cal L}_1=\{\ell_2=x_{00}\gs(\ell_1)+x_{01}x_0f_1=x_{00}(x_0e_2+x_1f_2)+x_{01}x_0f_1\},\nn\\
U_1=\{x_1\neq 0\}:&&{\cal L}_2/{\cal L}_1=\{\ell_2=x_{10}\gs(\ell_1)-x_{11}x_1e_1=x_{10}(x_0e_2+x_1f_2)-x_{11}x_1e_1\}.\nn\eea
On the intersection $x_0x_1\neq 0$, the $x_{ij}$'s are related as
\bea
x_{00}=x_{10},~~x_{01}=\frac{x_1^2}{x_0^2}x_{11}.\label{transition_Y2}\eea
This shows that $Y_2$ is the projectivisation
\bea Y_2 \cong \BB{P}({\cal O}_{\BB{P}^1}(0)\oplus{\cal O}_{\BB{P}^1}(-2)),\nn\eea
where the first (second) factor corresponds to $x_{00}$ ($x_{01}$).
This is precisely $\BB{F}_2$.
Apart from $Y_1=\BB{P}^1$, only $Y_2$ is a toric variety, with toric fans
\bea
\begin{tikzpicture}[scale=1]
\draw [step=0.3,thin,gray!40] (-0.8,-0.8) grid (0.8,0.8);
\draw [->] (0,0) -- (0.5,0);
\draw [->] (0,0) -- (0,0.5);
\draw [->] (0,0) -- (0,-0.5);
\draw [->] (0,0) -- (-0.3,-0.6);
\end{tikzpicture}\label{toric_fan_F2}\eea
The divisor $D_2$ has $z{\cal L}_2=0$, so $x_{00}=0$ and $x_{10}=0$. At $D_2$, one can verify ${\cal L}_2=z^{-1}(0)=\bra e_1,f_1\ket$. It is true in general that ${\cal L}_k=z^{-1}{\cal L}_{k-2}$ on $D_k$, which is a crucial fact to keep in mind.

The exceptional collection of $D^b(\BB{F}_2)$ was given in \cite{Hille_Perling_2011}, see also the lecture notes by \cite{Alastair_Craw}. We will get the exceptional collection for $Y_n$ in sec.\ref{sec_ECoDY}.
\end{example}

\subsection{Explicit Parametrisation of \texorpdfstring{$Y_n$}{Yn}}\label{sec_EPoY}
In a previous iteration of this paper, we had used explicit parametrisations to compute the sections of the various line bundles, but then we found a cleaner way, as was done in section \ref{sec_CS}, for instance. We keep this section nonetheless in case it is of use another day.
\begin{example}\label{ex_Y3}
From the picture \ref{mix_board}, we can visualise the parametrisation of $Y_n$.
First, the choice of ${\cal L}_k^{\perp}$ is not canonical, so we need to choose open covers for $Y_n$ and make choices for ${\cal L}_k^{\perp}$ differently. We take another look at $Y_{2,3}$.

Given $\ell_1=x_0e_0+x_1f_1$, over $U_0=\{x_0\neq 0\}$, resp. $U_1=\{x_1\neq 0\}$, the complement ${\cal L}_1^{\perp}$ can be picked as $x_0f_1$ resp. $-x_1e_1$. So over $U_0$, $\ell_2=x_{00}\gs(\ell_1)+x_{01}f_1$ and over $U_1$, $\ell_2=x_{10}\gs(\ell_1)-x_{11}x_1e_1$, giving us two more homogeneous coordinates.
These are the ones given in ex.\ref{ex_Y2}.

Going over to $Y_3$, then we have four opens $U_{ij}$, $i,j=0,1$,
\bea U_{ij}=\{x_ix_{ij}\neq0\}.\nn\eea
Over each open $\ell_3=x_{ij0}\gs(\ell_2)+x_{ij1}{\cal L}_2^{\perp}$.
The choices for ${\cal L}_2^{\perp}$ differ on each $U_{ij}$
\bea U_{00}=\{x_0x_{00}\neq 0\}:&&{\cal L}_2^{\perp}=x_{00}{\cal L}_1^{\perp}=x_{00}x_0f_1,\nn\\
U_{01}=\{x_0x_{01}\neq 0\}:&&{\cal L}_2^{\perp}=-x_{01}\gs(\ell_1),\nn\\
U_{10}=\{x_1x_{10}\neq 0\}:&&{\cal L}_2^{\perp}=x_{10}{\cal L}_1^{\perp}=-x_{10}x_1e_1\nn\\
U_{11}=\{x_1x_{11}\neq 0\}:&&{\cal L}_2^{\perp}=-x_{11}\gs(\ell_1)\nn.\nn\eea
To get some uniformity of notation, one can also set $e_1=\gs(\ell_0)$ and $f_1=\ell_0^{\perp}$ formally and so $\ell_1=x_0\gs(\ell_0)+x_1\ell_0^{\perp}$.
Now the pattern should be clear, at level 4, we have e.g. over $U_{010}=\{x_0x_{01}x_{010}\neq 0\}$, we have
\bea &&\ell_4=x_{0100}\gs(\ell_3)+x_{0101}{\cal L}^{\perp}_3,\nn\\
&&{\cal L}^{\perp}_3=x_{010}{\cal L}_2^{\perp}=-x_{010}x_{01}\gs(\ell_1).\nn\eea
\end{example}
Summarising the example, we have now an easy proposition:
\begin{proposition}\label{prop_exp_para}
  We cover $Y_n$ with $2^n$ open sets $U_{i_1,\cdots,i_n}$, with $i_k=0,1$.
  \bea U_{i_1i_2,\cdots,i_n}=\{x_{i_1}x_{i_1i_2}\cdots x_{i_1\cdots i_n}\neq 0\}.\nn\eea
  The parametrisation of a line $\ell_n\in {\cal L}_n/{\cal L}_{n-1}$ over $U_{i_1,\cdots,i_{n-1}}\subset Y_{n-1}$ is
  \bea \ell_n=x_{i_1i_2,\cdots,i_{n-1}0}\gs(\ell_{n-1})-x_{i_1i_2,\cdots,i_{n-1}1}x_{i_1i_2,\cdots,i_{n-1}}\cdots x_{i_1i_2,\cdots,i_k}\gs(\ell_{k-1}),\nn\eea
  where $i_k$ is the rightmost index that is 1. In case $k=1$ and $i_1=1$, then $\gs(\ell_0)$ is taken to be $e_1$. While if no $i_k$ is 1, we take $\gs(\ell_{-1})$ to be $-f_1$.
\end{proposition}
\begin{proof}
 Recall ${\cal L}_n$ is spanned by $a\gs(\ell_{n-1})+b{\cal L}_{n-1}^{\perp}$ as well as $\ell_{n-1},\cdots,\ell_1$, and $z^{-1}({\cal L}_{n-1})$ is spanned by $\gs(\ell_{n-1}),\cdots,\gs(\ell_1)$ and $z^{-1}(0)$.
 If $a\neq 0$, we can represent ${\cal L}_n^{\perp}$ with $a{\cal L}^{\perp}_{n-1}$. If $b\neq 0$, then $-b\gs(\ell_{n-1})$ can be chosen as the representative.

 Now, over $U_{i_1,\cdots,_{n-1}}$, we write
 \bea \ell_n=x_{i_1,\cdots,i_{n-1},0}\gs(\ell_{n-1})+x_{i_1,\cdots,i_{n-1},1}{\cal L}_{n-1}^{\perp}.\nn\eea
 To parametrise ${\cal L}_{n-1}^{\perp}$, we need to differentiate between $i_{n-1}=1$ or not. If yes, then 
 \[{\cal L}_{n-1}^{\perp}=-x_{i_1,\cdots,i_{n-1}}\gs(\ell_{n-2}),\] 
 if not we need to keep writing out $\ell_{n-2}$ and inspect if $i_{n-2}=1$ and so on.
\end{proof}

With the explicit parametrisations, we can get the transition functions between two opens. We gave the transition for $Y_2$ in \eqref{transition_Y2}, now we go to the next case $Y_3$:
\begin{example}\label{ex_transition_Y3}
  Continuing from Ex.\ref{ex_Y3}, we have the transition functions
  \bea U_{i0}\cap U_{i1}:&&x_{i10}=x_{i00},~~x_{i11}=\frac{x_{i0}^2}{x_{i1}^2}x_{i01},\nn\\
  U_{00}\cap U_{10}:&&x_{100}=x_{000},~~~x_{101}=\frac{x_0^2}{x_1^2}x_{001}-\frac{x_{11}^2}{x_{10}^2}\frac{x_1}{x_0}x_{100},\nn\\
  U_{01}\cap U_{11}:&&x_{110}=x_{010},~~~
  x_{111}=\frac{x_1^2}{x_0^2}x_{011}-\frac{x_1}{x_0}x_{110}.\nn\eea
Over $U_{00}\cap U_{10}$, we have two parametrisations of $\ell_3$
\bea U_{00}:&&\ell_3=x_{000}\gs(\ell_2)+x_{001}x_{00}x_0f_1=x_{000}\gs(x_{00}\gs(\ell_1)+x_{01}x_0f_1)+x_{001}x_{00}x_0f_1\nn\\
U_{10}:&&\ell_3'=x_{100}\gs(\ell'_2)-x_{101}x_{10}x_1e_1=x_{100}\gs(x_{10}\gs(\ell_1)-x_{11}x_1e_1)-x_{101}x_{10}x_1e_1.\nn\eea
These two need only agree up to ${\cal L}_2$. We conclude $x_{000}=x_{100}$ (recall that $x_{00}=x_{10}\neq 0$, $x_{01}=(x_1^2/x_0^2)x_{11}$ on $U_{00}\cap U_{10}$) and
\bea && x_{000}\gs(x_{01}x_0f_1+x_{11}x_1e_1)+x_{001}x_{00}x_0f_1+x_{101}x_{10}x_1e_1=0\;\opn{mod}\;{\cal L}_2\nn\\
\To && x_{000}\frac{x_{11}}{x_{00}}\frac{x_1}{x_0}\gs(x_1f_1+x_0e_1)+x_{001}x_0f_1+x_{101}x_1e_1=0\;\opn{mod}\;{\cal L}_2\nn\\
\To && x_{000}\frac{x^2_{11}}{x^2_{00}}\frac{x_1}{x_0}x_1e_1+x_{001}x_0f_1+x_{101}x_1e_1=0\;\opn{mod}\;{\cal L}_1,\nn\eea
which leads to the relation above and also that $\ell'_3$ and $\ell_3$ differ by multiples of $\ell_1,\ell_2$.

Over $U_{01}\cap U_{11}$, we have $x_{00}=x_{10}$ (but they may be zero), $x_{01}=(x_1^2/x_0^2)x_{11}\neq 0$ and
\bea U_{01}:&&\ell_3=x_{010}\gs(\ell_2)-x_{011}x_{01}\gs(\ell_1)=x_{010}\gs(x_{00}\gs(\ell_1)+x_{01}x_0f_1)-x_{011}x_{01}\gs(\ell_1)\nn\\
U_{11}:&&\ell'_3=x_{110}\gs(\ell'_2)-x_{111}x_{11}\gs(\ell_1)=x_{110}\gs(x_{10}\gs(\ell_1)-x_{11}x_1e_1)-x_{111}x_{11}\gs(\ell_1).\nn\eea
Note that $\ell_2-\ell'_2=x_{01}(x_0/x_1)\ell_1$, so
\bea \ell_3-\ell'_3=(x_{010}-x_{110})\gs(\ell_2)+x_{110}x_{01}\frac{x_0}{x_1}\gs(\ell_1)+(x_{111}x_{11}-x_{011}x_{01})\gs(\ell_1),\nn\eea
which is zero on the nose, rather than up to ${\cal L}_2$, and we get the stated relations.
\end{example}

\begin{example} As we defined $D_1$ to be ${\cal L}_1=\{f_1\}$, so that for all $D_k$, the place function is given by $x_{i_1,\cdots,i_{k-1},0}$.
For example, for $Y_2=\BB{F}_2$, we have $D_1=\{x_0=0\}$ and $D_2=\{x_{i0}=0\}$. These divisors are toric: $D_1$ corresponds to the ray pointing southwest in \eqref{toric_fan_F2}, while $D_2$ is the south pointing ray.
For other $Y_n$, they are no longer toric.
\end{example}

\subsection{\texorpdfstring{$\BB{C}^*$}{Cstar}-Actions}\label{sec_equivariance}
\begin{definition}($\BB{C}^*$-Equivariance)
 There are two $\BB{C}^*$ actions on $Y_n$: the principle and secondary $\BB{C}^*$. 
 The principle one acts on $e_i,f_i$ with weight $-i$. The secondary one acts only on $f_i$ with weights $-1$.
\end{definition}
The reason for the assignment of the weight is: since $Y_n$ is defined using the action $z$ that shifts $e_i,f_i$ down to $e_{i-1},f_{i-1}$. To make $z$ equivariant, the weight difference between $e_i,e_{i+1}$ must be equal to that between $f_i,f_{i+1}$. So the possible weights are: $[p,p+n,p+2n,\cdots]$ for $[e_1,e_2,\cdots]$ and $[q,q+n,q+2n,\cdots]$ for $[f_1,f_2,\cdots]$. But an overall shift of weights of all $e_i,f_i$ is trivial, so we are left with two possibilities.
We pick the weights of the principal action to be $p=q=n=-1$, and the secondary one to be $p=n=0$, $q=-1$.
\begin{lemma}\label{lem_weights}
 Let the weights of the principle and secondary $\BB{C}^*$ action be $w_t$ and $w_s$, then we have
 \bea &&w_t(x_{i_1\cdots,i_{n-1},0})=1,~~~~
 w_t(x_{i_1\cdots,i_{n-1},1})=\left\{\begin{array}{cc}-w_t(x_{i_1\cdots,i_{k-1},1})+2-n+k & k>0 \\
           2-n                            & k=0
                                        \end{array}\right.\\
&& w_s(x_{i_1\cdots,i_{n-1},0})=0,~~~~
 w_s(x_{i_1\cdots,i_{n-1},1})=\left\{\begin{array}{cc}-w_t(x_{i_1\cdots,i_{k-1},1}) & k>0 \\
           1                            & k=0
                                        \end{array}\right..\nn\eea
where $i_k$ is the right most index among $\{i_1,\cdots, i_{n-1}\}$ that is 1.

In particular
\bea w_t(\frac{x_{[n-1]0}}{x_{[n-1]1}})=1-(n\opn{mod}2),~~~w_s(\frac{x_{[n-1]0}}{x_{[n-1]1}})=(-1)^n.\nn\eea
\end{lemma}
\begin{proof}
 For $w_t$, the $e_i,f_i$ have weights $-i$. In $\ell_1=x_0e_1+x_1f_1$, to keep $\ell_1$ invariant, we must have $w_t(x_i)=1$.

 The rest follows from the explicit parametrisation given in prop.\ref{prop_exp_para}. Note $\gs(\ell_p)$ has principal weight $-1$ always, hence $x_{i_1\cdots,i_{n-1},0}$ is given weight 1. Further the sum of weights of $x_{i_1\cdots,i_{n-1},1}$, $x_{i_1\cdots,i_{n-2},0}$, ..., $x_{i_1\cdots,i_k,0}$ and $x_{i_1\cdots,i_{k-1},1}$ must be 1, giving our result.
 
 For the secondary action, only $f_i$ have weight $-1$, and so $\gs(\ell_p)$ has weight 0, and the rest of the argument is the same.
\end{proof}
\begin{example}\label{ex_weights}
 We list the weights for coordinates of $Y_3$,
 \bea \begin{array}{|ccc|ccc|}
 \hline
    & w_t & w_s & & w_t & w_s \\
    \hline
      x_0 & 1 & 0 & x_1 & 1 & 1 \\ \hline
      x_{00} & 1 & 0 & x_{01} & 0 & 1\\ \hline
      x_{10} & 1 & 0 & x_{11} & 0 & -1 \\ \hline
      \end{array},~~~\begin{array}{|ccc|ccc|}
 \hline
      x_{000} & 1 & 0 & x_{001} & -1 & 1 \\ \hline
      x_{010} & 1 & 0 & x_{011} & 1 & -1 \\ \hline
      x_{100} & 1 & 0 & x_{101} & -1 & -1 \\ \hline
      x_{110} & 1 & 0 & x_{111} & 1 & 1 \\
      \hline
 \end{array}.\nn\eea
      One can verify that the transitions in Ex.\ref{ex_transition_Y3} preserve the weights.
\end{example}

Finally we prove lemma \ref{lem_trivial_lb}, which is the load bearing lemma for all the computations of sections, pushdowns or pullbacks.
We also encourage the reader to use the explicit parametrisation for small $Y_n$ in sec.\ref{sec_EPoY} to explicitly verify the claims of the lemma.

\begin{proof} of lemma \ref{lem_trivial_lb}.

First recall that the divisors $D_k$ were defined in \eqref{divisor_D} as the locus $\{z{\cal L}_k\subset {\cal L}_{k-2}\}$.
\begin{enumerate}
  \item The correspondence between line bundles and divisors is standard: one takes a rational section of the line bundle, its zeros and poles give the divisor. Our setting is equivariant, so the section has to be equivariant too. 
      First for $Y_1$, choose $\ell_1$ so that its $e_1$ component is 1. This section has a pole at the point ${\cal L}_1=\opn{span}f_1$, i.e. we get ${\cal O}(-D_1)$. Finally we need to take into account the weight of the place function, which is $t^{-1}$ for $1/x_0$. 
      
      Next we construct sections for ${\cal E}_k={\cal L}_k/{\cal L}_{k-1}$ recursively. 
      For $\ell_2\in {\cal L}_2/{\cal L}_1$, we demand $z(\ell_2)=\ell_1$, which surely has a pole at ${\cal L}_2/{\cal L}_1={\cal L}^{\perp}_1$ (defined in \eqref{def_perp}), that is, at $D_2$. Thus this section has a pole at $D_1$ and $D_2$, and we get ${\cal E}_2\cong t^{-2}{\cal O}(-D_1-D_2)$, where the $t^{-2}$ again comes from the weights of the place functions $1/x_0,1/x_{00}$. One may continue and define a section of $\ell_k\in {\cal L}_k/{\cal L}_{k-1}$ such that $z(\ell_k)=\ell_{k-1}$ and so on.
  \item %

Given a point $\{x_1,\cdots,x_{n-2}\}\in Y_{n-2}$, consider its pre-image under $q_{k+1}$ or $q_{k+2}$, i.e. a lift
\bea \{x_1,\cdots x_{k-1},x_k,x,z^{-1}x_k,z^{-1}x_{k+1},\cdots,z^{-1}x_{n-2}\}\in D_{k+1}\cap D_{k+2}\subset Y_n\label{preimg_q}\eea
where $x$ is such that $zx\in\opn{span}\bra x_1,\cdots,x_{k-1}\ket$. As already pointed out in the proof of lem.\ref{lem_fibration}, the span of $x_1,\cdots,x_k,x$ equals $z^{-1}\opn{span}\bra x_1,\cdots,x_{k-1}\ket$, hence contains the kernel of $z$. This shows that the resulting flag is unambiguous and lives in $Y_n$ (in fact in $D_{k+1}\cap D_{k+2}$). Conversely, take any $\{y_1,\cdots,y_k,\udl{y_{k+1}},\udl{y_{k+2}},y_{k+3},\cdots\}\in D_{k+1}\cap D_{k+2}$. By definition $y_{k+2}$ is in the complement of $\opn{span}\bra y_1,\cdots,y_{k+1}\ket$ in $z^{-1}\opn{span}\bra y_1,\cdots,y_k\ket=\opn{span}\bra y_1,\cdots,y_{k+1},z^{-1}y_k\ket$, so may be written as $\gs y_k$, i.e. uniquely determined from $y_k$. Its image under $q_{k+1}$ is \bea \{y_1,\cdots,y_{k-1},zy_{k+2},zy_{k+3},\cdots\}=\{y_1,\cdots,y_{k-1},y_k,zy_{k+3},\cdots\},\nn\eea
which lifts according to \eqref{preimg_q} to $\{y_1,\cdots,y_{k-1},y_k,y,y_{k+2},\cdots\}$, where $y$ must be $y_{k+1}$. The argument for $q_{k+2}$ is similar. 

In contrast, assuming $l\geq k+2$, we want to show that
\bea 
\begin{tikzpicture}
  \matrix (m) [matrix of math nodes, row sep=2em, column sep=3em]
    { D_{k+1}\cap D_{l+1} & D_{k+1} \\
      D_{l-1} & Y_{n-2} \\ };
  \path[->, font=\scriptsize]
  (m-1-1) edge node[above] {$q_{l+1}$} (m-1-2)
  (m-1-1) edge node[left] {$q_{k+1}$} (m-2-1)
  (m-1-2) edge node[right] {$q_{k+1}$} (m-2-2)
  (m-2-1) edge node[above] {$q_{l-1}$} (m-2-2);
\end{tikzpicture}\nn\eea
is a pullback. To head off confusion, the bottom left and top right are in $Y_n$, while the top left is considered as in $Y_{n+2}$.

First let $\{x_1,\cdots,x_k,\udl{x_{k+1}},\cdots\}$ and $\{y_1,\cdots,y_{l-2},\udl{y_{l-1}},\cdots\}$ parametrise a point in $D_{k+1}$ and $D_{l-1}$ respectively. They project to the same point in $Y_{n-2}$:
\bea \{x_1,\cdots,x_{k-1};zx_{k+2},\cdots zx_{l-1};zx_l,\cdots,zx_n\}\sim\{y_1,\cdots,y_{k-1};y_k,\cdots,y_{l-3};zy_l,\cdots,zy_n\}\label{same_flag}\eea
where by $\sim$, we mean they represent the same flag. We claim that a point in the fibre product can be parametrised as a flag
\bea \{x_1,\cdots,x_{k-1};x_k,\udl{x_{k+1}},x_{k+2},\cdots,x_{l-1},z^{-1}y_{l-2},\udl{z^{-1}y_{l-1}}, z^{-1}y_l,\cdots,z^{-1}y_n\}.\label{temporary_param}\eea
Remember that $\opn{span}\bra x_1,\cdots, x_{k+1}\ket=z^{-1}\opn{span}\bra x_1,\cdots, x_{k-1}\ket$ contains the kernel of $z$ due to the special position $x_{k+1}$ is in, so as far as the flag is concerned, there is no ambiguity in writing $z^{-1}(\cdots)$ in \eqref{temporary_param}. We check that the flag above lives in $Y_{n+2}$. In fact for $i=l-2,\cdots,n$, we have 
\bea zz^{-1}y_i&=&y_i\in z^{-1}\opn{span}\bra y_1,\cdots,y_{i-1}\ket
=z^{-1}\opn{span}\bra x_1,\cdots, x_{k-1},zx_{k+2},\cdots,zx_{l-1},y_{l-2},\cdots,y_{i-1}\ket\nn\\
&=&\opn{span}\bra x_1,\cdots, x_{k-1},x_k,\udl{x_{k+1}},x_{k+2},\cdots,x_{l-1},z^{-1}y_{l-2},\cdots,z^{-1}y_{i-1}\ket\nn\eea
and so the incidence condition \eqref{eq_def_Yn} holds. Also note
\bea zz^{-1}y_{l-1}&=&y_{l-1}\in z^{-1}\opn{span}\bra y_1,\cdots,y_{l-3}\ket\stackrel{\eqref{same_flag}}{=}z^{-1}\opn{span}\bra x_1,\cdots,x_{k-1},zx_{k+2},\cdots,zx_{l-1}\ket\nn\\
&=&\opn{span}\bra x_1,\cdots,x_{k-1},x_k,\underline{x_{k+1}},x_{k+2},\cdots,x_{l-1}\ket,\nn\eea 
showing that \eqref{temporary_param} lives in the intersection $D_{k+1}\cap D_{l+1}$.

Now we need to check that $D_{k+1}\cap D_{l+1}$ maps under $q_{k+1},q_{l+1}$ to $D_{k+1}$ and $D_{l-1}$ as expected, which would show that $D_{k+1}\cap D_{l+1}$ is the fibre product.
Its image under $q_{k+1}$ reads
\bea \{x_1,\cdots,x_{k-1};y_k,\cdots,y_{l-3},y_{l-2},\udl{y_{l-1}},y_l,\cdots,y_n\},\nn\eea
but we may replace the $x$'s with the corresponding $y$'s due to \eqref{same_flag}.
Its images under $q_{l+1}$ reads 
\bea \{x_1,\cdots,x_{k-1};x_k,\udl{x_{k+1}},x_{k+2},\cdots,x_{l-1},y_l,\cdots,y_n\},\nn\eea
where we can also replace the $y$'s with $x$'s, again due to \eqref{same_flag} and that $\bra x_1,\cdots,x_{k-1};x_k,\udl{x_{k+1}}\ket$ contains the kernel of $z$.

\item 
We consider the line bundle with divisor $t^2{\cal O}(D_{k+1}+D_{k+2})\cong{\cal E}^{\vee}_{k+2}\otimes{\cal E}_k$ restricted to $D_{k+1}\cap D_{k+2}$. %
To show that it is trivial, we need a nowhere vanishing section for $\hom({\cal E}_{k+2}|_{D_{k+1}\cap D_{k+2}},{\cal E}_k|_{D_{k+1}\cap D_{k+2}})$. The map $z$ will do exactly this. First $z$ does map ${\cal L}_{k+2}$ to ${\cal L}_k$ at $D_{k+2}$. When further restricted to $D_{k+1}$, it sends ${\cal L}_{k+1}$ to ${\cal L}_{k-1}$, and hence it induces a well defined map ${\cal E}_{k+2}\to {\cal E}_k$. It is also never zero since, as we remarked earlier, ${\cal L}_{k+2}$ will not be mapped to ${\cal L}_{k-1}$. 
Finally, the $\BB{C}^*$-weight of this section is that of $z$ i.e. $t$, we get $t^2{\cal O}(D_{k+1}+D_{k+2})|_{D_{k+1}\cap D_{k+2}}\cong t{\cal O}_{D_{k+1}\cap D_{k+2}}$.

\item If $\iota:\;D_{k+1}\to Y_n$ is the inclusion, then $\iota^*{\cal E}_{k+1}$ is the line ${\cal L}_k^{\perp}$. From the short exact sequence $0\to {\cal E}_k\to z^{-1}({\cal L}_{k-1})/{\cal L}_{k-1}\to {\cal L}_k^{\perp}\to0$, we have an isomorphism of line bundles ${\cal E}_k\otimes {\cal L}_k^{\perp}\cong \det(z^{-1}({\cal L}_{k-1})/{\cal L}_{k-1})=\gs({\cal E}_{k-1})\otimes {\cal L}_{k-1}^{\perp}\cong t^{-1}{\cal E}_{k-1}\otimes {\cal L}_{k-1}^{\perp}$. Now we can continue the isomorphism until we get $t^{1-k}e_1\otimes f_1\cong t^{-k-1}s^{-1}{\cal O}_{D_{k+1}}$. Finally, we use item 1 and get
    \bea t^{-2k-1}\iota^*{\cal O}(-D_1-\cdots-D_{k+1})\otimes{\cal O}(-D_1-\cdots-D_k)\cong \iota^*{\cal E}_{k+1}\otimes {\cal E}_k\cong t^{-k-1}s^{-1}{\cal O}_{D_{k+1}},\nn\eea
    which is the result to be proved. Ex.\ref{ex_sec_Y2} illustrates this point by exhibiting the transition function.
\end{enumerate}
\end{proof}

\section{Homological Algebra on \texorpdfstring{$D^b(Y_n)$}{the bounded derived category of Yn}}
We now describe the homological algebra of the derived category of the varieties $Y_n$, with a view to the use of tilting theory to convert geometric calculations into algebraic ones.
\subsection{Exceptional Collection of \texorpdfstring{$D^b(Y_n)$}{the bounded derived category of Yn}}\label{sec_ECoDY} 
\begin{definition}\label{def_exp_collection}
  Over $Y_n$, we set
  \bea L_{a_1,\cdots,a_n}:={\cal O}_{Y_n}(a_1D_1+\cdots+a_nD_n).\nn\eea
  We define a collection of $2^n$ line bundles $F_a$ with $a=0,\cdots,2^n-1$ as
  \bea F_a={\cal O}_{Y_n}(a_1D_1+\cdots+a_nD_n),~~~a_1=a,~a_2=\lfloor a/2\rfloor,~a_3=\lfloor a/4\rfloor,\cdots,a_n=\lfloor a/2^{n-1}\rfloor,\nn\eea
  that is $a_k=2a_{k+1}$ or $2a_{k+1}+1$.
\end{definition}
\begin{example}
  For $Y_2=\BB{F}_2$, the collection is
  \bea \bra L_{0,0},L_{1,0},L_{2,1},L_{3,1}\ket,\nn\eea
  while for $Y_3$ it is
  \bea \bra L_{0,0,0},L_{1,0,0},L_{2,1,0},L_{3,1,0},L_{4,2,1},L_{5,2,1},L_{6,3,1},L_{7,3,1}\ket.\nn\eea
\end{example}
\begin{theorem}\label{thm_exp_collection_Yn}
  The collection in def.\ref{def_exp_collection} is a full exceptional collection of $D^b(Y_n)$.
\end{theorem}
\begin{proof}
  We first prove 
  \bea \ext^{\sbullet}(F_a,F_b)=0,~~\textrm{for}~ a>b.\nn\eea
  For the computation of derived homs between line bundles, we will use the Leray spectral sequence in its simplest form (see the proof of Theorem 2.8 in \cite{Costa2011DerivedCO}).

  As the $F_a$'s are line bundles, computing $\ext(F_n,F_m)$ amounts to computing $R\pi_*(F_n^{\vee}\otimes F_m)$, where $\pi:\,Y_n\to *$. But pushing down to a point can be done in stages: we denote with $\pi_k:\, Y_k\to Y_{k-1}$ and so $R\pi_*=R\pi_{1*}\cdots R\pi_{n*}$.

  Let $L$ be a line bundle, since $\pi_n$ is a fibration, we get that $R\pi_{n*}L$ is a vector bundle over $Y_{n-1}$. The subsequent pushdowns have the same property. Besides, each pushdown is along a $\BB{P}^1$, over which a vector bundle splits as line bundles. Thus we will get $R\pi_*(L)=0$ if at any stage of the pushdown, we encounter only line bundles of degree $-1$, since ${\cal O}(-1)$ has no cohomology over $\BB{P}^1$. This is the basic strategy for showing the vanishing of Ext.
  
  If $a>b$, then $F_a^{\vee}\otimes F_b$ is a line bundle of the form ${\cal O}(\cdots-D_k)$ with $k$ being the right most digit of $a$ (as a binary number) that is bigger than $b$. Then the pushdown $\pi_{k+1,*}\cdots\pi_{n,*}$ does nothing and $\pi_{k*}$ gives zero. This proves the vanishing of $\ext^{\sbullet}(F_a,F_b)$, for $a>b$.
  
  \smallskip
  
  To state that the collection is exceptional, we need further prove
  \bea \ext^{>0}(F_a,F_b)=0,~~\textrm{for}~ a\leq b,~~~~\ext^0(F_a,F_a)=\BB{C}.\nn\eea
  The second condition is trivial, but the first needs more work. We will prove in prop.\ref{prop_sec_L} that the line bundle $L_{a_1\cdots a_n}$ has no higher sheaf cohomology if $a_k\geq 2a_{k+1}-1$. 
  The reader may first look at the explicit calculation of $Y_2$ in Ex.\ref{ex_sec_Y2} to see the significance of our choice $F_a$. Now for the $\hom^{\sbullet}(F_a,F_b)$. Writing $L_{\vec a}$ for $F_a$, the hom is $H^{\sbullet}(Y_n,L_{\vec b-\vec a})$. We have $(b_k-a_k)-2(b_{k+1}-a_{k+1})=(b_k-2b_{k+1})-(a_k-2a_{k+1})\geq-1$. So the higher cohomologies are absent.
  
  Finally, the fullness of the collection is proved by Orlov, as Theorem~2.6 in \cite{Orlov_1993}. Translated into our setting, Orlov's result says that $\bra \pi_n^*D^b(Y_{n-1})\otimes {\cal E}_n,\pi_n^*D^b(Y_{n-1})\ket$ generates $D^b(Y_n)$, where ${\cal E}_n$ is the canonical line bundle ($\pi_n^*$ actually embeds $D^b(Y_{n-1})$ as a fully faithful subcategory). One is free to tensor ${\cal E}_n$ with a line bundle on $Y_{n-1}$, thus $\bra \pi^*D^b(Y_{n-1}),\pi^*D^b(Y_{n-1})\otimes F_{2^{n-1}+1}\ket$ also generates $D^b(Y_n)$. An induction on $n$ now finishes the proof.
\end{proof}
\begin{example}\label{ex_sec_Y2}
 For $Y_2$ we will the sheaf cohomology for ${\cal O}(pD_1+qD_2)$, $p,q\geq0$, using the explicit parametrisation of $Y_n$. We use the stepwise pushdown used in the proof of the last proposition. Working in the open $U_0$, the pushdown along $\BB{P}^1_2$ gives the sections
 \bea \pi_{2*}{\cal O}(pD_1+qD_2)=\opn{span}\bra \big(\frac{x_{01}}{x_{00}}\big)^m,~m=0,\cdots,q\ket,\nn\eea
 which is a vector bundle over $Y_1$. On the intersection $U_0\cap U_1$, its transition function is
 \bea \big(\frac{x_{01}}{x_{00}}\big)^m= \big(\frac{x_{11}}{x_{10}}\big)^m\big(\frac{x_1}{x_0}\big)^{2m}.\nn\eea
 Thus each term gives a line bundle of degree $-2m$, therefore if $p-2m\geq-1$, the next pushdown $\pi_{1*}$ will not produce higher cohomology. The assumption $p\geq 2q-1$ would then guarantee the absence of higher cohomology. The sections are
 \bea \Gc(Y_2,{\cal O}(pD_1+qD_2))=\opn{span}\bra \big(\frac{x_{01}}{x_{00}}\big)^m\big(\frac{x_1}{x_0}\big)^n,~m=0,\cdots,q, ~~n=0,\cdots,p-2m\ket,~~p\geq 2q-1.\nn\eea

In contrast, we investigate ${\cal O}(-2D_2)$. The first pushdown gives $\pi_{2*}L=(x_{01}/x_{00})^{-1}[-1]$ i.e. at degree 1. This is a line bundle of degree $+2$ over $Y_1$, so we get
\bea H^1(Y_2,{\cal O}(-2D_2))=\opn{span}\bra \frac{x_{00}}{x_{01}}\big(\frac{x_1}{x_0}\big)^n,~~n=0,1,2\ket.\nn\eea
To summarise, if $p\geq 2q-1$, the line bundle ${\cal O}(pD_1+qD_2)$ has no higher cohomology. This motivated our choice of the collection in def.\ref{def_exp_collection}.
\end{example}

Before finishing the proof of Thm.\ref{thm_exp_collection_Yn}, we give a nice way to compute the cohomology of $L_{a_1,a_2,\cdots}$, without using the explicit parametrisation in the last example.
\begin{example}\label{ex_L_sheaf_coh}
  Let $L=L_{a_1,\cdots,a_n}$, we compute its sheaf cohomology using repeated pushdown $Y_n\stackrel{\pi_n}{\to}Y_{n-1}\cdots\stackrel{\pi_1}{\to}*$. If $a_n=0$ then $\pi_{n*}L_{a_1,\cdots,a_n}=L_{a_1,\cdots,a_{n-1}}$. If $a_n>0$, then we use the exact sequence 
  \bea 0\to F(-D)\to F \to F|_{D}\to 0\label{quotient}\eea
  to lower $a_n$ until 0, while if $a_n<0$, we use the same sequence to raise it under we reach $-1$, at which point the pushdown $\pi_{n*}$ gives zero. 
  
  Consider $L_{01}$, then \eqref{quotient} gives (use item 4 lem.\ref{lem_trivial_lb} in the second step)
  \bea \pi_{2*}L_{01}&=&\pi_{2*}\opn{cone}(L_{01}|_{D_2}[-1]\to L_{00})
  =\pi_{2*}\opn{cone}(\frac{s}{t}L_{-20}|_{D_2}[-1]\to L_{00})=\opn{cone}(\frac{s}{t}L_{-2}[-1]\to L_{0}),\nn\\
  \pi_{1*}\pi_{2*}L_{01}&=&\opn{cone}(\frac{s}{t}\pi_{1*}L_{-2}[-1]\to \pi_{1*}L_{0}).\nn\eea
  We have $\pi_{1*}L_0=1$, while
  \bea \pi_{1*}L_{-2}=\pi_{1*}\opn{cone}(L_{-1}|_{D_1}\to L_{-1})[-1]=\opn{cone}(s^{-1}\to 0)[-1]=s^{-1}[-1].\nn\eea
  So in total we get
  \bea H^{\sbullet}(Y_2,L_{01})=t^{-1}[-1]+[0].\nn\eea
\end{example}
The proof of the next proposition is based on the same algorithm
\begin{proposition}\label{prop_sec_L}
Assume the sequence satisfies $a_1\geq 2a_2-1,\cdots, a_{n-1}\geq 2a_n-1$. Then the line bundle $L_{a_1,\cdots,a_n}$ has no higher sheaf cohomology while $H^0(Y_n,L_{a_1,\cdots,a_n})$ has the following basis
  \bea &&\Gc(Y_n, L_{a_1,\cdots,a_n})\cong \bigoplus_{b_n=0}^{a_n}\bigoplus_{b_{n-1} =2c_{n-1}}^{a_{n-1}}\cdots\bigoplus_{b_2=2c_2}^{a_2}\bigoplus_{b_1=2c_1}^{a_1}
  L_{\udl{b_1},\cdots,\udl{b_n}},\label{iso_vec_space}\\
  &&\hspace{1.55cm}c_k=b_{k+1}-b_{k+2}+b_{k+3}-\cdots\pm b_n\nn\eea
where $L_{\udl{b_1},\cdots,\udl{b_n}}$ is the complex line obtained from restricting $L_{b_1,\cdots,b_n}$ to the point $D_1\cap \cdots \cap D_n$.
\end{proposition}
\begin{proof}
  Consider the exact sequence
  \bea 0\to F(-D)\to F \to F|_{D}\to 0.\eea
  If the two side terms have no higher cohomology, nor will the middle one. 
  Setting $F=L_{a_1,\cdots,a_n}$ and $D=D_n$, then $F(-D)=L_{a_1,\cdots,a_n-1}$ is covered by an induction on the $a_i$'s. The right sheaf is $L_{a_1,\cdots,a_n}|_{D_n}$. We first give it an easier notation\footnote{
  That is, $L_{a_1,\cdots,\udl{a_n}}=\iota_*\iota^*L_{a_1,\cdots,a_n}$, with $\iota:\,D_n\hookrightarrow Y_n$ being the closed embedding. We also allow ourselves the latitude to write $L_{a_1,\cdots,\udl{a_n}}$ as the restriction $\iota^*L_{a_1,\cdots,a_n}$. The reason is 1. $\iota_*$ does not do anything nor does it require deriving, and 2. sheaf cohomology of $\iota_*F$ over $Y_n$ is computed as that of $F$ on $D_n$.}
  \bea L_{a_1,\cdots,\udl{a_n}}:=L_{a_1,\cdots,a_n}|_{D_n}. \label{L_with_line}\eea
  To compute its cohomology, we use item 4 from of lem.\ref{lem_trivial_lb}
  \bea H^*(Y_n,L_{a_1,\cdots,a_n}|_{D_n})&=&H^*(D_n,\iota_n^*L_{a_1,\cdots,a_n})
  =H^*(D_n,\Big(\frac{s}{t^{n-1}}\Big)^{a_n}\iota_n^*L_{a_1-2a_n,\cdots,a_{n-1}-2a_n,0})\nn\\
  &=&H^*(Y_{n-1},\Big(\frac{s}{t^{n-1}}\Big)^{a_n}L_{a_1-2a_n,\cdots,a_{n-1}-2a_n}),\nn\eea
  reducing the calculation to $Y_{n-1}$. Now if $a_{n-1}\geq2a_n-1$, then the next pushdown $R\pi_{n-1}$ will not give higher cohomology. 
  
  Now using \eqref{quotient} repeatedly, and a double induction on $n$ and $a_n$, we see that when $a_i\geq2a_{i+1}-1$, there is no higher cohomology. 
\end{proof}
\begin{corollary}\label{cor_vanishing}
  The equivariant character of \eqref{iso_vec_space}, denoted $\chi(L_{a_1,\cdots, a_n})$ satisfies the recursion formula
\bea\chi(L_{a_1,\cdots,a_n})(s,t)=\frac{1}{1-s}\Big(\chi(L_{a_2,\cdots,a_n})(\frac{s}{t},t)-s^{a_1+1}\chi(L_{a_2,\cdots,a_n})(\frac{1}{st},t)\Big).  \label{recursion_formula}\eea
We have the vanishing result $\chi(L_{\cdots,2b-c-1,\udl{b},c,\cdots})=0$.
\end{corollary}
Note we will never actually use the equivariant character, so the proof can be safely skipped over.
\begin{proof}
  Maintaining the same assumption, we derive the equivariant character (same as $H^0$). From the above discussion
  \bea \chi(L_{a_1,\cdots,a_n})=\sum_{b_n=0}^{a_n}\chi(L_{a_1,\cdots,a_{n-1},\udl{b_n}})=\sum_{b_n=0}^{a_n}(\frac{s}{t^{n-1}})^{b_n}\chi(L_{a_1-2b_n,\cdots,a_{n-1}-2b_n}).\nn\eea
  Now apply \eqref{quotient} to $a_{n-1}-2b_n$, we write
  \bea \chi(L_{a_1,\cdots,a_n})&=&\sum_{b_n=0}^{a_n}(\frac{s}{t^{n-1}})^{b_n}\sum_{b_{n-1}=2b_n}^{a_{n-1}}\chi(L_{a_1-2b_n,\cdots,a_{n-2}-2b_n,\udl{b_{n-1}-2b_n}})\nn\\
  &=&\sum_{b_n=0}^{a_n}\sum_{b_{n-1}=2b_n}^{a_{n-1}}\frac{s^{-b_n+b_{n-1}}}{t^{(n-2)b_{n-1}-(n-3)b_n}}\chi(L_{a_1-2(b_{n-1}-b_n),\cdots,a_{n-2}-2(b_{n-1}-b_n)}).\nn\eea
  Now it is clear how to proceed, we set $c_k=b_{k+1}-b_{k+2}+\cdots\pm b_n$,
  \bea \chi(L_{a_1,\cdots,a_n})
  &=&\sum_{b_n=0}^{a_n}\sum_{b_{n-1}=2c_{n-1}}^{a_{n-1}}\frac{s^{c_{n-2}}}{t^{b_{n-1}+(n-3)c_{n-2}}}\chi(L_{a_1-2c_{n-2},\cdots,a_{n-2}-2c_{n-2}})\nn\\
  &=&\sum_{b_n=0}^{a_n}\sum_{b_{n-1}=2c_{n-1}}^{a_{n-1}}\sum_{b_{n-2}=2c_{n-2}}^{a_{n-2}}\frac{s^{c_{n-3}}}{t^{b_{n-2}+b_n+(n-4)c_{n-3}}}\chi(L_{a_1-2c_{n-3},\cdots,a_{n-3}-2c_{n-3}})\nn\\
  &=&\cdots\nn\\
  &=&\sum_{b_n=0}^{a_n} \sum_{b_{n-1}=2c_{n-1}}^{a_{n-1}}\cdots \sum_{b_2=2c_2}^{a_2}\frac{s^{c_1}}{t^{b_2+b_4+\cdots}}\chi(L_{a_1-2c_1}).\nn\eea
  Since the line $L_{\udl{1}}$ has weight $s$, we get finally
\bea \chi(L_{a_1,\cdots,a_n})
    =\sum_{b_n=0}^{a_n}\sum_{b_{n-1}=2c_{n-1}}^{a_{n-1}}\cdots\sum^{a_1}_{b_1=2c_1}s^{b_1-b_2+b_3-\cdots}t^{-b_2-b_4-\cdots}.\nn\eea
 
  To perform the summation, note the summand can be written as
  \bea s^{b_1-b_2+b_3-\cdots}t^{-b_2-b_4-\cdots}=s^{b_1-c_1}t^{-c_1-b_3-b_5-\cdots}\nn\eea
  and so the sum over $b_1$ gives
  \bea\chi(L_{a_1,\cdots,a_n})(s,t)&=&\sum_{b_n=0}^{a_n}\sum_{b_{n-1} =2c_{n-1}}^{a_{n-1}}\cdots\sum_{b_2=2c_2}^{a_2} \frac{s^{c_1}-s^{a_1-c_1+1}}{1-s}t^{-c_1-b_3-b_5-\cdots}\nn\\
  &=&\frac{1}{1-s}\sum_{b_n=0}^{a_n}\sum_{b_{n-1} =2c_{n-1}}^{a_{n-1}}\cdots\sum_{b_2=2c_2}^{a_2} \big((\frac{s}{t})^{c_1}-s^{a_1+1}(st)^{-c_1}\big)t^{-b_3-b_5-\cdots}\nn\\
  &=&\frac{1}{1-s}\Big(\chi(L_{a_2,\cdots,a_n})(\frac{s}{t},t)-s^{a_1+1}\chi(L_{a_2,\cdots,a_n})(\frac{1}{st},t)\Big).  \nn\eea
  This is the recursion relation, which curiously resembles a $q$-difference equation.
  
  The proof of the vanishing result is a direct application of the recursion formula, and is left as an exercise. The `reason' for this vanishing is given in prop.\ref{prop_pushdown_q}, namely the fibre degree is $-1$.
\end{proof}

Recall that we have a $\BB{P}^1$-fibration $\BB{P}^1\to D_{k+1}\stackrel{q}{\to} Y_{n-2}$, next we want to derive some formulae concerning the pushdown $Rq_*$ and pull back $q^*$.
\begin{proposition}\label{prop_pushdown_q}
  In the $\BB{P}^1$-bundle $q_{k+1}:\,D_{k+1}\to Y_{n-2}$, the fibrewise degree of $\iota_{k+1}^*L_{a_1,\cdots,a_n}$ (with $\iota_{k+1}$ being the inclusion $D_{k+1}\hookrightarrow Y_n$) is $a_k-2a_{k+1}+a_{k+2}$ (set $a_{k+2}=0$ if $k+1=n$).
  
  We have the formula for the pullback
  \bea &&q^*L_{b_1,\cdots,b_{n-2}}=t^{b_k}\iota^*L_{b_1,\cdots,b_{k-1},b_k,b_k,b_k,b_{k+1},\cdots,b_{n-2}},~~k+1<n-1,\nn\\
  &&q^*L_{b_1,\cdots,b_{n-2}}=t^{b_k}\iota^*L_{b_1,\cdots,b_{k-1},b_k,b_k,b_k},~~k+1=n-1,\nn\\
  &&q^*L_{b_1,\cdots,b_{n-2}}=\iota^*L_{b_1,\cdots,b_{n-2},0,0},~~k+1=n.\nn\eea
\end{proposition}
\begin{proof}
To begin, we recall the points in the proof of lem.\ref{lem_fibration} that are essential for understanding the projection.

Recall that $\ell_{k+1}$ is chosen from the 2-dim space $z^{-1}({\cal L}_k)/{\cal L}_k$. At $D_{k+1}$, we have the key feature ${\cal L}_{k+1}=z^{-1}{\cal L}_{k-1}$, or
in plain words, $\ell_k,\ell_{k+1}$ spans the 2-dim space $z^{-1}({\cal L}_{k-1})/{\cal L}_{k-1}$. It follows that the next choice of $\ell_{k+2}$ only depends on ${\cal L}_{k-1}$ but not $\ell_k,\ell_{k+1}$, the choice of $\ell_{k+3}$ only depends on ${\cal L}_{k-1},\ell_{k+2}$ etc. This exhibits the bundle structure $\BB{P}^1\to D_{k+1}\to Y_{n-2}$. Moreover, the projection map $q:\,\ell_{k+2}\to z\ell_{k+2}$ is injective, since the kernel of $z$ is in ${\cal L}_{k+1}$ to which $\ell_{k+2}$ is transverse.

  To compute the pushdown $q_*L_{a_1,\cdots,a_k,\udl{a_{k+1}},a_{k+2},\cdots}$ it is important to understand its degree along the $\BB{P}^1$-fibre. If $\deg\geq0$ (resp. $\deg<-1$), the pushdown is a vector bundle (resp. vector bundle shifted by $[-1]$) over $Y_{n-2}$, while $\deg=-1$ gives zero pushdown (which is the explanation of cor.\ref{cor_vanishing} as we promised). We can first use item 4 of lem.\ref{lem_trivial_lb} to convert $a_{k+1}D_{k+1}$ to $-2a_{k+1}(D_k+\cdots+D_1)$, and so we get $(a_k-2a_{k+1})D_k$ on the $\BB{P}^1$-fibre where $\ell_k$ belongs. This explains part of the claimed degree formula $a_k-2a_{k+1}+a_{k+2}$. To understand the $a_{k+2}$ term, note the divisor $D_{k+2}$ describes the locus where $z(\ell_{k+2})\| \ell_k\opn{mod}{\cal L}_{k-1}$. Indeed at $D_{k+2}$, $z\ell_{k+2}$ is in ${\cal L}_k$ and never in ${\cal L}_{k-1}$ and so must be $\ell_k$ up to ${\cal L}_{k-1}$.
  Now we fix $z(\ell_{k+2})$ on the base of the projection and let $\ell_k$ vary along the fibre, there is exactly one point where $z(\ell_{k+2})\|\ell_k$ mod ${\cal L}_{k-1}$, which explains the $+a_{k+2}$ term in the degree. 
  
  Knowing the fibrewise-degree, especially when the degree is $-1$, we can resort to an induction for computing $Rq_*L_{\vec a}$: 
  apply $Rq_*$ to short exact sequence $0\to L_{\vec a}(-D_k)\to L_{\vec a}\to L_{\vec a}|_{D_k}\to0$. The right most term can be computed using items 2,3,4 of lem.\ref{lem_trivial_lb}.   
  If $L_{\vec a}(-D_k)$ and $L_{\vec a}$ have nonnegative degrees along the fibre, we apply the same short exact sequence. Keep doing this until $L_{\vec a}(-D_k)$ has fibre degree $-1$ and so vanishes under $Rq_*$. We provide a concrete example of this procedure in ex\ref{ex_for_pullback}.
  
  With this knowledge we can compute the pullback. Since $q^*L_{b_1,\cdots,b_{n-2}}$ is a line bundle of fibrewise degree 0, whose pushdown gives back $L_{b_1,\cdots,b_{n-2}}$, and we get only one possible candidate: \bea t^{b_k}\iota^*L_{b_1,\cdots,b_{k-1},b_k,b_k,b_k,b_{k+1},\cdots,b_{n-2}},\nn\eea
  indeed from the short exact sequence above $t^{b_k}q_*L_{b_1,\cdots,b_{k-1},b_k,\udl{b_k},b_k,\cdots,b_{n-2}}=
  t^{b_k}q_*L_{b_1,\cdots,b_{k-1},\udl{b_k},\udl{b_k},b_k,\cdots,b_{n-2}}=q_*L_{b_1,\cdots,b_{k-1},\udl{0},\udl{0},b_k,\cdots,b_{n-2}}
  =L_{b_1,\cdots,b_{k-1},b_k,\cdots,b_{n-2}}$.
  
\end{proof}
\begin{example}\label{ex_for_pullback}
  We illustrate with the pushdown of $\iota_{k+1}^*{\cal O}(D_{k+1})$. Using the short exact sequence described in the proof of the last proposition
  \bea 0\to q_{k+1,*}\iota_{k+1}^*{\cal O}(D_{k+1})\to q_{k+1,*}\iota_{k+1}^*{\cal O}(D_k+D_{k+1})\to q_{k+1,*}\iota_{k,*}\iota_k^*\iota_{k+1}^*{\cal O}(D_k+D_{k+1})\to0.\nn\eea
  The rightmost term, by item 4 of lem.\ref{lem_trivial_lb}, is $t^{-1}{\cal O}_{Y_{n-2}}$, while the middle term vanishes due to fibre degree $-1$, hence the leftmost term is $t^{-1}{\cal O}_{Y_{n-2}}[-1]$.
\end{example}

Finally we record the Serre functors for $Y_n$.
\begin{proposition}
  The canonical bundle of $Y_n$ is
  \bea K_{Y_n}=s^nt^{-n(n-1)/2}\,L_{-2n,\cdots,-4,-2}.\nn\eea
\end{proposition}
Recall also that the canonical bundle of $D_k$ is given by the adjunction formula
\bea K_{D_k}=\iota^*K_{Y_n}(D_k),~~~\iota:\; D_k\hookrightarrow Y_n.\nn\eea
\begin{proof}
  For $Y_1$, the canonical bundle is of degree $-2$, i.e. $K\cong L_{-2}$, one needs only find out the weight. The line bundle $L_2\otimes K$ has a nowhere vanishing section $d(x_1/x_0)$ (see sec.\ref{sec_GoY}) of weight $s$ and so $K=sL_{-2}$.
  
  Assume that $K$ is obtained for $Y_n$, then for $Y_{n+1}$, consider $D_{n+1}$, the adjunction formula says 
  \[\iota_{n+1,*}K_{D_{n+1}}=K_{Y_{n+1}}(D_{n+1})|_{D_{n+1}},\]
  but $D_{n+1}=Y_n$ so $\iota_{n+1,*}K_{Y_n}=K_{Y_{n+1}}(D_{n+1})|_{D_{n+1}}$, thus 
  \bea K_{Y_{n+1}}|_{D_{n+1}}=s^nt^{-n(n-1)/2}\,L_{-2n,\cdots,-4,-2,\udl{-1}}=s^nt^{-n(n-1)/2}\frac{s}{t^n}\,L_{-2n-2,\cdots,-6,-4,\udl{-2}}.\nn\eea
  We get $K_{Y_{n+1}}=s^{n+1}t^{-(n+1)n/2}L_{-2n-2,\cdots,-4,-2}$.
\end{proof}
\begin{corollary}(lem 4.4 \cite{Cautis2007KnotHV})
  From the Serre functors we get the two sided adjoints to $q_k^*$ and $\iota_{k*}$
  \bea 
  &&tRq_*{\cal O}_D(D)[1] \dashv q^* \dashv Rq_* \dashv t{\cal O}_D(D)q^*[1],\nn\\
  &&\iota_{k*}{\cal O}_{D_k}(-D_k)[1]\dashv L\iota^*_k \dashv \iota_{k*} \dashv L\iota^* {\cal O}(D_k)[-1]\nn\eea
\end{corollary}
In the equation above, a line bundle such as ${\cal O}_D$ represents the operation given by tensoring with it.
\begin{proof}
  Starting from $q^*$, it has a natural right adjoint. To get its left adjoint, one can follow the standard manipulation to use Serre dual to place it at the `correct slot'
  \bea \hom_D(F,q^*G)&=&\hom_D(q^*G,F\otimes K_D[n-1])^{\vee}=\hom_{Y_{n-2}}(G,Rq_*(F\otimes K_D)[n-1])^{\vee}\nn\\
  &=&\hom_{Y_{n-2}}(Rq_*(F\otimes K_D)[1],G\otimes K_{Y_{n-2}}),\nn\eea
  hence the left adjoint acting on $F$ is $K_{Y_{n-2}}^{-1}\otimes Rq_*(F\otimes K_D)[1]$.
  We can simplify using the projection formula
  \bea q^*K_{Y_{n-2}}^{-1}&=&s^{2-n}t^{(n-2)(n-3)/2}q^*L_{2(n-2),\cdots,{2(n-k)},\cdots,2}\nn\\
  &\stackrel{\ref{prop_pushdown_q}}{=}&s^{2-n}t^{(n-2)(n-3)/2}t^{2(n-k)}L_{2(n-2),\cdots,2(n-k),\udl{2(n-k)},2(n-k),\cdots,2}\nn\\
  &\stackrel{\ref{lem_trivial_lb}}{=}&s^{2-n}t^{(n-2)(n-3)/2}t^{2(n-k)}(st^{-k+1})^{-2}L_{2n,\cdots,2(n+2-k),\udl{2(n+1-k)},2(n-k),\cdots,2}\nn\\
  &=&tL\iota^*K_{Y_n}^{-1}=tK_D^{-1}(D).\nn\eea
  From this we get the stated formula for the left adjoint for $q^*$. The calculation for the other adjoints are similar. 
\end{proof}
From these adjoints we can get the adjoint functors for the cup and cap. 
\begin{corollary}(Lem.4.4 of \cite{Cautis2007KnotHV})\label{cor_2_side_adj}
We have the two-sided adjoints
\bea t^{1/2}{}_k\CAP[1] \dashv {}^k\CUP \dashv t^{-1/2}{}_k\CAP[-1].\nn\eea  
\end{corollary}

\subsection{Tilting Bundle}\label{sec_TB}
Let $X$ be a smooth variety, a tilting bundle (or object) $T$ is a generator of $D^b(X)$ that satisfies $\opn{Ext}^i(T,T)=0$ for $i>0$. We let $\Gl=\opn{Ext}^0(T,T)$ be the endomorphism algebra of $T$, let further $D^b(\opn{mod}_{\Gl})$ be the bounded derived category of right-$\Gl$-modules, then we have a pair of derived equivalences
\bea& \Psi:~~D^b(\opn{mod}_{\Gl})\rightleftharpoons D^b(X):~~\Phi\label{equivalence}\\
&\Psi(-)=-\otimes_{\Gl}^{\BB{L}}T,~~~\Phi(-)=\BB{R}\!\hom(T,-).\nn\eea
The composition $\Psi\Phi$ is quasi-isomorphic to the identity functor: $F\simeq \Psi\Phi(F)$ for any $F\in D^b(X)$. This key point here is often phrased as saying that we have a resolution of the diagonal ${\cal O}_{\Gd}$ in terms of $T$
\bea T^{\vee}\!\otimes^{\BB{L}}_{\Gl}T\simeq {\cal O}_{\Gd}~\in D^b(X\times X).\label{res_diag}\eea
Then applying $\Psi\Phi$ is just applying the Fourier--Mukai transform with ${\cal O}_{\Gd}$ as the kernel, but this is the identity transform, except now $\Psi\Phi(F)$ is resolved in terms of objects that constitute $T$.
\begin{example}\label{ex_Gl_P^1}
  The simplest example is $\BB{P}^1=\BB{P}(V)$ for $V$ a vector space of dimension 2.
  Let $Q$ be the tautological quotient bundle, fitting in the exact Euler sequence
  \bea 0\to {\cal O}(-1)\to V\otimes{\cal O}\to Q\to0.\nn\eea
  We have the resolution of ${\cal O}_{\Gd}$ as ${\cal O}\boxtimes{\cal O}$ module
  \bea 0\to {\cal O}(-1)\boxtimes Q^*\to{\cal O}\boxtimes{\cal O}\to {\cal O}_{\Gd}\to0.\label{Orlov_res}\eea
  To describe the middle map, pick basis $\{e_i\}$ for $V$ and dual basis $\{e^i\}$ for $V^*$. The middle map goes as: $\BB{C}\to e^i\boxtimes e_i$, then evaluates $e^i$ on ${\cal O}(-1)$ and $Q^*$ on $e_i$. The image vanishes precisely on the diagonal, hence the exactness of the sequence.
  
  For us, it is more convenient to write another resolution of ${\cal O}_{\Gd}$, which is of course quasi-isomorphic to \eqref{Orlov_res} (but not chain homotopic)
  \bea \begin{tikzpicture}
  \matrix (m) [matrix of math nodes, row sep=1.2em, column sep=4em]
    {   {\cal O}(-1)\boxtimes V^*\boxtimes {\cal O} & {\cal O}\boxtimes {\cal O} \\
     & {\cal O}(-1)\boxtimes {\cal O}(1) \\   };
  \path[->,font =\scriptsize]
  (m-1-1) edge node[above] {$e^i\boxtimes1$} (m-1-2)
  (m-1-1) edge node[below] {$-1\boxtimes e^i$} (m-2-2);
  \node at(2.1,0.05) {$\oplus$};
  
  \node at(4.25,-0.05) {$\simeq {\cal O}_{\Gd}$};
\end{tikzpicture}.\label{Beilinson_res}\eea
To see that they are quasi-isomorphic, take the dual of the Euler-sequence above
\bea 0\to Q^* \to V^*\otimes{\cal O}\to {\cal O}(1)\to0.\nn\eea
So we can replace the southeast arrow in \eqref{Beilinson_res} with the kernel ${\cal O}(-1)\boxtimes Q^*$.
But as the Euler sequence does not split, this replacement is not effected via a chain homotopy.
  
The new resolution is more in line with the tilting picture given above: let $T={\cal O}\oplus{\cal O}(1)$ in \eqref{equivalence}.
The quiver algebra $\Gl=\opn{End}(T,T)$ reads
\bea \begin{tikzpicture}
  \matrix (m) [matrix of math nodes, row sep=1.2em, column sep=5em]
    {   e_1 & e_0 \\ };
  \path[->,font =\scriptsize, bend right = 20]
  (m-1-2) edge node[above] {$z^0$} (m-1-1);
  \path[->,font =\scriptsize, bend left = 20]
  (m-1-2) edge node[below] {$z^1$} (m-1-1);
\end{tikzpicture}\nn\eea
where the nodes $e_0/e_1$ correspond to the objects ${\cal O}(0/1)$, while the arrows are the \emph{irreducible} maps between the objects. In this case $\hom({\cal O},{\cal O}(1))=V^*$ is 2-dimensional, generated by the two homogeneous coordinates $z^0,z^1$, leading to the two arrows. Now the derived tensor product $T^{\vee}\otimes_{\Gl}^{\BB{L}}T$ gives the complex \eqref{Beilinson_res}.

Finally we write \eqref{Beilinson_res} in a \emph{matrix form} for later use
\bea 
\begin{tikzpicture}
  \node at (0,0) {${\cal O}_{\Gd}=\Big[{\cal O},~{\cal O}(-1)\Big]\Gd_1 
\begin{bmatrix}
  {\cal O} \\ {\cal O}(1)
\end{bmatrix}~=$}; 
\node at (0,-.6) {\phantom{x}};
\end{tikzpicture}
\begin{tikzpicture}
    \node at (0,0) {$\Big[{\cal O},~~{\cal O}(-1)\Big]$}; 
    \node at (0,-.6) {\phantom{x}};
\end{tikzpicture}
\begin{tikzpicture}
\matrix (m) [left delimiter={[},right delimiter={]},matrix of math nodes, row sep=1.4em, column sep=.7em]
    { \boxtimes  & 0 \\
      \boxtimes V^*[1]\boxtimes  & \boxtimes \\ };
      \path[->,font =\scriptsize]
      (m-2-1) edge[red] (m-1-1)
      (m-2-1) edge (m-2-2);
\end{tikzpicture}
\begin{tikzpicture}
    \node at (0,0) {$\begin{bmatrix}
  {\cal O} \\ \\ {\cal O}(1)
\end{bmatrix}$};
\node at (0,-.6) {\phantom{x}}; 
\end{tikzpicture}\label{style}\eea
where the red arrow is induced from ${\cal O}(-1)\otimes V^*\to {\cal O}$ and the black one $V^* \otimes {\cal O}\to {\cal O}(1)$.
In the future we omit the tensor product or simply write $1$ for $\boxtimes$. 

We illustrate working with these matrices by computing the kernels of the Serre dual and its inverse
\bea \Gs=\Gd_*\go[1],~~~\Gs^{-1}=\Gd_*\go^{-1}[-1],~~~{\rm where}~\go=s{\cal O}(-2).\label{Serre_kernel}\eea
This is done by the following `matrix multiplication'
\bea \Gs=\begin{bmatrix}
  {\cal O} & {\cal O}(-1)
\end{bmatrix}
\begin{bmatrix}
  1 & \\ V^*[1] & 1
\end{bmatrix}
\begin{bmatrix}
  \pi_*(\go)[1] & \pi_*(\go(-1))[1] \\
  \pi_*(\go(1))[1] & \pi_*(\go)[1] \\
\end{bmatrix}\begin{bmatrix}
  1 & \\ V^*[1] & 1
\end{bmatrix}
\begin{bmatrix}
  {\cal O} \\ {\cal O}(1)
\end{bmatrix}\nn\eea
where the $i$-$j$ entry of the middle matrix is obtained as $\hom({\cal O}(j-1),\Gs {\cal O}(i-1))=\pi_*\go(i-j)[1]$.

We work on the middle three terms
\bea &&
\begin{bmatrix}
  1 & \\ V^*[1] & 1
\end{bmatrix}
\begin{bmatrix}
  1 & V \\
  0 & 1 \\
\end{bmatrix}\begin{bmatrix}
  1 & \\ V^*[1] & 1
\end{bmatrix}=\begin{bmatrix}
  1 & \\ V^*[1] & 1
\end{bmatrix}
\begin{bmatrix}
  1\ot V\otimes V^*[1] & V \\ V^*[1] & 1
\end{bmatrix}\nn\\
&=&\begin{bmatrix}
  1 & \\ V^*[1] & 1
\end{bmatrix}
\begin{bmatrix}
  V_{ad}[1] & V \\ V^*[1] & 1
\end{bmatrix}=\begin{bmatrix}
  V_{ad}[1] & V \\ V^*\otimes V_{ad}[2]\to V^*[1] & V^*\otimes V[1]\to 1
\end{bmatrix}\nn\\
&=&\begin{bmatrix}
  V_{ad}[1] & V \\ s^{-1}S^3V^*[2] & V_{ad}[1]
\end{bmatrix}\nn\eea
where $V,V^*,V_{ad}$ are treated as the (anti-)fundamental and adjoint representations of $SU(2)$.
We also used chain homotopies to cancel isomorphic terms in the complex.\footnote{Called Gaussian eliminations in \cite{bar2007fast} lem.4.2. One has to watch out for the zig-zag corrections.} From the final matrix form, one gets the kernel for $\Gs$
\bea \begin{tikzpicture}
  \matrix (m) [matrix of math nodes, row sep=.4em, column sep=0.5em]
    { & {\cal O}\otimes V_{ad}\otimes {\cal O} & \\
    {\cal O}(-1)\otimes s^{-1}S^3V^*\otimes {\cal O} & & {\cal O}\otimes V \otimes {\cal O}(1) \\
    & {\cal O}(-1)\otimes V_{ad}\otimes {\cal O}(1) & \\};
    \path[->, font=\scriptsize]
  (m-2-1) edge (m-1-2)
  (m-2-1) edge (m-3-2)
  (m-1-2) edge (m-2-3)
  (m-3-2) edge (m-2-3);
  \node at (-5.5,0) {\small$\Gs\sim$};
  \node at (3.5,-.5) {\small$\stackrel{\uparrow}{\deg0}$};
\end{tikzpicture}\label{kernel_Gs}\eea
where all maps can be phrased in terms of $SU(2)$ representations, e.g. the two right arrows are: $V_{ad}\subset V\otimes V^*$, then multiply the $V^*$ factor to the left or right. The two left arrows are: $s^{-1}S^3V^*\subset V^*\otimes V_{ad}$, then multiply $V^*$ to the left or right.
Similarly the inverse Serre functor has kernel
\bea \begin{tikzpicture}
  \matrix (m) [matrix of math nodes, row sep=.3em, column sep=0.5em]
    { {\cal O}\otimes{\cal O} & \\
     & {\cal O}\otimes V\otimes {\cal O}(1) \\
     {\cal O}(-1)\otimes{\cal O}(1) & \\};
    \path[->, font=\scriptsize]
  (m-1-1) edge (m-2-2)
  (m-3-1) edge (m-2-2);
  \node at (-3.5,0) {\small$\Gs^{-1}\sim$};
  \node at (1,-.5) {\small$\stackrel{\uparrow}{\deg1}$};
  \node at (-1.3,0) {\small$\oplus$};
\end{tikzpicture}\label{kernel_Gs_inv}\eea
with the maps go as $\BB{C}\subset V\otimes V^*$, then again multiply $V^*$ to left or right.

\end{example}
In \cite{Orlov_1993} Orlov gave a \emph{fibrewise} version of the above resolution of the diagonal.
Let $E\stackrel{\pi}{\to}X$ be a rank 2 vector bundle, which is the only case we need. One can write down a resolution of the diagonal $\Gd\subset \BB{P}(E)\times \BB{P}(E)$
  \bea 0\to {\cal O}(-1)\boxtimes Q^*\to{\cal O}\boxtimes{\cal O}\to {\cal O}_{\Gd}\to0\nn\eea
  that looks identical to the $\BB{P}^1$ case earlier except now $\boxtimes$ stands for: 
  \bea A\boxtimes B:=p_1^*A\otimes_{{\cal O}_X} p_2^*B,~~~p_{1,2}:\;\BB{P}(E)\times_X\BB{P}(E)\to\BB{P}(E).\label{boxtimes}\eea
  We also have the fibrewise version of \eqref{Beilinson_res}, essentially just swapping out $\boxtimes$ for its fibrewise version\eqref{boxtimes}, and $V^*$ for $\pi_*{\cal O}(1)$.
  
  That we can build bi-module resolution via the $\BB{P}^1$-fibration is crucial for us, since our space $Y_n$ is a repeated $\BB{P}^1$-fibration, and now we can apply the tilting step by step
  \bea Y_n\stackrel{\pi_n}{\to} Y_{n-1}\stackrel{\pi_{n-1}}{\to} \cdots\to *\nn\eea
  using \eqref{Beilinson_res} for each step. We denote the resolution in the matrix form as in \eqref{style}
  \bea {\cal O}_{\Gd_{Y_n}}=[L_{\sbullet}^{\vee}]\Gd_n [L_{\sbullet}],\label{resolution_Delta_n}\eea
  with $L_{\sbullet}$ being the exceptional collection of thm.\ref{thm_exp_collection_Yn}. 
  \begin{lemma}
    Writing the resolution of the structure sheaf of the diagonal $\Gd_{Y_n}\subset Y_n\times Y_n$ as \eqref{resolution_Delta_n}, the matrix $\Gd_n$ can be obtained iteratively from $\Gd_{n-1}$
\bea \Gd_n=
    \begin{bmatrix}
      \Gd_{n-1} & 0 \\ 
      \Gd_{n-1}H\Gd_{n-1}[1] & \Gd_{n-1}
    \end{bmatrix},~~~~
    H_{IJ}=\pi_*(L_I\otimes_{{\cal O}_{Y_{n-1}}}\pi_{n*}(L_{2^{n-1}\cdots21})\otimes_{{\cal O}_{Y_{n-1}}}L_J^{\vee}),\nn\eea
where $L_I,L_J$ run over the exceptional collection for $Y_{n-1}$. 
  \end{lemma}
The maps in the matrix are also described recursively, given in the proof.
\begin{proof}
    Assume we have obtained $\Gd_{n-1}$, that is, a derived tensor product $A\otimes^{\BB{L}}_{{\cal O}_{Y_{n-1}}}B$ can be written using \eqref{resolution_Delta_n}
    \bea A\otimes^{\BB{L}}_{{\cal O}_{Y_{n-1}}}B=[A\otimes_{{\cal O}_{Y_{n-1}}}L_{\sbullet}^{\vee}]\Gd_{n-1} [L_{\sbullet}\otimes_{{\cal O}_{Y_{n-1}}}B].\label{used_I}\eea
    where $L_{\sbullet}$ runs over the exceptional collection for $Y_{n-1}$.
    
    For $Y_n$, we consider the fibration $\BB{P}^1\to Y_n\stackrel{\pi_n}{\to} Y_{n-1}$. 
    To write the fibrewise version of \eqref{Beilinson_res}, we need two line bundles of fibre degree 0 and 1, which we choose to be 
    $L_{0,\cdots,0}$, and $L_{2^{n-1},\cdots,1}$. Then we have the quasi-isomorphism
    \bea \begin{tikzpicture}
  \matrix (m) [matrix of math nodes, row sep=1em, column sep=0.5em]
    { L_{2^{n-1}\cdots1}^{\vee}  \otimes_{{\cal O}_{Y_{n-1}}} \pi_{n*}L_{2^{n-1}\cdots1}\otimes_{{\cal O}_{Y_{n-1}}}L_{0\cdots0}  
    & L_{0\cdots0}^{\vee}\otimes_{{\cal O}_{Y_{n-1}}} L_{0\cdots0} \\
    & L_{2^{n-1}\cdots1}^{\vee}\otimes_{{\cal O}_{Y_{n-1}}}L_{2^{n-1}\cdots1} \\};
    \path[->, font=\scriptsize]
  (m-1-1) edge[red] (m-1-2)
  (m-1-1) edge (m-2-2);
  \node at (-7,0) {${\cal O}_{\Gd_{Y_n}}\sim$};
  \node at (3,0) {$\oplus$};
\end{tikzpicture}\nn\eea
where the red/black map is multiplying $\pi_{n*}L_{2^{n-1}\cdots1}$ to the left or right.
In matrix form, this is identical to \eqref{style}, but with $V^*$ replaced with a rank 2 vector bundle $\pi_{n*}L_{2^{n-1}\cdots1}$ over $Y_{n-1}$
\bea 
\begin{tikzpicture}
    \node at (-3,0) {${\cal O}_{\Gd_{Y_n}}\sim$};
    \node at (0,0) {$\Big[L_{0\cdots0}^{\vee},~~L_{2^{n-1}\cdots1}^{\vee}\Big]$}; 
    \node at (0,-.7) {\phantom{x}};
\end{tikzpicture}
\begin{tikzpicture}
\matrix (m) [left delimiter={[},right delimiter={]},matrix of math nodes, row sep=1.4em, column sep=.7em]
    { {\cal O}_{Y_{n-1}}  & 0 \\
      \pi_{n*}L_{2^{n-1}\cdots1}[1]  & {\cal O}_{Y_{n-1}} \\ };
      \path[->,font =\scriptsize]
      (m-2-1) edge[red] (m-1-1)
      (m-2-1) edge (m-2-2);
      \node at (2.5,-0.1) {$\otimes$};
      \node at (-2.5,-0.1) {$\otimes$};
\end{tikzpicture}
\!\!
\begin{tikzpicture}
    \node at (0,0) {$\begin{bmatrix}
  L_{0\cdots0} \\ \\ L_{2^{n-1}\cdots1}
\end{bmatrix}$};
\node at (0,-.6) {\phantom{x}}; 
\end{tikzpicture},\nn\eea
where $\otimes$ is over ${\cal O}_{Y_{n-1}}$. The red/black map is again multiplying $\pi_{n*}L_{2^{n-1}\cdots1}$ to the left/right.

Now we iterate the same steps by writing the tensor product $\otimes_{{\cal O}_{Y_{n-1}}}$ using the resolution of ${\cal O}_{\Gd_{Y_{n-1}}}$ as in \eqref{used_I}. 
For example 
\bea L^{\vee}_{0\cdots0}\otimes_{{\cal O}_{Y_{n-1}}}L_{0\cdots0}=[L^{\vee}_{I0}]\Gd_{n-1}[L_{J0}].\nn\eea
This gives the diagonal blocks. For the lower left entry, resolve $\pi_{n*}L_{2^{n-1}\cdots1}$ as 
\bea L_I^{\vee}\otimes(\Gd_{n-1})_{IJ}\underbrace{\pi_*(L_J\otimes_{{\cal O}_{Y_{n-1}}}\pi_{n*}L_{2^{n-1}\cdots1}\otimes_{{\cal O}_{Y_{n-1}}}L^{\vee}_M)}_{H}(\Gd_{n-1})_{MN}\otimes L_N.\nn\eea
From this we get the lower left $\Gd_{n-1}H\Gd_{n-1}[1]$ block. The maps upwards is the multiplication of $\Gd H$ to the left, and the one rightwards multiplies $H\Gd$ to the right.
\end{proof}
We summarise the above iteration procedure as: to obtain $\Gd_n$, start with the $2\times 2$ matrix
\bea 
\begin{bmatrix}
  {\cal O}_{Y_{n-1}}  & 0 \\
    \pi_{n*}L_{2^{n-1}\cdots1}[1] & {\cal O}_{Y_{n-1}}
\end{bmatrix}\nn\eea
where all entries are vector bundles over $Y_{n-1}$. For each entry, say $E$, replace it with the $2\times2$ matrix
\bea \begin{bmatrix}
  {\cal O}_{Y_{n-2}}  & 0 \\
    \pi_{n-1,*}L[1] & {\cal O}_{Y_{n-2}}
\end{bmatrix}
\begin{bmatrix}
  \pi_{n-1,*}E  & \pi_{n-1,*}(L^{\vee}\otimes_{{\cal O}_{Y_{n-1}}}E) \\
    \pi_{n-1,*}(L\otimes_{{\cal O}_{Y_{n-1}}}E) & \pi_{n-1,*}E
\end{bmatrix}
\begin{bmatrix}
  {\cal O}_{Y_{n-2}}  & 0 \\
    \pi_{n*}L[1] & {\cal O}_{Y_{n-2}}
\end{bmatrix},\label{eq_summary}\eea
where $L=L_{2^{n-2}\cdots1}$ and the unwritten tensor products are over ${\cal O}_{Y_{n-2}}$. Then one repeats with $L=L_{2^{n-3}\cdots1}$ etc. As we said above, keeping  track of the maps in $\Gd_n$ is laborious, but on the other hand they are usually self-evident, up to sign.

We will demonstrate the process by iterating twice, thereby expressing $\Gd_n$ in terms of $\Gd_{n-2}$. According to our summary \eqref{eq_summary}, we will take $E$ to be either ${\cal O}_{Y_{n-1}}$ or $\pi_{n*}L_{2^{n-1},\cdots,1}[1]$. %
  Let $I,J$ be the multi-indices $I=(2^{n-2},2^{n-3},\cdots,1)$, $J=(2^{n-3},2^{n-4},\cdots,1)$ with $I=(2J,1)$. 
  For $E=\pi_{n*}L_{2I,1}$, we need to compute
  \bea \Gd_{n-1}H\Gd_{n-1}=
  \begin{bmatrix}
    {\cal O}_{Y_{n-2}} & \\ \pi_{n-1,*}L_I[1] & {\cal O}_{Y_{n-2}}
  \end{bmatrix}
  \begin{bmatrix}
    \pi_{n-1,*}\pi_{n*}L_{2I,1} & \pi_{n-1,*}\pi_{n*}L_{I,1} \\
    \pi_{n-1,*}\pi_{n*}L_{3I,1} & \pi_{n-1,*}\pi_{n*}L_{2I,1} 
  \end{bmatrix}\begin{bmatrix}
    {\cal O}_{Y_{n-2}} & \\ \pi_{n-1,*}L_I[1] & {\cal O}_{Y_{n-2}}
  \end{bmatrix}.\nn\eea
  Multiplying the right two matrices will require computing two cones
  \bea \opn{cone}(\pi_{n-1,*}\pi_{n*}L_{2I,1}\otimes_{{\cal O}_{Y_{n-2}}}\pi_{n-1,*}L_I\to \pi_{n-1,*}\pi_{n*}L_{3I,1}),~~~A\nn\\
  \opn{cone}(\pi_{n-1,*}\pi_{n*}L_{I,1}\otimes_{{\cal O}_{Y_{n-2}}}\pi_{n-1,*}L_I\to \pi_{n-1,*}\pi_{n*}L_{2I,1}),~~~B.\nn\eea
  All the maps here are canonical, so the calculation is mechanical. But we take a shortcut here by noting that the multiplication of the two right matrices is what one would do to resolve $\pi_{n*}L_{2I,1}$ and $
  \pi_{n*}L_{3I,1}$ using ${\cal O}_{Y_{n-2}}$ and $L_I$. Indeed the first/second row are the `coefficients' of ${\cal O}_{Y_{n-2}}$ and $L_I$ in the resolution of $\pi_{n*}L_{2I,1}$ and $
  \pi_{n*}L_{3I,1}$.
  
  Knowing this, we write calculate $\pi_{n*}L_{K,1}$ (where $K=2I,3I$) as 
  \bea \pi_{n,*}L_{K,1}=\pi_{n,*}\opn{cone}(L_{K,\udl{1}}[-1]\to L_{K,0})\stackrel{lem.\ref{lem_trivial_lb}}{=}\opn{cone}(\frac{s}{t^{n-1}}L_{K-2}[-1]\to L_K)\nn\eea
  where by $K-2$ we mean subtracting 2 from every entry of $K$. For $K=2I$, the final ($(n-1)^{th}$) index of $L_K$ is 2, we manipulate it so it is always 0 or 1
  \bea L_{2I}=L_{4J,2}&=&\opn{cone}(L_{4J,\udl{2}}[-1]\to L_{4J,1})=\opn{cone}(\frac{s}{t^{n-2}}L_{4J-2,\udl{1}}[-1]\to L_{4J,1})\nn\\
  &=&\opn{cone}(\frac{s}{t^{n-2}}\opn{cone}(L_{4J-2,0}\to L_{4J-2,1})[-1]\to L_{4J,1}).\nn\eea
  Inserting this back into $\pi_{n,*}L_{K,1}$ gives us the first row of the matrix mentioned above
  \bea 
  \begin{bmatrix}
    \frac{s}{t^{n-2}}L_{4J-2}([1]+\frac{1}{t}) & \frac{s}{t^{n-2}}L_{2J-2}\oplus L_{2J}
 \end{bmatrix} \nn\eea
 the first entry being the `coefficient' of ${\cal O}_{Y_{n-2}}$ and second that of $L_I$.
 Doing the same for $K=3I$
  \bea L_{6J,3}&=&\opn{cone}(\frac{s^2}{t^{2n-4}}L_{6J-4,\udl{1}}[-1]\to L_{6J,2})
  =\opn{cone}(\frac{s^2}{t^{2n-4}}L_{6J-4,\udl{1}}[-1]\to \opn{cone}(\frac{s}{t^{n-2}}L_{6J-2,\udl{1}}[-1]\to L_{6J,1}))\nn\\
  &=&\opn{cone}(\frac{s^2}{t^{2n-4}}\opn{cone}(L_{6J-4,0}\to L_{6J-4,1})[-1]\to \opn{cone}(\frac{s}{t^{n-2}}\opn{cone}(L_{6J-2,0}\to L_{6J-2,1})[-1]\to L_{6J,1})).\nn\eea
  This gives the second row
  \bea 
  \begin{bmatrix}
    \frac{s^2}{t^{2n-4}}L_{6J-4}[1]\oplus\frac{s}{t^{n-2}}L_{6J-2}[1] & \frac{s^2}{t^{2n-4}}L_{4J-4}\oplus \frac{s}{t^{n-2}}(1+\frac1t)L_{4J-2}\oplus L_{4J}
  \end{bmatrix}\nn\eea
  This leads to
  \bea \Gd_{n-1}H\Gd_{n-1}=\begin{bmatrix}
    {\cal O}_{Y_{n-2}} & \\ \pi_{n-1,*}L_I[1] & {\cal O}_{Y_{n-2}}
  \end{bmatrix}\begin{bmatrix}
  \frac{s}{t^{n-2}}L_{4J-2}([1]+\frac1t) & \frac{s}{t^{n-2}}L_{2J-2}\oplus L_{2J} \\
    \frac{s^2}{t^{2n-4}}L_{6J-4}[1]\oplus\frac{s}{t^{n-2}}L_{6J-2}[1] & \frac{s^2}{t^{2n-4}}L_{4J-4}\oplus \frac{s}{t^{n-2}}(1+\frac1t)L_{4J-2}\oplus L_{4J}
  \end{bmatrix}\nn\eea
  The final matrix multiplication is again mechanical, we note $\pi_{n-1,*}L_I=\opn{cone}(s/t^{n-2}L_{2J-2}[-1]\to L_{2J})$, leading to Gaussian eliminations
  \bea \Gd_{n-1}H\Gd_{n-1}=
  \begin{bmatrix}
    \frac{s}{t^{n-2}}L_{2I-2}([1]+\frac{1}{t}) & \frac{s}{t^{n-2}}L_{2J-2}\oplus L_{2J} \\
    \frac{s^2}{t^{2n-3}}L_{6J-4}[1]+\frac{s}{t^{n-1}}L_{6J-2}[1] & \frac{s}{t^{n-2}}L_{2I-2}([1]+\frac{1}{t})
  \end{bmatrix}.\nn\eea
  We remark that, in the lower left and top right entries, there are maps between the two terms, as the object $s^{-1}t^{n-2}L_{2[n-2]}$ always has a 1-dim cohomology at degree 1 and weight 0, for all $n>3$.
\begin{example}\label{ex_Gl_Y2}  
  As a concrete example, set $n=2$ for $Y_2$, so $I=1$ and $J=\emptyset$, then the matrix $\Gd_2$ in \eqref{style} reads
  \bea \Gd_2=
\begin{bmatrix}
  1 & 0 & 0 & 0 \\ (1+s)[1] & 1 & 0 & 0 \\
  s([1]+\frac{1}{t})[1] & (1+s)[1] & 1 & 0\\
  \frac{s}{t}(1+s)[2] & s([1]+\frac{1}{t})[1] & (1+s)L_1[1] & 1
\end{bmatrix}\nn.\eea
One can keep track of the maps between the entries in the above iterative process, but by writing the above matrix slightly differently, the maps become more self-explanatory
\bea {\cal O}_{\Gd_{Y_2}}\sim\Big[
  L^{\vee}_{00} ~~ L^{\vee}_{10} ~~ L^{\vee}_{21} ~~ L^{\vee}_{31}
\Big]
\begin{bmatrix}
  1 & 0 & 0 & 0 \\ V^*[1] & 1 & 0 & 0 \\
  \wedge^2V^*[2]\oplus\frac{s}{t}[1] & V^*[1] & 1 & 0\\
  \frac{s}{t}V^*[2] & \wedge^2V^*[2]\oplus\frac{s}{t}[1] & V^*[1] & 1
\end{bmatrix}
\begin{bmatrix}
  L_{00} \\ L_{10} \\ L_{21} \\ L_{31}
\end{bmatrix}\label{res_Y2}\eea
We point out some maps, looking at the first column, we have
\bea \frac{s}{t}\underbrace{V^*\otimes L_{31}^{\vee}}_{\deg-2}\to \frac{s}{t} \underbrace{L_{21}^{\vee}}_{\deg-1},~~~~\frac{s}{t}\underbrace{\wedge^2V^*\otimes L_{21}^{\vee}}_{\deg-2}\to \frac{s}{t} \underbrace{V^*L_{10}^{\vee}}_{\deg-1},~~~~\frac{s}{t}\underbrace{\otimes L_{21}^{\vee}}_{\deg-1}\to \underbrace{L_{00}^{\vee}}_{\deg0},~~~~\underbrace{V^*\otimes L_{10}^{\vee}}_{\deg-1}\to \underbrace{L_{00}^{\vee}}_{\deg0}.\nn\eea
In short, there are arrows going up or right between two entries whenever their degrees differ by 1 and permitted by equivariance. %
\end{example}

\begin{example}\label{ex_Gl_Y3}
For $n=3$, we set $I=(2,1)$ and $J=1$, and ${\cal O}_{\Gd_{Y_3}}$ is resolved as
\bea {\cal O}_{\Gd_{Y_3}}\sim 
\begin{bmatrix}
  L_{000}^{\vee} & L_{210}^{\vee} & L_{421}^{\vee} & L_{631}^{\vee}
\end{bmatrix}
\begin{bmatrix}
  1 & \\
  \pi_{2*}L_{21}[1] & 1 & \\
  \frac{s}{t}L_{2}[2]+\frac{s}{t^2}L_{2}[1] & \pi_{2*}L_{21}[1] & 1 & \\ 
  \frac{s}{t^2}\pi_{2*}L_{41}[2] & \frac{s}{t}L_{2}[2]+\frac{s}{t^2}L_{2}[1] & \pi_{2*}L_{21}[1] & 1
\end{bmatrix}\begin{bmatrix}
  L_{000} \\ L_{210} \\ L_{421} \\ L_{631}
\end{bmatrix}\nn\eea
where all tensor products (unwritten) are over ${\cal O}_{Y_1}$. Carrying out the iterative step once more, we get
\bea &&{\cal O}_{\Gd_3}\sim
\begin{bmatrix}
  L_{\sbullet\sbullet\sbullet}^{\vee}
\end{bmatrix}
\underbrace{\begin{bmatrix}
  \Gd_2 & 0 \\
  \nabla_2[1] & \Gd_2
\end{bmatrix}}_{\Gd_3}
\begin{bmatrix}
  L_{\sbullet\sbullet\sbullet}
\end{bmatrix}\nn\\
&&\nabla_2[1]=
\begin{bmatrix}
  \frac{s^2}{t}[3]+\frac{s^2}{t^2}[2] & \frac{s}{t}V^*[2]+\frac{s}{t^2}V^*[1] & s[2]+\frac{s}{t}[1] & V^*[1] \\
  0 & \frac{s^2}{t}[3]+\frac{s^2}{t^2}[2] & \frac{s}{t}V^*[2] & s[2]+\frac{s}{t}[1] \\
  (\frac{s^2}{t^2}S^2V^*+\frac{s^3}{t^3})[3] & (\frac{s}{t^2}S^3V^*+\frac{s^2}{t^3}V^*)[2] & \frac{s^2}{t}[3]+\frac{s^2}{t^2}[2] & \frac{s}{t}V^*[2]+\frac{s}{t^2}V^*[1] \\
  \frac{s^3}{t^2}V^*[4] & (\frac{s^2}{t^2}S^2V^*+\frac{s^3}{t^3})[3] & 0 & \frac{s^2}{t}[3]+\frac{s^2}{t^2}[2] 
\end{bmatrix}\nn\eea
\end{example}
One can also apply Orlov's resolution `from the middle': since $Y_n|_{D_{k+1}}\to Y_{n-2}$ is a $\BB{P}^1$ bundle, ex.\ref{ex_Gl_Y2} applies. For this we need to choose two line bundles of degree 0 and 1 along the $\BB{P}^1$-fibre. The choices are not unique, for example we can pick $L_0={\cal O}_{D_{k+1}}$ and $L_1=\iota^*_{k+1}L_{-[k+1]}$ (see prop.\ref{prop_pushdown_q}).
We get the resolution for the diagonal $\Gd_{D_{k+1}}$
\bea \begin{tikzpicture}
  \matrix (m) [matrix of math nodes, row sep=1.2em, column sep=4em]
    {   L_1^{\vee}\boxtimes_{{\cal O}_{Y_{n-2}}}(q_*L_1)\boxtimes_{{\cal O}_{Y_{n-2}}}L_0 & L_0^{\vee}\boxtimes_{{\cal O}_{Y_{n-2}}} L_0 \\
     & L_1^{\vee}\boxtimes_{{\cal O}_{Y_{n-2}}}L_1 \\   };
  \path[->,font =\scriptsize]
  (m-1-1) edge[red] (m-1-2)
  (m-1-1) edge (m-2-2);
  \node at(2.7,0.05) {$\oplus$};
  \node at(4.75,-0.05) {$\simeq {\cal O}_{\Gd_D}$};
\end{tikzpicture}\label{Beilinson_res_halfway}\eea
where $q=q_{k+1},\iota=\iota_{k+1}, D=D_{k+1}$. We will use this in the proof of prop.\ref{prop_Kauffman}.

\subsection{Comparison with \texorpdfstring{\cite{King1997TiltingBO}}{[Kin97]} and \texorpdfstring{\cite{BUTLER1999323}}{[BK99]}}\label{sec_CwKaB}
So far we have obtained the resolution of ${\cal O}_{\Gd}$ as ${\cal O}_{Y_n\times Y_n}$ modules via an iterative method. In terms of the equivalent module category $D^b(\opn{mod}_{\Gl})$, this amounts to resolving $\Gl$ as a projective $\Gl$-$\Gl$ module, which can be done using the methods of \cite{King1997TiltingBO} and \cite{BUTLER1999323}. We want to compare the two methods, in fact, the resolutions of ${{\cal O}_{\Gd}}$ for $Y_2,Y_3$ were first found using methods of this section.

The main reason for including this comparison is for the possible deformations of the tangle invariants. 
In \cite{BarmeierWang}, a method was developed for computing the deformation quantisation (to all orders in $\hbar$) of the sheaf category of a variety starting from a certain minimal resolution. We envision applying this method to deform $D^b(Y_n)$ or $D^b(\opn{mod}_{\Gl_n})$, and hence the resulting tangle invariants, in the same spirit of \cite{Def_qnt}. An important feature of this method is that one needs only check the associativity of the deformed quiver algebra on the words contained in the degree $-3$ term of the resolution of $\Gl$ as projective $\Gl$ bi-module. We would therefore potentially obtain 3-parameter (i.e. $t$, $q$, $\hbar$) tangle invariants.

\begin{example}\label{ex_quiver_Y2}
Next we compare the bi-module resolution of ${\cal O}_{\Gd}$ to King's method (\cite{King1997TiltingBO} sec.5) of obtaining a resolution of ${\cal O}_{\Gd}$ by finding a bi-module resolution of the algebra
\bea \Gl=\opn{End}(T),~~~T=L_{00}\oplus L_{10}\oplus L_{21}\oplus L_{31}.\nn\eea
First, the algebra $\Gl$ is a quiver algebra, with idempotents $e_i,i=0,1,2,3$ corresponding to the four line bundles. 
The arrows between nodes are the irreducible maps in between the various line bundles
 \bea \begin{tikzpicture}
  \matrix (m) [matrix of math nodes, row sep=1.2em, column sep=3.5em]
    {   L_{3,1} & L_{2,1} & L_{1,0} & L_{0,0} \\ };
  \path[<-,font =\scriptsize, thick]
  (m-1-1) edge node[above] {$c_i$} (m-1-2)
  (m-1-2) edge node[above] {$b_i$} (m-1-3)
  (m-1-3) edge node[above] {$a_i$} (m-1-4);

  \path[<-,font =\scriptsize, bend right = 25]
  (m-1-1) edge node[below] {$u_2$} (m-1-3);

  \path[<-,font =\scriptsize, bend left = 25]
  (m-1-2) edge node[above] {$u_1$} (m-1-4);
\end{tikzpicture}\label{quiver_Y2}\eea
where the maps are $a_1=b_1=c_1=1$, $a_2=b_2=c_2=x_1/x_0$ and $u_1=u_2=x_{01}/x_{00}$, with the relations
\bea t_i:&&c_iu_1-u_2a_i=0,~~i=1,2,\nn\\
r:&&b_2a_1-b_1a_2=0,\nn\\
s:&&c_2b_1-c_1b_2=0.\label{set_of_rel_Y2}\eea
As an illustration, we calculated already in ex.\ref{ex_sec_Y2} that $\hom(L_{0,0},L_{3,1})$ is generated by
\bea (\frac{x_1}{x_0})^{0,1,2,3},~~~\frac{x_{01}}{x_{00}}(\frac{x_1}{x_0})^{0,1}.\nn\eea
These six maps are represented by the paths $c_1b_1a_1$, $c_1b_1a_2$, $c_1b_2a_2$, $c_2b_2a_2$, $c_1d_1$ and $c_2d_1$ (read the paths from right to left).

One can now write down a resolution of $P_{\sbullet}\to \Gl$ as a $\Gl$ bi-module.
For $\Gl$ of gl.dim 2, there are only $P_{2,1,0}$. The $P_1$ term is generated by length 1 paths, while $P_2$ is generated by the relations. Thus $P_{\sbullet}$ reads
\bea \begin{tikzpicture}
  \matrix (m) [matrix of math nodes, row sep=0.1em, column sep=2em]
    { & \Gl e_2\otimes u_1\otimes e_0\Gl &  \\
       & & \Gl e_0\otimes e_0\Gl \\
      \Gl e_2\otimes (b_2a_1-b_1a_2) \otimes e_0\Gl &\Gl e_1\otimes \vec a\otimes e_0\Gl &  \\
      & &\Gl e_1\otimes e_1\Gl \\
      \Gl e_3\otimes (\vec cu_1-u_2\vec a) \otimes e_0\Gl &\Gl e_2\otimes \vec b \otimes e_1\Gl &  \\
      & &\Gl e_2\otimes e_2\Gl \\
      \Gl e_3\otimes (c_2b_1-c_1b_2) \otimes e_1\Gl &\Gl e_3\otimes  \vec c \otimes e_2\Gl &  \\
      & &\Gl e_3\otimes e_3\Gl \\
      & \Gl e_3\otimes u_2\otimes e_1\Gl & \\};
  \path[->,font =\scriptsize]
  (m-3-2) edge (m-2-3)
  (m-3-2) edge (m-4-3)
  (m-5-2) edge (m-4-3)
  (m-5-2) edge (m-6-3)
  (m-7-2) edge (m-6-3)
  (m-7-2) edge (m-8-3)
  (m-1-2) edge (m-2-3)
  (m-1-2) edge (m-6-3)
  (m-9-2) edge (m-8-3)
  (m-9-2) edge (m-4-3)
  (m-3-1) edge (m-3-2)
  (m-3-1) edge (m-5-2)
  (m-7-1) edge (m-7-2)
  (m-7-1) edge (m-5-2)
  (m-5-1) edge (m-1-2)
  (m-5-1) edge (m-3-2)
  (m-5-1) edge (m-7-2)
  (m-5-1) edge (m-9-2);
\end{tikzpicture}.\label{bi_mod_res_Y2}\eea
It is best if we illustrated the differentials with examples:
\bea 
&\partial(\Gl e_2\otimes u_1\otimes e_0\Gl)=\Gl u_1e_0\otimes e_0\Gl-\Gl e_2\otimes e_2u_1\Gl,\nn\\
&\partial(\Gl e_3\otimes(c_2b_1-c_1b_2)\otimes e_1\Gl)
=\Gl c_2e_2\otimes b_1\otimes e_1\Gl+\Gl e_3\otimes c_2\otimes e_2b_1\Gl-
\Gl c_1e_2\otimes b_2\otimes e_1\Gl-\Gl e_3\otimes c_1\otimes e_2b_2\Gl.\nn\eea
From this resolution, by replacing $\Gl e_i$ with $F_i^{\vee}$ and $e_i\Gl$ with $F_i$ (see def.\ref{def_exp_collection}), we get again the resolution of ${\cal O}_{\Gd}\in D^b(Y_2\times Y_2)$
\bea \begin{tikzpicture}
  \matrix (m) [matrix of math nodes, row sep=0.1em, column sep=2em]
    { & L_{-2,-1}\otimes \frac{s}{t}\otimes L_{00} &  \\
       & & L_{00}\otimes L_{00} \\
      L_{-2,-1}\otimes s \otimes L_{00} &L_{-10}\otimes V^*\otimes L_{00} &  \\
      & &L_{-10}\otimes L_{10} \\
      L_{-3,-1}\otimes \frac{s}{t} V^* \otimes L_{00} &L_{-2,-1}\otimes V^* \otimes L_{10} &  \\
      & & L_{-2,-1}\otimes L_{21} \\
      L_{-3,-1}\otimes s \otimes L_{10} &L_{-3,-1}\otimes V^* \otimes L_{21} &  \\
      & &L_{-3,-1}\otimes L_{31} \\
      &L_{-3,-1}\otimes \frac{s}{t} \otimes L_{10} & \\};
  \path[->,font =\scriptsize]
  (m-3-2) edge (m-2-3)
  (m-3-2) edge (m-4-3)
  (m-5-2) edge (m-4-3)
  (m-5-2) edge (m-6-3)
  (m-7-2) edge (m-6-3)
  (m-7-2) edge (m-8-3)
  (m-1-2) edge (m-2-3)
  (m-1-2) edge (m-6-3)
  (m-9-2) edge (m-8-3)
  (m-9-2) edge (m-4-3)
  (m-3-1) edge (m-3-2)
  (m-3-1) edge[red] (m-5-2)
  (m-7-1) edge (m-7-2)
  (m-7-1) edge (m-5-2)
  (m-5-1) edge (m-1-2)
  (m-5-1) edge (m-3-2)
  (m-5-1) edge (m-7-2)
  (m-5-1) edge (m-9-2);
\end{tikzpicture}\nn\eea 
In terms of $SU(2)$ representations, the maps should be self-evident, e.g. the one marked in red goes as: $s\BB{C}\hookrightarrow sV\otimes V^*\cong V^*\otimes V^*$ and then multiply $V^*$ to the right.
Writing this complex as a matrix, we get exactly the matrix in \eqref{res_Y2}. But it is clear that the matrix form is far more manageable.
\end{example}

\smallskip

Next we record the quiver algebra for $Y_3$:
\bea \begin{tikzpicture}
  \matrix (m) [matrix of math nodes, row sep=6em, column sep=3em]
    { L_{731} & L_{631} & L_{521} & L_{421} & L_{310} & L_{210} & L_{100} & L_{000} \\ };
  \path[<-,line width = 1pt, font =\scriptsize]
  (m-1-1) edge node[above] {$\tilde c_i$} (m-1-2)
  (m-1-2) edge node[above] {$\tilde b_i$} (m-1-3)
  (m-1-3) edge node[above] {$\tilde a_i$} (m-1-4)
  (m-1-4) edge node[above] {$d_i$} (m-1-5)
  (m-1-5) edge node[above] {$c_i$} (m-1-6)
  (m-1-6) edge node[above] {$b_i$} (m-1-7)
  (m-1-7) edge node[above] {$a_i$} (m-1-8);

  \path[<-,font =\scriptsize, bend right = 25]
  (m-1-1) edge node[below] {$\tilde u_2$} (m-1-3)
  (m-1-3) edge node[below] {$\tilde v$} (m-1-5)
  (m-1-5) edge node[below] {$u_2$} (m-1-7);
  \path[<-,font =\scriptsize, bend left = 25]
  (m-1-2) edge node[above] {$\tilde u_1$} (m-1-4)
  (m-1-4) edge node[above] {$v$} (m-1-6)
  (m-1-6) edge node[above] {$u_1$} (m-1-8);

  \path[<-,font =\scriptsize, bend right = 35, line width =1pt]
  (m-1-2) edge node[below] {$\tilde z_i$} (m-1-5);
  \path[<-,font =\scriptsize, bend left = 35, line width =1pt]
  (m-1-4) edge node[above] {$z_i$} (m-1-7);
\end{tikzpicture}\label{quiver_Y3}\eea
Here each length 1 path is an irreducible map between the nodes. For example $\hom(L_{000},L_{210})$ has sections of weights $1,s,s^2,s/t$. The first three are sections of $L_{200}$, so they represent reducible paths: $L_{000}\to L_{100}\to L_{210}$. But the fourth section does not factor through $L_{100}$ and give the irreducible path $u_1$. Similarly $\hom(L_{100},L_{421})=\bra 1,s,s^2,s^3,s/t,s^2t,s/t^2,s^2/t^2\ket$. The first six are paths that factor through intermediate nodes, and the last two give the arrows named $z_1,z_2$. 
The irreducible paths correspond to the monomials at degree $-1$ of the matrix $\Gd_3$ above.

The path algebra obeys a total of 26 relations, divided into various types
\begin{enumerate}
\item Six involve two successive short thick arrows, e.g. $R_{24}=d_2c_1-d_1c_2=0$, where the indices indicate the start and end node. They all have weight $s$.

\item Ten involve two thin curved arrows and two short thick arrows:
\bea \begin{tikzpicture}
  \matrix (m) [matrix of math nodes, row sep=6em, column sep=2em]
    { L_{310} & L_{210} & L_{100} & L_{000} \\ };
  \path[<-,line width = 1pt, font =\scriptsize]
  (m-1-1) edge node[above] {$c_i$} (m-1-2)
  (m-1-3) edge node[above] {$a_i$} (m-1-4);

  \path[<-,font =\scriptsize, bend left = 25]
  (m-1-2) edge node[above] {$u_1$} (m-1-4);
  \path[<-,font =\scriptsize, bend right = 25]
  (m-1-1) edge node[below] {$u_2$} (m-1-3);
  
  \node at (5,0) {$R_{03}^i=u_2a_i-c_iu_1=0$};
\end{tikzpicture}\nn\eea
 $R_{03}^i=u_2a_i-c_iu_1=0$, $R_{14}^i=vb_i-d_iu_2$, $R_{25}^i=\tilde vc_i-\tilde a_iv$, $R_{36}^i=\tilde u_1d_i-\tilde b_i\tilde v$ and $R_{47}^i=\tilde u_2\tilde a_i-\tilde c_i\tilde u_1$. They have weights $s/t$ and $s^2/t$ for $i=1,2$.
\item Four that involve one thick curved arrow
\bea \begin{tikzpicture}
  \matrix (m) [matrix of math nodes, row sep=6em, column sep=2em]
    { \phantom{X} & \phantom{X} & L_{521} & L_{421} & L_{310} & L_{210} & L_{100} & \phantom{X} \\ };
  \path[<-,line width = 1pt, font =\scriptsize]
  (m-1-3) edge node[above] {$\tilde a_i$} (m-1-4);
  
  \path[<-,font =\scriptsize, bend right = 25]
  (m-1-3) edge node[below] {$\tilde v$} (m-1-5)
  (m-1-5) edge node[below] {$u_2$} (m-1-7);
  
  \path[<-,font =\scriptsize, bend left = 35, line width =1pt]
  (m-1-4) edge node[above] {$z_i$} (m-1-7);
  
  \node at (7,0){$R_{15}=\tilde vu_2+\tilde a_1z_2-\tilde a_2z_1$.};
\end{tikzpicture}\nn\eea
There are four such, with the expression
\bea R_{04}=vu_1+z_2a_1-z_1a_2,~R_{15}=\tilde vu_2+\tilde a_1z_2-\tilde a_2z_1,~R_{26}=\tilde u_1v+\tilde z_2c_1-\tilde z_1c_2,~R_{37}=\tilde u_2\tilde v+\tilde c_1\tilde z_2-\tilde c_2\tilde z_1.\nn\eea
They have weights $s^2/t^2$.

\item Two that involve two thick curved arrows
\bea \begin{tikzpicture}
  \matrix (m) [matrix of math nodes, row sep=6em, column sep=1.5em]
    { \phantom{L} & L_{631} & L_{521} & L_{421} & L_{310} & L_{210} & L_{100} & \phantom{L} \\ };
  \path[<-,font =\scriptsize, bend right = 25]
  (m-1-5) edge node[below] {$u_2$} (m-1-7);
  \path[<-,font =\scriptsize, bend left = 25]
  (m-1-2) edge node[above] {$\tilde u_1$} (m-1-4);
  
  \path[<-,font =\scriptsize, bend right = 35, line width =1pt]
  (m-1-2) edge node[below] {$\tilde z_i$} (m-1-5);
  \path[<-,font =\scriptsize, bend left = 35, line width =1pt]
  (m-1-4) edge node[above] {$z_i$} (m-1-7);
\end{tikzpicture}\nn\eea
which read $R_{16}^i=\tilde z_iu_2-\tilde u_1z_i$, $i=1,2$ and have weights $s^2/t^3,s^3/t^3$ for $i=1,2$.
\item Final four that involve three paths, e.g.
\bea \begin{tikzpicture}
  \matrix (m) [matrix of math nodes, row sep=6em, column sep=1.5em]
    { \phantom{L} & L_{6,3,1} & L_{5,2,1} & L_{4,2,1} & L_{3,1,0} & L_{2,1,0} & L_{1,0,0} & \phantom{L} \\ };
  \path[<-,line width = 1pt, font =\scriptsize]
  (m-1-2) edge node[above] {$\tilde b_i$} (m-1-3)
  (m-1-3) edge node[above] {$\tilde a_i$} (m-1-4)
  (m-1-5) edge node[above] {$c_i$} (m-1-6)
  (m-1-6) edge node[above] {$b_i$} (m-1-7);
  \path[<-,font =\scriptsize, bend right = 35, line width =1pt]
  (m-1-2) edge node[below] {$\tilde z_i$} (m-1-5);
  \path[<-,font =\scriptsize, bend left = 35, line width =1pt]
  (m-1-4) edge node[above] {$z_i$} (m-1-7);
\end{tikzpicture}\nn\eea
The four relations read
\bea R_{16}^3=\tilde z_1c_1b_1-\tilde b_1\tilde a_1z_1,~R_{16}^4=\tilde z_1c_2b_1-\tilde b_2\tilde a_1z_1,~ R_{16}^5=\tilde z_2c_2b_1-\tilde b_1\tilde a_2z_2,~R_{16}^6=\tilde z_2c_2b_2-\tilde b_2\tilde a_2z_2.\nn\eea
They have weights $s^{1,2,3,4}/t^2$ respectively.
\end{enumerate}
The 26 relations correspond to the monomials at degree $-2$ of the matrix $\Gd_3$ above.
It is gratifying to see that our matrix method above reproduces all the relations, which is no easy matter to spot from the quiver.

The matrix entries in $\Gd_3$ of degree $-3,-4$ have interpretations as relations of relations, etc. We will not list these higher relations, but are happy to provide them if requested.

\end{document}